\documentclass{amsart}
\usepackage{amssymb}
\usepackage{amscd}
\usepackage{graphicx}
\usepackage[cmtip,all]{xy}

\copyrightinfo{2007}{American Mathematical Society}

\newtheorem{theorem}{Theorem}[section]
\newtheorem{lemma}[theorem]{Lemma}
\newtheorem{proposition}[theorem]{Proposition}
\newtheorem{corollary}[theorem]{Corollary}
\newtheorem{cordef}[theorem]{Corollary/Definition}

\theoremstyle{definition}
\newtheorem{definition}[theorem]{Definition}
\newtheorem{notation}[theorem]{Notation}
\newtheorem{example}[theorem]{Example}

\theoremstyle{remark}
\newtheorem{remark}[theorem]{Remark}

\numberwithin{equation}{section}

\newcommand{\bA}{\mathbb{A}}
\newcommand{\bP}{\mathbb{P}}
\newcommand{\bR}{\mathbb{R}}

\newcommand{\bQ}{\mathbb{Q}}
\newcommand{\bZ}{\mathbb{Z}}
\newcommand{\bC}{\mathbb{C}}

\newcommand{\sH}{\mathrm{H}}
\newcommand{\calC}{\mathcal{C}}

\newcommand{\calL}{\mathcal{L}}
\newcommand{\calM}{\mathcal{M}}
\newcommand{\calO}{\mathcal{O}}
\newcommand{\calP}{\mathcal{P}}

\newcommand{\calT}{\mathcal{T}}

\newcommand{\calS}{\mathcal{S}}
\newcommand{\calG}{\mathcal{G}}
\newcommand{\calH}{\mathcal{H}}

\newcommand{\calY}{\mathcal{Y}}
\newcommand{\calD}{\mathcal{D}}

\newcommand{\Ext}{\mathrm{Ext}}

\newcommand{\Sym}{\mathrm{Sym}}

\newcommand{\SL}{\mathrm{SL}}

\newcommand{\Gm}{\mathbb{G}_m}
\newcommand{\Proj}{\mathrm{Proj}}
\newcommand{\Sing}{\mathrm{Sing}}

\newcommand{\lc}{\mathrm{lct}}

\newcommand{\Pic}{\mathrm{Pic}}

\newcommand{\rank}{\mathrm{rank}}
\newcommand{\scl}{\widetilde{S}_{(C,L)}}
\newcommand{\gquot}{/\!\!/}
 
\newcommand{\Supp}{\mathrm{Supp}}

\DeclareMathOperator{\Aut}{Aut}

\DeclareMathOperator{\Sat}{Sat}
\DeclareMathOperator{\Spec}{Spec}

\DeclareMathOperator{\mult}{mult}

\DeclareMathOperator{\rk}{rk}
\DeclareMathOperator{\Div}{div}

\newcommand{\Es}{\widetilde{E_7}}
\newcommand{\Ee}{\widetilde{E_8}}

\begin{document}
\title{Deformations of singularities and variation of GIT quotients}
\author{Radu Laza}
\address{Department of Mathematics, University of Michigan, Ann Arbor, MI 48109}
\email{rlaza@umich.edu}

\subjclass[2000]{Primary 14J17, 14B07, 32S25; Secondary 14L24}
\date{}
\bibliographystyle{amsplain}

\begin{abstract}
We study the deformations of the minimally elliptic surface singularity $N_{16}$. A standard argument reduces the study of the deformations of $N_{16}$ to the study of the moduli space of pairs $(C,L)$ consisting of a plane quintic curve and a line. We construct this moduli space in two ways: via the periods of $K3$ surfaces and by using geometric invariant theory (GIT).  The GIT construction depends on the choice of the linearization. In particular, for one choice of linearization we recover the space constructed via $K3$ surfaces and for another we obtain the full deformation space of $N_{16}$. The two spaces are related by a series of explicit flips. In conclusion, by using the flexibility given by  GIT and the standard tools of Hodge theory,  we obtain a good understanding of the deformations of $N_{16}$.
\end{abstract}

\maketitle

\section{Introduction}\label{sectintro} 
Singularities and their deformations have always played a central role in algebraic geometry, being of  fundamental importance in several branches of the field, such as the classification of surfaces, the minimal model program, and the compactification problem for moduli spaces.  The easiest and the first to be understood were the {\it simple singularities} (also known as duVal singularities).  Work of many mathematicians, including Brieskorn, Pinkham, and Looijenga, have extended the results for simple singularities to the next level of complexity, the {\it unimodal singularities} (simple elliptic, cusp and triangle). The focus of this paper is a detailed study of a new class of singularities: the minimal-elliptic surface singularity $N_{16}$. The class $N_{16}$ sits immediately after the simple and unimodal singularities in Arnold's hierarchy of singularities, and its understanding is essential to any attempt of studying deformations of singularities more complex than the unimodal ones. 

\smallskip

The most effective tool for the study of the deformations of the unimodal singularities is the theory of deformations with $\bC^*$-action  of Pinkham \cite{pinkham0}. The starting point is that most of the unimodal singularities have a good $\bC^*$-action such that the induced action on the tangent space to the deformations has  all but one of the weights negative. The non-negative direction is topologically trivial and can  be ignored. On the other hand, Pinkham's theory says that the deformations in the negative direction can be globalized and interpreted as a moduli space of certain pairs (see \cite[Appendix]{looijengatriangle}). Thus, the deformation problem is essentially reduced to  a moduli problem for which  standard algebro-geometric tools are available. For example, in the case of the triangle singularities, the pairs $(S,H)$  under consideration consist of  a $K3$ surface $S$ and a divisor $H$  such that $H$ forms a fixed configuration of rational curves. Since the moduli space of (lattice  polarized) $K3$ surfaces is well understood, one obtains a wealth of information about the deformations of the triangle singularities (e.g. \cite{pinkhammonodromy, pinkhamduality, novaacta, looijengatriangle}). 

\smallskip

The case of the singularity $N_{16}$ is similar,  but there is one important difference: due to the increase in modality, the zero weight directions are non-trivial and cannot be ignored. A partial modular interpretation exists for the zero weight direction as well, but in contrast to the pure negative weight situation, the globalization is no longer guaranteed. As explained below, we solve this globalization problem by using GIT. Once this is done, we use the resulting moduli space to get a good hold on the structure of the deformation space of $N_{16}$. A short description of the content of the paper is given below. Further details and statements of the main results are contained in the introductions of the individual sections.

\smallskip

The singularity $N_{16}$ is a double cover of the cone over $5$ points in $\bP^1$. Consequently, following the method of Pinkham, we can essentially identify the deformations of non-positive weight of $N_{16}$ with the moduli space of pairs $(C,L)$ consisting of a plane quintic curve and a line such that the intersection is transversal.   A natural approach is to construct this moduli space as a geometric invariant theory (GIT)  quotient.  We start by studying the moduli space of pairs $(C,L)$, where  $C$ is a plane curve of degree $d$  and $L$ is a line. The moduli of such pairs is then $X\gquot G$, where $X\cong |dL|\times |L|$ is the parameter space for pairs and $G=\SL(3)$ acts naturally via the diagonal action. By definition, the construction depends on the choice of a linearization $\calL\in \Pic^G(X)$, which is parameterized by a single rational parameter $t\in \bQ_+$,   the slope of $\calL$. For each such choice, we obtain a moduli space of pairs $\mathcal{M}(t)$. This type of situation was analyzed in general circumstances by Thaddeus \cite{thaddeus} and Dolgachev-Hu \cite{dolgachevhu}. In particular, it is known that there exists only a finite number of non-isomorphic  quotients $\calM(t)$ related by explicit birational transformations.

\smallskip

In section \ref{sectgen}, we establish a number of general  qualitative results (valid for any degree $d$) on  the dependence of   the GIT stability for degree $d$ pairs $(C,L)$  on the parameter $t$. Namely, there are two main results here. The first result, the interpolation theorem (theorem \ref{mainthm1}), says that the stability  at $t=0$, $1$ and $\frac{d}{2}$ is equivalent to the stability  of $C$ as a degree $d$ plane curve, of $C+L$ as a degree $d+1$  curve,  and  of the intersection $C\cap L$ as a $d$-tuple of points in $\bP^1$ respectively. The second result (theorem \ref{mainthm2}) relates the stability of the pair $(C,L)$ for the slope $t$ with the singularities of the divisor pair $(\bP^2,\frac{3}{d+t}(C+tL))$. Namely, if the pair  $(\bP^2,\frac{3}{d+t}(C+tL))$ is log canonical, then  $(C,L)$ is semistable for the slope $t$. We note that this is a relative version of earlier results of Hacking \cite[\S10]{hacking} and Kim-Lee \cite{kimlee}. As a consequence of these two results, the dependence of the stability condition on the parameter $t$ can be roughly stated as saying  that as  $t$ increases from $0$ to $\frac{d}{2}$ the curve $C$ is allowed to have more complicated singularities, but we require stronger transversality conditions on the intersection $C\cap L$ (see Ex. \ref{exampled3}).

\smallskip

 Also in section \ref{sectgen}, we note that the variation of GIT for the pairs $(C,L)$ is closely related to the deformations of the cones over $d$-tuples of points in $\bP^1$ (see  \S\ref{sectpinkham}). The basic observation is that, due to the interpolation theorem, the minimal closed orbits at $t=\frac{d}{2}$ are the pairs $(C,L)$ with $C$  a cone and $L$  a line not passing through the vertex of $C$. Furthermore, if $C$ is not a cone and $L$ is transversal, the pair $(C,L)$ is stable at $t=\frac{d}{2}-\epsilon$ for $\epsilon$ small. It follows then that the variation of GIT morphism $\calM(\frac{d}{2}-\epsilon)\to \calM(\frac{d}{2})$ is a global object associated to the natural retraction map $S_{\le 0}\to S_0$ modulo the $\bC^*$-action, where $S_{\le 0}$ and $S_0$ denote the deformations of the cone of non-positive and zero weight respectively. The fibers of $\calM(\frac{d}{2}-\epsilon)\to \calM(\frac{d}{2})$ are (at least generically) weighted projective spaces corresponding to the negative weight deformations modulo the $\bC^*$-action. In other words, the variation of GIT quotients as $t$ increases from $\frac{d}{2}-\epsilon$ to $\frac{d}{2}$ is  the standard globalization of Pinkham in a relative version over the zero weight deformations (i.e. the equisingular stratum).  

\smallskip

In section \ref{sectdeg5}, we do a detailed analysis of the stability condition for the degree $5$ case. The general results of section \ref{sectgen} specialize nicely in this situation, greatly simplifying the analysis.  We restrict here to the discussion of the special role played by the slope $t=1$. First, from theorem \ref{mainthm1}, it follows that a pair $(C,L)$ is (semi)stable at $t=1$ if and only if $C+L$ is (semi)stable as a plane sextic. Thus, a pair $(C,L)$ such that  $C+L$ is reduced and has at worst simple singularities is stable for slope $t=1$ (cf. Shah \cite{shah}). We then note that the variation of GIT quotients at $t=1$ corresponds to a natural division of singularities in three large classes. Namely,  assume for simplicity that $L$ is transversal to $C$, then the pair $(C,L)$ is semistable for some $t<1$ if and only if $C$ has at worst simple singularities. For $t=1$, $C$ is allowed to have simple elliptic singularities. Finally, the pairs with $C$ having worst singularities become semistable only for some $t>1$.  By taking the double cover of $\bP^2$ branched along $C+L$, this division of singularities in three types with respect to the stability condition is conceptually explained by theorem \ref{mainthm2} and the well-known division of the surface singularities: canonical (the rational double points), log canonical (the simple elliptic and cusp singularities), and  not log canonical. For us, this division is relevant due to a theorem of Shah \cite{shahinsignificant} which says that the log canonical surface singularities are insignificant limit singularities. As a consequence, we can relate the GIT construction of sections \ref{sectgen} and \ref{sectdeg5} to the Hodge theoretical construction of section \ref{sectk3}.

\smallskip

 In section \ref{sectk3}, we note that there exists a simple alternative construction of the moduli space of degree $5$ pairs. Namely, we view a degree $5$ pair $(C,L)$ as a plane sextic $B=C+L$, and  associate to it the double cover $S_{(C,L)}$ of $\bP^2$ branched along $B$. Assuming that $(C,L)$ is a GIT stable pair at $t=1$, it follows that $S_{(C,L)}$ has only rational double points. Thus, its desingularization $\widetilde{S}_{(C,L)}$ is a degree two $K3$ surface. The special nature of the sextic $B$ imposes conditions on the Neron-Severi lattice of $\widetilde{S}_{(C,L)}$. It follows that $\widetilde{S}_{(C,L)}$ is an $M$-polarized $K3$ surface (see \cite{mirrork3}) for a certain rank $6$ hyperbolic lattice $M$. The moduli space of $M$-polarized $K3$ surfaces is well known to be  locally symmetric of type $\calD/\Gamma$ for appropriate choices of a type IV domain $\calD$ and  of an arithmetic group $\Gamma$ acting on $\calD$. The main result of  section \ref{sectk3}, theorem \ref{mainthm4},  says  then that this construction gives an isomorphism   $\calM(1)\cong (\calD/\Gamma)^*$ between the GIT quotient $\calM(1)$ and  the Baily-Borel compactification  $(\calD/\Gamma)^*$  of $\calD/\Gamma$. In conclusion, we obtain a dual description for $\calM(1)$. As explained below, this fact has numerous consequences on the structure of the deformations of $N_{16}$.
 
 \smallskip
 
 In the last section, we discuss the implications of the results described above on the structure of  deformations of the singularity $N_{16}$. Similarly to the situation of unimodal singularities, we analyze the structure of the discriminant hypersurface in the versal deformation of $N_{16}$ and the possible combinations of  singularities occurring in a nearby fiber. Our main conclusion is that essentially all the results of Pinkham, Looijenga, and Brieskorn for unimodal singularities have a natural counterpart in our situation (see theorems \ref{thmstructure}, \ref{simplecombinations}, and \ref{thmstrata}). 
 
 \smallskip
 
 There are two main ideas involved in the study of the deformations of $N_{16}$. First, from the results of section \ref{sectgen}, the correct global object associated to the deformations of non-positive weight of $N_{16}$ is the fibration $\calM(\frac{5}{2}-\epsilon)\to\calM(\frac{5}{2})$. However, by itself, this description does not give much. Thus, the second main idea is  to gain information about the deformation space by exploiting the description of the moduli space of degree $5$ pairs as the quotient $\calD/\Gamma$. The basic observation that makes this useful  is that Hodge theory transforms many questions about singularities into purely arithmetic statements. There exists, however, a  disadvantage to this approach. Namely,  the Hodge theoretical arguments work well only for the singularities of finite monodromy, i.e. the simple singularities. This means that  $\calD/\Gamma$ gives a good description only for the simple singularity stratum in the deformation space of $N_{16}$. From our point of view, this is completely natural and it is easily rectified. We recall that $\calD/\Gamma$ and the deformation space of $N_{16}$ correspond to the GIT quotients $\calM(1)$  and $\calM(\frac{5}{2}-\epsilon)$ respectively. Thus, the  information that is missing from the $\calD/\Gamma$ description can be recovered by following the series of  explicit birational modifications that relate  $\calM(1)$  and $\calM(\frac{5}{2}-\epsilon)$. Geometrically, the variation of GIT quotients $\calM(1)\dashrightarrow \calM(\frac{5}{2}-\epsilon)$  can be interpreted as introducing, one at a time, the non-simple singularity strata in the deformation of $N_{16}$. 
 
 \smallskip
 
One interesting aspect about the dual construction (GIT/Hodge theory) of this paper is that the flips that transform $\calM(1)$ into $\calM(\frac{5}{2}-\epsilon)$ can be interpreted also in terms of arithmetic arrangements of hyperplanes (N.B. a priori they have only a GIT meaning).  As explained in \S\ref{sectriangle}, this is closely related to  Looijenga's construction \cite{looijengacompact}. Essentially, the GIT approach is dual to the approach of \cite[\S10]{looijengacompact}. For the deformations of the triangle singularities, the two approaches coincide. For $N_{16}$ the situation is less clear, but conjecturally we should again have a coincidence. 

\smallskip

We close by noting that the techniques developed in this paper can be applied to other classes of singularities as well. Specifically, we have in mind the threefold singularity $O_{16}$, the cone over a cubic surface. This case would be the first example of a detailed study of the deformation space for a genuine threefold singularity  (not a suspension of a surface singularity). The two key ingredients of our method -- the flexibility given by GIT and the explicit model obtained via the period map -- have natural counterparts for $O_{16}$. Namely, the GIT analysis adapts well in higher dimensions, and the construction of section \ref{sectk3} can be done by using  cubic fourfolds instead of  K3 surfaces.  Details will appear elsewhere.

\subsubsection*{Acknowledgment} This paper is a revised version of the author's thesis.  I am grateful to my advisor, Robert Friedman, for his guidance and help.  I would also like to thank Michael Thaddeus for teaching me about the variation of GIT quotients, and Igor Dolgachev and Rob Lazarsfeld for helpful comments and suggestions. The referee's comments  helped me better organize and clarify the material.

\subsection*{Notations and Conventions}
For the  basic GIT notions and notations, we follow Mumford \cite{GIT}. 
The conventions and notations for singularities are those of Arnold et al. \cite{agv1}. The only notable difference is the use of $\widetilde{E_r}$ for $r=6,7,8$ to denote the simple elliptic (parabolic) singularities.  In addition to the simple ($A_n$, $D_n$, and $E_r$), simple elliptic ($\widetilde{E_r}$), and cusp singularities ($T_{p,q,r}$), we are concerned with the triangle (exceptional unimodal) singularities: $Z_{11}$, $Z_{12}$, $W_{12}$, and $W_{13}$ (\cite[pg. 247]{agv1}). 
The singularity $N_{16}$ is a trimodal singularity with  normal form:
$$N_{16} : f_5(x, y) + z^2,$$
where $f_5$ is a homogeneous degree $5$ polynomial. The singularity $N_{16}$ deforms only to the simple and unimodal singularities listed above. 

Since the deformation spaces of a singularity and of its suspensions can be identified, all singularities are considered up to stable equivalence (i.e. up to adding squares of new variables). In particular, depending on the context,  $N_{16}$ refers to either a surface singularity or a  curve singularity.

\section{Variation of GIT quotients for pairs $(C,L)$}\label{sectgen}
In this section, we construct the moduli space of pairs $(C,L)$ consisting of a plane curve of degree $d$ and a line by using the geometric invariant theory (GIT).  We then  relate this construction to the deformations of ordinary multiplicity $d$ points. 

\begin{definition}
We call a pair $(C,L)$ consisting of a plane curve $C$ of degree $d$ and a line $L\subset \bP^2$ a {\it degree $d$ pair}. Two such pairs are {\it equivalent} if they are projectively equivalent.
\end{definition}

The natural GIT set-up for the study of the moduli of pairs is that of the group $G=\textrm{SL}(3)$ acting diagonally on the parameter space $X$  of degree $d$ pairs, where $$X=\bP(\sH^0(\bP^2,\calO_{\bP^2}(d)))\times\bP(\sH^0(\bP^2,\calO_{\bP^2}(1)))\cong\bP^N\times \bP^2$$ 
and $N={d+2 \choose 2}-1$. The GIT quotient $X\gquot G$ depends on the choice of an ample $G$-linearized line bundle $\calL\in\Pic^G(X)$. Namely, by definition, we have
\begin{equation*}
X\gquot_{\calL} G=\Proj \bigoplus_{n\ge 0} \sH^0(X,\calL^{\otimes n})^G
\end{equation*}
The dependence of the quotient on the choice of the linearization was analyzed in general circumstances by Thaddeus \cite{thaddeus} and Dolgachev-Hu \cite{dolgachevhu}. In our situation, since $\Pic^G(X)\cong \Pic(X)\cong \bZ\times \bZ$, the results of \cite{dolgachevhu,thaddeus} say that the quotient $X\gquot_{\calL} G$ depends on a single rational parameter $t\in \bQ_+$, the slope of  $\calL$, and that only finitely many non-isomorphic quotients are actually obtained. 

\begin{definition}
An ample linearization $\mathcal{L}\in \textrm{Pic}^G(X)$ is said {\it to be of slope $t\in \mathbb{Q}_+$} 
if $\mathcal{L}\cong\pi_1^*\calO_{\bP^N}(a)\otimes \pi_2^*\calO_{\bP^2}(b)$ with $t=\frac{b}{a}$. We denote by $X^s(t)$ and
$X^{ss}(t)$ the sets of stable points and semistable points respectively. The corresponding GIT quotient is denoted $\calM(t)$ or $X\gquot_t G$. A point $x\in X^{ss}(t)$ will be called  
{\it $t$-semistable} or {\it semistable at $t$}  (and similarly for stable and unstable points). All these notions depend only on $t$.
\end{definition}

\begin{remark}
The definitions make sense also for the two extremal cases $\calL=\pi_1^*\calO_{\bP^N}(1)$ and $\calL=\pi_2^*\calO_{\bP^2}(1)$. We talk about linearizations of slope $0$ and $\infty$. It is immediate that $\calM(0)$ is isomorphic to the moduli space of degree $d$ plane curves and that $\calM(\infty)=\emptyset$.  
\end{remark}

From the general results of the theory of variation of GIT quotients, it follows that   there exists a finite number of critical slopes, say $0=t_0<t_1<\dots<t_n<\infty$, such that:
\begin{itemize}
\item[i)] $\calM(t)\neq \emptyset$ iff $t\in[t_0,t_n]$;
\item[ii)] $\calM(t)$ is birational to $\calM(t')$ for all $t,t'\in (t_0,t_n)$ (N.B. $\calM(t_0)$ and $\calM(t_n)$ are lower dimensional); 
\item[iii)] $\calM(t)\cong \calM(t')$ for $t,t'\in(t_i,t_{i+1})$ (for $i=0,\dots,n-1$);
\item[iv)] For small $\epsilon>0$, we have the following commutative diagram
\begin{equation*}
\xymatrix@R=.25cm{
{\calM(t_i-\epsilon)}\ar@{>}[ddr]_{f_{-}}\ar@{<-->}^{f}[rr]&   &{\calM(t_i+\epsilon})\ar@{>}[ddl]^{f_{+}}\\
                                 &   &                   \\
                                 &\calM(t_i)&                  
}
\end{equation*}
with $f_-$, $f_+$ birational morphisms, and $f$ a flip (for $i=1,\dots,n-1$). Furthermore, for the extremal critical values, there exist fibration morphisms $\calM(t_0+\epsilon)\to \calM(t_0)$ and $\calM(t_n-\epsilon)\to \calM(t_n)$. 
\end{itemize}
The standard terminology is to call the critical slopes $t_i$ {\it walls}, the intervals $(t_i,t_{i+1})$ {\it chambers}, and the birational modifications from iv) {\it wall crossings}. 

In the context of these general results,  we are interested in the following questions: 
\begin{enumerate}
\item[(Q1)] Given a pair $(C,L)$, for which values of the slope $t\in \mathbb{Q}_+$ is $(C,L)$ a semistable pair?
\item[(Q2)] Find the critical values $t_1, \dots, t_n$ and describe their geometric relevance.
\item[(Q3)] Describe the wall crossing that occurs for slope $t_i$.
\end{enumerate}
In this section, we establish a series of qualitative answers (valid for all degrees) to these questions. The degree $5$ case is  then considered in detail in section \ref{sectdeg5}.

The first general result about the stability of degree $d$ pairs is to identify the boundary walls $t_0=0$ and $t_n=\frac{d}{2}$, and to describe the stability condition for the slopes $0$, $1$ and $\frac{d}{2}$. 

\begin{theorem}\label{mainthm1}
Let $(C,L)$ be a degree $d$ pair. Then, there exists an interval (possibly empty) 
$[\alpha,\beta]\subset[0,\frac{d}{2}]$ such that $(C,L)$ is $t$-semistable if and only if $t\in[\alpha,\beta]$. 
Furthermore,
\begin{itemize}
\item[i)] $\alpha=0$ if and only if $C$ is a semistable degree $d$ plane curve. If $C$ is stable as a plane curve, we also have $\beta>0$.
\item[ii)] $1\in [\alpha,\beta]$ if and only if $C+L$ is semistable as a degree $d+1$ plane curve. If $C+L$ is stable as a plane curve, then $\alpha<1<\beta$.
\item[iii)] $\beta=\frac{d}{2}$ if and only if $C\cap L$ forms a semistable $d$-tuple of points in $L\cong \bP^1$. If $C\cap L$ is stable as a $d$-tuple and $C$ is not a cone, then $\alpha<\frac{d}{2}$.
\end{itemize} 
\end{theorem}

The theorem says that $t$ interpolates between two conditions of stability: the stability
of degree $d$ curves in $\bP^2$ (at $t=0$) and the stability of $d$-tuples of points in $\bP^1$ (at $t=\frac{d}{2}$). This can be rephrased as saying that, as $t$ increases, we allow $C$ to be more singular, but  require stronger transversality conditions for $C\cap L$. More precisely, we have (see also \cite[\S10]{hacking}, \cite{kimlee}):
\begin{theorem}\label{mainthm2}
Let $(C,L)$ be a degree $d$ pair. If the pair $(\bP^2,\frac{3}{d+t} (C+t L))$ is log canonical, then $(C,L)$ is $t$-semistable. 
\end{theorem}

The following example illustrates  theorems \ref{mainthm1} and \ref{mainthm2}. 

\begin{example}\label{exampled3} 
For degree $3$ pairs, the critical slopes are $t_0=0$, $t_1=\frac{3}{5}$, $t_2=1$  and $t_3=\frac{3}{2}$. The stability of a degree $3$ pair $(C,L)$ is described by the following rules:
\begin{enumerate}
\item[i)] If  $L$ passes through a singular point of $C$, then the pair is $t$-unstable for all $t>0$.
\item[ii)] Assume that $L$ does not meet $C$ in a singular point. Then the pair $(C,L)$ is $t$-(semi)stable iff 
 $t\in (\alpha,\beta)$  ($t\in [\alpha,\beta]$ respectively), where $\alpha$ and $\beta$ are given by 
$$
\alpha=\begin{cases}
0           &\textrm{ if }C \textrm{ has at worst }A_1 \textrm{ singularities}\\
\frac{3}{5} &\textrm{ if }C \textrm{ has an }A_2 \textrm{ singularity}\\
1           &\textrm{ if }C \textrm{ has an }A_3 \textrm{ singularity}\\
\frac{3}{2} &\textrm{ if }C \textrm{ has a }D_4 \textrm{ singularity}                  
\end{cases}
\textrm{ and }
\beta=\begin{cases}
\frac{3}{5} &\textrm{ if }L \textrm{ is inflectional to }C\\
1           &\textrm{ if }L \textrm{ is tangent to }C\\
\frac{3}{2} &\textrm{ if }L \textrm{ is transversal to }C                  
\end{cases}
$$
\end{enumerate}
For degree $3$ pairs, the converse of  \ref{mainthm2} also holds (compare with Ex. \ref{nonanalytic}). For example, the case $\alpha=\frac{3}{5}$ is equivalent to saying that the log canonical threshold of a cusp is $\frac{5}{6}$.
\end{example}

The proofs of \ref{mainthm1} and \ref{mainthm2} are given in \S\ref{proofmainthm1} and \S\ref{proofmainthm2} respectively. Additionally, we discuss in this section the algorithmic determination of the critical slopes (\S\ref{sectalgo}), and the relation to the theory of deformations with $\bC^*$-action (\S\ref{sectpinkham}).

\subsection{The numerical criterion for pairs}\label{numcriterion} The main tool of investigating the dependence of the stability condition  on the choice of linearization is  the Hilbert-Mumford numerical criterion (\cite[Thm. 2.1]{GIT}): {\it a point $x\in X$ is stable (semistable) with respect to a linearization $\mathcal{L}\in \Pic^G(X)$  if and only if  $\mu^\mathcal{L}(x,\lambda)>0$ 
(resp. $\mu^\mathcal{L}(x,\lambda)\ge0$) for every nontrivial $1$-PS $\lambda$ of $G$}, where $\mu^\mathcal{L}(x,\lambda)$ is the numerical function of Mumford (\cite[Def. 2.2]{GIT}). 

The functorial properties of $\mu^\mathcal{L}(x,\lambda)$  give the following identity:
\begin{equation*}
\mu^{\mathcal{O}(a,b)}(x,\lambda)=a \mu^{\mathcal{O}_{\bP^2}(1)}(c,\lambda)+ b \mu^{\mathcal{O}_{\bP^2}(1)}(l,\lambda)=a(\mu(c,\lambda)+\frac{b}{a}
 \mu(l,\lambda))
 \end{equation*}
where $\mathcal{O}(a,b)=\pi_1^*\mathcal{O}_{\bP^N}(a)\otimes\pi_2^*\mathcal{O}_{\bP^2}(b)$, $x=(c,l)\in X\cong \bP^N\times \bP^2$, and $\mu(c,\lambda):=\mu^{\mathcal{O}(1)}(c,\lambda)$ is the standard numerical function used for degree $d$ plane curves. Since the numerical criterion tests only for the sign of $\mu^{\mathcal{O}(a,b)}(x,\lambda)$, we can normalize by dividing by $a$. Thus, we test the (semi)stability of a point $x$ by using the function:
\begin{equation*}
\mu^t(x,\lambda):=\mu(c,\lambda)+t \mu(l,\lambda)
\end{equation*} 
where $t=\frac{b}{a}\in \mathbb{Q}_+$ is the slope of the linearization $\calL=\mathcal{O}(a,b)$. In particular, since $\mu^t(x,\lambda)$ is linear in $t$, we obtain the following corollary of the numerical criterion: 

 \begin{cordef}\label{intstab}
 For every $x\in X$ there exists a finite (possibly empty) interval  $[\alpha,\beta]$ ($\alpha,\beta\in\bQ$) such that
 \begin{enumerate}
 \item[i)] $x$ is semistable for $t$ $\iff$ $t\in  [\alpha,\beta]\cap (\bQ_+\cup\{0\})$.
 \item[ii)] if $x$ is stable for some $t$, then it is stable for all  $t\in  (\alpha,\beta)\cap \bQ_+$.
 \end{enumerate}
 We will call the interval $[\alpha,\beta]$  {\bf the interval of stability} of the point $x$.
 \end{cordef}
 \begin{proof}
 Fix $x\in X$ and associate to a $1$-PS $\lambda$ the closed interval  $I_{\lambda}:=\{t\in [0,\infty) \mid \mu^t(x,\lambda)\ge 0\}$. By the numerical criterion, $x$ is $t$-semistable iff $t\in \cap_{\lambda} I_\lambda=: [\alpha,\beta]$. The rationality of the endpoints and the second statement follow easily from the fact that, for the application of the numerical criterion for $x$, one only needs to consider a finite number of subgroups $\lambda$ (see \ref{remfinite}).
 \end{proof}
 
\subsubsection{The numerical criterion and configurations of monomials}\label{sectconfigurarion}
For applications, it is important to make $\mu^t(x,\lambda)$ explicit. As is customary, we choose coordinates such that the $1$-PS $\lambda$ is 
diagonal, i.e. $\lambda$ is given by
\begin{equation*}
s\in \Gm\xrightarrow{\lambda} \textrm{diag}(s^{r_0},s^{r_1},s^{r_2})\in G= \SL(3)
\end{equation*}
for some weights  $r_0,r_1,r_2\in\bZ$ (not all zero) with  $r_0+r_1+r_2=0$.
We assume additionally that $r_0\ge r_1\ge r_2$. With respect to 
these coordinates, a point $x\in X$ is represented by 
two homogenous forms $c=\sum_{i,j} c_{ij}x_0^{d-i-j}x_1^ix_2^j$ 
and $l=a_0x_0+a_1x_1+a_2x_2$ of degrees $d$ and $1$ respectively. Then,  $\mu^t(x,\lambda)=\mu(c,\lambda)+t \mu(l,\lambda)$ is computed by (\cite[pg. 81]{GIT}):
\begin{equation}\label{temp3}
\mu(c,\lambda)=\max\left\{(d-(i+j)r_0+i r_1+j r_2 \mid \textrm{for all  } i,j \textrm{ such that } c_{ij}\neq 0\right\}
\end{equation}
and
\begin{equation}\label{temp4}
\mu(l,\lambda)=\max\left\{r_i\mid a_{i}\neq 0\right\}.
\end{equation}

We note that the function $\mu^t(x,\lambda)$  depends on the slope $t\in \bQ_{+}$ and two other ingredients: the weights of $\lambda$ and the monomials occurring with non-zero coefficient in $c$ and $l$.  It is convenient to further normalize $\mu^t(x,\lambda)$ as follows:

\begin{definition}\label{defnormalize}
Let $\lambda$ be $1$-PS.  Fix coordinates such that $\lambda$ is diagonal with weights $r_0,r_1,r_2$ as above. For any $x=(c,l)\in X$, denote by $\Xi=(\Xi_d,\Xi_1)$ the set of monomial occurring with non-zero coefficient in $c$ and $l$, and call it {\it the associated configuration of monomials}. We denote $r:=\frac{r_1}{r_0}$ ({\it the normalized weight of $\lambda$})  and  $||\lambda||:=r_0$ ({\it the norm of $\lambda$}).  We then define
\begin{equation*}
\mu^t(\Xi,r):=\frac{\mu^t(x,\lambda)}{||\lambda||}
\end{equation*}
(and similarly  $\mu(\Xi_i,r)$ for $i=1,d$).
\end{definition}

\begin{remark}\label{remfinite} 
We make the following simple observations about the previous definition:
\begin{itemize}
\item[(1)] As suggested by notation,  $\mu^t(\Xi,r)$ depends only on $\Xi$, $r$, and $t$. The function $\mu^t(\Xi,r)$ is linear in $t$, and piecewise linear in $r$ .
\item[(2)] For a fixed degree, there are only finitely many possibilities for $\Xi$. Since $r_0+r_1+r_2=0$ and $r_0\ge r_1\ge r_2$ we get $r\in [-\frac{1}{2},1]$.
\item[(3)] Given a degree $d$ pair $(C,L)$ we say that {\it $\Xi$ is associated to $(C,L)$} if there exists a choice of coordinates such that $\Xi$ is associated to the defining equations $(c,l)$ of the pair. The numerical criterion can be restated as: {\it a pair $(C,L)$ is $t$-semistable iff for any configuration $\Xi$ associated to the pair we have $\min_{r\in[-\frac{1}{2},1]}\mu^t(\Xi,r)\ge 0$.}
\item[(4)] In particular, the finiteness results of the theory of variation of GIT quotients are easily obtained in our situation. For example, it follows that for a given pair $(C,L)$, one needs to apply the numerical criterion only for a finite number of $1$-parameter subgroups $\lambda_1,\dots,\lambda_k$ (independent of $t$). 
\end{itemize}
\end{remark}

From (\ref{temp3}) and (\ref{temp4}), we note that  $\mu^t(\Xi,r)$ depends only on the ``maximal (or support) monomials''  
of $\Xi$. To make this precise, we introduce the following notions (see \cite[Ch. 7]{mukai}):

\begin{definition}\label{order}
Let $m=x_0^ax_1^bx_2^c$ be a monomial, and $r$ a rational number.  We denote $\langle m , r\rangle:=a+b r-c (1+r)$ and define a partial ordering on the set of  monomials of a given degree by $m>m'$ iff $m\neq m'$ and $\langle m , r\rangle\ge\langle m' , r\rangle$ for all $r\in[-\frac{1}{2},1]$.  For a set of degree $k$ monomial $\Xi_k$ we define {\it the support} $\Supp(\Xi_k)$ as the subset of maximal monomials with respect to this partial ordering. Similarly, for $\Xi=(\Xi_d,\Xi_1)$ we set $\Supp(\Xi):=(\Supp(\Xi_d),\Supp(\Xi_1))$.
\end{definition} 

\subsubsection{Algorithmic Description of Stability}\label{sectalgo} 
\begin{figure}[htb]
\begin{center}
\scalebox{.50}{\includegraphics{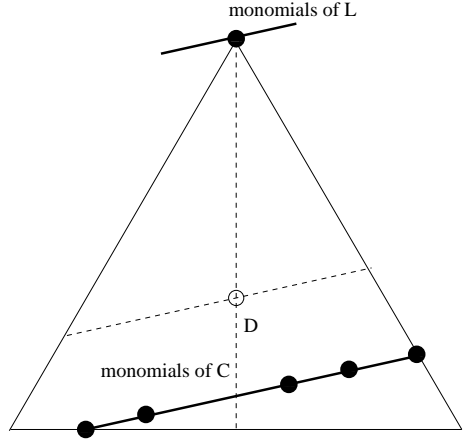}}
\end{center}
\caption{Configurations of monomials and stability}
\label{convexspan}
\end{figure}

For a fixed slope $t$, the analysis of stability for pairs is roughly equivalent to the analysis of stability for degree $d+t$ curves (see \S\ref{caset1} below). This case is well known (\cite[\S1.9]{mumford}, \cite[\S7.2]{mukai}). Thus, we are essentially done once we have reduced the analysis  to a finite number of critical values of $t$. This is achieved by the following algorithm:

\begin{itemize}
\item[Step 1] ({\it Find the critical slopes $t_i$}): For a given degree $d$, the number of configurations $\Xi$ is finite. For each $\Xi$, $\mu^t(\Xi,r)$ is continuous, piecewise linear in $r$, with critical points depending only on $\Xi$. It suffices to apply the numerical criterion only for the critical points $\{r_0=-\frac{1}{2},\dots, r_k=1\}$. Given $\Xi$ and a critical point $r_i$, $\mu^t(\Xi,r_i)$ is linear in $t$, changing sign at most once. Thus, the set of critical slopes is included  in the set $T:=\{t\mid \mu^t(\Xi,r_i)=0, \mu^{t-\epsilon}(\Xi,r_i)\neq0 \textrm{ for some } \Xi \textrm{ and } r_i\}$. A post-processing step removes the irrelevant slopes from $T$ (e.g. $t\in T$ with $t>\frac{d}{2}$). The effectivity follows from the fact that once $\Xi$ is fixed, the algorithm is linear, and  the outcome depends only on $\Supp(\Xi)$. The number of possible supports is of order $O(2^d)$.

\item[Step 2] ({\it For each $t_i$ find the maximal relevant configurations}): This is the standard analysis of stability for plane curves (compare \cite[1.11]{mumford} with figure \ref{convexspan}). Essentially, one has to consider only the configurations $\Xi$ such that $\mu^{t_i}(\Xi, r)=0$ for some $r\in [-\frac{1}{2},1]$.

\item[Step 3] ({\it Interpret geometrically the results of Step 2}): Again, this is the standard analysis of stability (e.g. \cite[pg. 81--82]{GIT}). The only slight difference is that we have to analyze in addition to $C$ also the relative position of $L$ (see remark \ref{geominterpretation}).
\end{itemize}
The first two steps are purely combinatorial, and are easily implemented. In contrast, the last step requires a careful geometric analysis, which is possible only for low degree pairs. 

\begin{remark}\label{graphrule} 
We can visualize the stability condition by drawing  the monomials of a configuration $\Xi=(\Xi_d,\Xi_1)$ in a triangle in the plane as in Mumford \cite[\S 1.9]{mumford}.  For the monomials in $\Xi_1$, we use the vertices of the triangle in the obvious way. The stability of $\Xi$ depends only on the support monomial of $\Xi_1$ and the boundary $\Gamma$  of the convex span of the support monomials of $\Xi_d$ (i.e. the Newton diagram). The configuration is unstable at $t=0$ (the plane curve case) iff the center of the triangle lies above the boundary $\Gamma$. More generally for slope $t$, we have the same rule, but $\Gamma$ is translated with $\frac{2t}{3}$ units in the direction of the support monomial of $\Xi_1$. An example is given in  figure \ref{convexspan}. There, $\Xi$ is $t$-semistable iff $t\ge\frac{3D}{2}$.
\end{remark}

\begin{remark}\label{geominterpretation} For the geometric analysis of stability, an  important role is played by the observation that a $1$-PS $\lambda$ singles out a bad flag $p_\lambda\in L_\lambda$ (\cite[pg. 82]{GIT}). With our convention on coordinates, we have $p_\lambda=(1:0:0)$ and $L_\lambda\equiv(x_2=0)$. A simple useful fact  is that the relative position of the flag $(p_\lambda,L_\lambda)$ with respect 
to $L$ determines the support of $\Xi_1$:
$$\Supp(\Xi_1)=
\begin{cases}
x_0 & \textrm{ iff } p_\lambda\notin L\\
x_1 & \textrm{ iff } p_\lambda \in L \text{ but }  L\not\equiv L_\lambda\\
x_2 & \textrm{ iff } L\equiv L_\lambda
\end{cases}.
$$ 
\end{remark}

\subsection{Proof of Theorem \ref{mainthm1}}\label{proofmainthm1} The theorem follows from Cor. \ref{intstab} and the analysis of stability condition at the slopes $0$, $1$, and $\frac{d}{2}$.  The notations are those of \S\ref{numcriterion}. In particular, $x=(c,l)$ is a point representing the pair $(C,L)$.

\subsubsection{Slope $t=0$ case} Since $\mu^0(x,\lambda)=\mu(c,\lambda)$, the statement i) of \ref{mainthm1} follows immediately from the numerical criterion. \qed

\begin{remark}\label{structure0}\label{caset0}
The linearization of slope $0$ defines the projection $X=\bP^N\times \bP^2\xrightarrow{\pi_1} \bP^N$. Thus, $\calM(0)$ is precisely the GIT quotient corresponding to degree $d$ plane curves. The variation of quotients morphism 
$\calM(\epsilon)\xrightarrow{\pi} \calM(0)$ can be interpreted as the forgetful map $(C,L)\to C$ which, over the locus of curves with trivial stabilizer,  is  a $\bP^2$-bundle in the \'etale topology.
\end{remark}

\subsubsection{Slope $t=1$ case}\label{caset1} Assume more generally that $t$ is any positive integer. We claim that {\it $x$ is $t$-semistable iff  $c\cdot l^t$ is semistable as a degree $d+t$ homogeneous form}.  This follows from the numerical criterion and the identity:
\begin{equation}\label{equationdt}
\mu^t(x,\lambda)=\mu(c\cdot l^t,\lambda). 
\end{equation}

Alternatively, we can see ii) of \ref{mainthm1} more intrinsically  as follows. Let $V$ be the standard representation of $G=\SL(3)$. By definition, we have $X=\bP(\Sym^d(V^*))\times \bP(V^*)$. The choice of linearization $\calL=\calO(1,t)$ of slope $t$ gives the embedding: 
\begin{equation*}
X=\bP(\Sym^d(V^*))\times \bP(V^*)\hookrightarrow\bP\left(\Sym^d(V^*)\otimes \Sym^t(V^*)\right)
\end{equation*}
(a composition of the Veronese and Segre embeddings). By definition, the stability with respect to the linearization $\calL$ is the stability with respect to the linear action of $G$ on $\bP\left(\Sym^d(V^*)\otimes \Sym^t(V^*)\right)$. The representation of $G$ on $\Sym^d(V^*)\otimes \Sym^t(V^*)$ is reducible with the top summand $\Sym^{d+t}(V^*)$ determined by the  multiplication map: 
\begin{eqnarray*}
\Sym^d(V^*)\otimes \Sym^t(V^*)&\xrightarrow{\pi}& \Sym^{d+t}(V^*)
\end{eqnarray*}
The conclusion now follows by noting that the affine cone over the image of $X$ consists of pure tensors, none lying in the kernel of the projection $\pi$. \qed 

\begin{remark}
The morphism $X\xrightarrow{j} \bP\left(\Sym^{d+t}(V^*)\right)$ constructed above is a normalization onto the image.  The statement descends also to the  GIT quotients.
\end{remark}

\subsubsection{Slope $t=\frac{d}{2}$ case}\label{cased2}
We note first that there are no semistable points for $t>\frac{d}{2}$.   
\begin{lemma}\label{lemmastabd2}
If $t>\frac{d}{2}$ then $X^{ss}(t)=\emptyset$. Similarly, $X^{s}(\frac{d}{2})=\emptyset$.
\end{lemma}
\begin{proof}
Let $x=(c,l)\in X$. Choose coordinates such that $l=x_2$.  Let  $\lambda$ be the diagonal $1$-PS of weights $r_0=r_1=1$ and $r_2=-2$. From (\ref{temp3}) and (\ref{temp4}), we get  
$$\mu(l,\lambda)=r_2=-2 \textrm{ and }\mu(c,\lambda)\le d,$$
and then
$$\mu^t(x,\lambda)=\mu(c,\lambda)+t\mu(l,\lambda)<0,$$ 
i.e. $x$ is unstable.
\end{proof}

On the other hand, there exist $t$-semistable points for $t=\frac{d}{2}$. 
\begin{lemma}\label{lemmadiscriminant}
Let $(C,L)$ be a degree $d$ pair. Assume that $L$ is transversal to $C$, then $(C,L)$ is semistable at $t=\frac{d}{2}$.
\end{lemma}
\begin{proof}
We consider the following discriminant hypersurface in $X$: 
$$\Sigma_1=\{(C,L)\mid L \textrm{ is tangent to } C \textrm{ or worse}\}.$$
Clearly, $\Sigma_1$ defines a $G$-invariant divisor in $X$.  Thus, $\Sigma_1$ gives an invariant section $\sigma_1\in H^0(X,\calO(a,b))^G$ for some $(a,b)$. The invariant section $\sigma_1$ is non-vanishing exactly when $L$ is transversal to $C$. Therefore, such pairs will be semistable at $t=\frac{b}{a}$. It remains to find $(a,b)$, the bidegree of $\Sigma_1$.  Fixing $C$ generic, we find that $b$ is the degree of the dual curve $\check C$, i.e. $b=d(d-1)$. Similarly, $a=2(d-1)$ is the degree of the discriminant for degree $d$ binary forms. We obtain $t=\frac{d}{2}$ and the lemma follows.
\end{proof}

We now conclude the proof of theorem \ref{mainthm2} by the following lemma.
\begin{lemma}\label{stabd2}
A pair $(C,L)$ is semistable at $t=\frac{d}{2}$ if and only if $L$ is not a component of $C$ and $\mult_p(C\cap L)\le \frac{d}{2}$ for every $p\in C\cap L$.
\end{lemma}
\begin{proof}
The ``only if'' part follows as in lemma \ref{lemmastabd2}. Namely, given a point $p$ with $k=\mult_p(C\cap L)> \frac{d}{2}$,  we choose coordinates such that $p=(1,0,0)$ and $L$ is given by $(x_2=0)$. It is then easy to see that there exists a  choice of weights (e.g. $r =\frac{k-2}{k+1}$) such that the resulting diagonal $1$-PS $\lambda$ destabilizes the pair $(C,L)$.

Conversely, assume that $L\not \subset C$ and $\mult_p(C\cap L)\le \frac{d}{2}$ for all $p$. We claim that the pair is semistable at $\frac{t}{2}$. If we suppose not, then we obtain a contradiction as follows. Choose  a destabilizing $1$-PS $\lambda$. Normalizing as in \ref{defnormalize}, we replace the pair $(C,L)$ by a configuration of monomials $\Xi=(\Xi_d,\Xi_1)$ and $\lambda$ by its normalized weight $r\in[-\frac{1}{2},1]$. The assumption that $\lambda$ destabilizes $(C,L)$ is equivalent to 
\begin{equation}\label{temp5}
\mu^{\frac{d}{2}} (\Xi,r)=\mu(\Xi_d,r)+\frac{d}{2}\mu(\Xi_1,r)<0.
\end{equation}

The proof now consists in analyzing the inequality (\ref{temp5}) and deducing geometric consequences.  Depending on the monomials occurring in $\Xi_1$ we distinguish $3$ cases: 
\begin{itemize}
\item[Case] ($\Supp(\Xi_1)=\{x_0\}$): From (\ref{temp4}) and the normalization procedure, we get $\mu(\Xi_1,r)\equiv 1$. Thus, 
\begin{equation}\label{temp6}
\mu(\Xi_d,r)<-\frac{d}{2} \textrm{ for some } r\in \left[-\frac{1}{2},1\right]
\end{equation}
The function $\mu(\Xi_d,r)$ is computed by (\ref{temp3}) with the weights  normalized by $r_0=1$, $r_1=r$, and $r_2=-1-r$. A simple analysis gives that all the degree $d$ monomials contained in $\Xi_d$ must be divisible by $x_2^{\lfloor \frac{d}{2}\rfloor}+1$. Geometrically, this means that $C$ contains a line with multiplicity strictly larger than $\frac{d}{2}$, contradicting the assumption on the intersection $C\cap L$.
\end{itemize}

 The argument in the remaining cases is similar. We obtain the following contradictions to the hypothesis: 
\begin{itemize}
\item[Case] ($\Supp(\Xi_1)=\{x_1\}$): $L$ passes through a point $p$ with $\mult_p(C)>\frac{d}{2}$.
\item[Case] ($\Supp(\Xi_1)=\{x_2\}$): $L$ is tangent with multiplicity larger than $\frac{d}{2}$ to $C$.
\end{itemize}
\end{proof}

From the previous lemma and general results of the variation of GIT quotients, we obtain the following result regarding the structure of the GIT quotient at $t=\frac{d}{2}$.
\begin{corollary}\label{structured2}
The space $\calM(\frac{d}{2})$ is naturally isomorphic to $\Sym^d(\bP^1)\gquot \SL(2)$, the moduli of unordered $d$-tuples in $\bP^1$.  The variation of quotients morphism $\calM(\frac{d}{2}-\epsilon)\xrightarrow{\pi} \calM(\frac{d}{2})$
can be interpreted as the forgetful map $(C,L)\to C\cap L\subset L\cong \bP^1$. Furthermore, for $d\ge 5$  $\calM(\frac{d}{2}-\epsilon)\to \calM(\frac{d}{2})$ is generically a weighted projective bundle in the \'etale topology.
\end{corollary}
\begin{proof}
We define a map $\Sym^d(\bP^1)\to X=\bP(\Sym^d)\times \bP^2$ by associating to a $d$-tuple of points in $\bP^1$ the projective cone $C$ over it together with a line $L$ not passing through the vertex. We make two basic observations:
\begin{itemize}
\item[(1)] The pair $(C,L)$ is semistable at $t=\frac{d}{2}$ (cf. Lemma \ref{stabd2}). The orbit of any pair $(C',L')$ with $C'\cap L'\cong C\cap L$ (as $d$-tuples) contains in its closure the orbit of $(C,L)$.
\item[(2)] The pair $(C,L)$ is stabilized by a $1$-PS $\lambda$ of weights $(1,1,-2)$. Generically, for $d\ge 5$, the stabilizer of $(C,L)$ is precisely $\lambda$.
\end{itemize}
The first item guarantees that the induced map $\Sym^d(\bP^1)\gquot \SL(2)\xrightarrow{j}  \calM(\frac{d}{2})$ is well defined and surjective. From (2) and a theorem of Luna \cite[Main Thm.]{luna} it follows that $j$ is actually a finite morphism. Clearly, $j$ has degree $1$. Thus, since both the source and the target of $j$ are normal varieties,  the morphism $j$ is an isomorphism.

The statement about the structure of the variation of quotients morphism follows from Thaddeus \cite[Thm. 5.6]{thaddeus} and Dolgachev-Hu \cite[Thm. 4.2.7]{dolgachevhu}.
\end{proof}

\bigskip

\bigskip

\subsection{Relation to the log canonical threshold}\label{secbounds}
\subsubsection{Proof of Theorem \ref{mainthm2}}\label{proofmainthm2}
The case $t=0$ was previously established by Hacking \cite[\S 10]{hacking} and  Kim-Lee \cite{kimlee}. The general case is essentially the same. Namely, in \S\ref{caset1}, we noted that the identity (\ref{equationdt}) implies that a pair $(C,L)$ is $t$-semistable if and only if $C+tL$ satisfies the numerical criterion for degree $d+t$ curves (N.B. the condition is purely numerical, so it makes sense for $t\in \bQ$). Assume that $(C,L)$ is $t$-unstable. By the numerical criterion, we can find coordinates $(x_0:x_1:x_2)$ on $\bP^2$ and relatively prime integral weights $w_1$, $w_2$ such that with respect to the affine coordinates $x=\frac{x_1}{x_0}$ and $y=\frac{x_2}{x_0}$ at $p=(1:0:0)$ we have 
$$w(f)>\frac{d+t}{3}(w_1+w_2),$$
where $f(x,y)$ is the equation of $C+tL$ and $w(f)$ denotes the order of $f$ at $p$. By considering the weighted blow-up  of $\bP^2$ with respect to  the weights $w_1$, $w_2$, we obtain via a standard discrepancy computation (see \cite[6.38]{kollarrational}) that 
$$a(E,\bP^2,\frac{3}{d+t}(C+tL))=\left(w_1+w_2-1-\frac{3}{d+t}\cdot w(f)\right)<-1,$$
where $E$ is the exceptional divisor of the blow-up and $a(\dots)$ denotes the discrepancy of  $E$ (see \cite[\S2.3]{kollarmori}). We conclude that $(\bP^2,C+tL)$ is not log canonical. \qed

\subsubsection{The stability threshold}\label{sectbounds} 
Let $(C,L)$ be a degree $d$ pair and $[\alpha,\beta]\subset [0,\frac{d}{2}]$ its interval of stability. As a consequence of theorem \ref{mainthm2}, it follows that the endpoints $\alpha$ and $\beta$ of the interval of stability are essentially determined by the worst singularity of $C$ and by  the point of highest multiplicity for the  intersection $C\cap L$  respectively. To simplify the computations involved in the complete analysis of stability for a given degree  $d$ (e.g. $d=5$ in section \ref{sectdeg5}), it is convenient to make the statement of theorem \ref{mainthm2} more explicit. Specifically, we measure the effect of the singularities of $C$ on the stability of the pair $(C,L)$ by the following notion: 

\begin{definition}\label{defstab}
Let $C$ be a plane curve of degree $d$ and $p\in C$ a singular point. We define {\it the stability threshold of $p\in C$}, denoted by $t_p(C)$, to be
\begin{equation}\label{deftp1}
t_p(C):=-\inf_{\lambda \textrm{ adapted to }p}\frac{\mu(c,\lambda)}{||\lambda||}
\end{equation}
where  a $1$-PS $\lambda$ is {\it adapted to } $p$ iff $p=p_\lambda$, the norm of $\lambda$ is $||\lambda||:=r_0$ (the highest weight of $\lambda$), and $c$ is the equation of $C$. 
\end{definition}

Note that the stability threshold does not depend on $L$. We then have: 
\begin{lemma}
Let $(C,L)$ be a degree $d$ pair with interval of stability $[\alpha,\beta]$. Then
$$\alpha\le\max \left\{\max_{p\in \Sing(C)}t_p(C),0\right\}$$
with equality if $L$ is transversal to $C$. \qed
\end{lemma}

The stability threshold of a singular point is related to a well-known invariant of the singularity, the log canonical threshold (see \cite{kollarpairs} and \cite[\S6.5]{kollarrational}).
\begin{definition}\label{deflinquas}
Let $p\in C$ be a singular point. We say that the singularity at $p$ is {\it linearly  semi-quasihomogeneous} if there exists a choice of  homogeneous coordinates $(x_0:x_1:x_2)$ and a choice of weights $w_1$ and $w_2$ such that the associated affine equation $f(x,y)$ of $C$ at $p$ is semi-quasihomogeneous (i.e. the leading term $f_w$ defines an isolated singularity at the origin).
\end{definition}

\begin{corollary}\label{logcan}
The following inequality holds
\begin{equation}\label{formulalog}
t_p(C)\le \frac{3}{\lc_p(C)}-d
\end{equation}
where $\lc_p(C)$ denotes the log canonical threshold of $p\in C$. The equality holds if 
the singularity at $p$ is linearly semi-quasihomogeneous. 
\end{corollary}
\begin{proof}
The claim follows by same argument as \ref{mainthm2}. If $p$ is a semi-quasihomogeneous singularity, the log canonical threshold is computed by $\lc_p(C)=\frac{w_1+w_2}{w(f)}$ (\cite[Prop. 8.14]{kollarpairs}). A similar formula holds for $t_p(C)$ under the assumption of linearity.
\end{proof}

For low degree pairs,  the equality in (\ref{formulalog}) almost always holds. We note, however, that starting with degree $4$ there are examples in which this fails:
\begin{example}[The stability threshold is not a local analytic invariant]\label{nonanalytic}
 Consider the following singularities of plane quartics: $C_1: (x^2+xy^3=0)$ and $C_2: ((x-y^2)^2-x^2y^2=0)$. Analytically, they both define an $A_5$ singularity at origin $p$. Thus, the log canonical threshold in both cases is $\frac{2}{3}$. It is easy to compute that $t_p(C_1)=\frac{1}{2}$ and respectively $t_p(C_2)=0$. It follows that the equality in (\ref{formulalog}) holds only in the first case.  The two situations are not distinguished locally, but globally. Namely, $C_1$ consists of a line tangent with multiplicity $3$ to a cubic, and $C_2$ is the union of two conics tangent with multiplicity $3$ in the point $p$.
\end{example}

The following bounds for the stability threshold in terms of the multiplicity of the singularity are easily established. 
\begin{proposition}\label{boundsmult}
Assume that $p$ is a point of multiplicity $k$ for a degree $d$ curve $C$. Then the 
following bounds hold:
$$\frac{3k}{2}-d\le t_p(C)\le 3k-d.$$
Moreover
\begin{itemize}
\item[i)] $t_p(C)=\frac{3k}{2}-d$ iff every line in the tangent cone at $p$ has multiplicity at most $\frac{k}{2}$;
\item[ii)]  $t_p(C)\le 3k\cdot\frac{d-1}{d+k-2}-d$ if $p$ is an isolated singularity;
\item[iii)] $t_p(C)=3k-d$ iff $C$ contains a line with multiplicity $k$ passing through $p$. \qed
\end{itemize} 
\end{proposition}

\begin{remark}\label{remarkcompute}
We close by noting that the stability threshold is easily computed in practice. Namely, let $p$ be an isolated singular point of multiplicity $k$ for $C$. By looking at the tangent cone, we have that either every tangent line occurs with multiplicity at most $\frac{k}{2}$  in the tangent cone, or there exists a unique {\it special tangent} $L_0$ with multiplicity larger than $\frac{k}{2}$.
The first case is covered by the previous proposition. In the second case, it is easily seen that the $1$-PS $\lambda$ that computes the stability threshold $t_p(C)$ has the property that $p_{\lambda}=p$ and $L_{\lambda}=L_0$. By choosing coordinates such that $p=(1:0:0)$ and the special tangent $L_0$ is given by $(x_2=0)$, we obtain a configuration  $\Xi=(\Xi_d,\Xi_1)$ of monomials. The stability threshold is then given by $t_p(C)=-\min_{r\in[-\frac{1}{2},1]} \mu(\Xi_d,r)$, which is easily computed (Remark \ref{graphrule}). The essential observation here is that, while $\Xi$ depends on the choice of coordinates, its support $\Supp(\Xi)$ does not.
\end{remark}
 
\subsubsection{Upper bounds for the interval of stability}\label{sectbounds2}
We can apply similar considerations for the analysis of the failure of stability due to a point $p\in C\cap L$. Note that  the singularity at $p$ of the $\bQ$-divisor $C+tL$ becomes worse as $t$ increases.  We obtain the following bounds for the end-point $\beta$ of the interval of stability.
\begin{proposition}\label{upbound}
Let $p$ be the point with the highest multiplicity in the intersection $C\cap L$. 
Then, the upper bound $\beta$ of the interval of stability of the pair $(C,L)$  satisfies the following estimates:
\begin{itemize}
\item[i)] if $L$ is a component of $C$ then $\beta\le \frac{d-3}{2}$;
\item[ii)] if $L$ is not a component of $C$ and $\mult_p(C\cap L)\le \frac{d}{2}$ then  $\beta=\frac{d}{2}$;
\item[iii)]  If $L$ is not a component of $C$ and $\mult_p(C\cap L)=k> \frac{d}{2}$ then
$$\frac{d}{2}-\frac{3}{2}(2k-d)\le \beta \le \frac{d}{2}-\frac{3(2k-d)}{2(2k-1)}.$$
\end{itemize} 
Furthermore, the above estimates are sharp. \qed
\end{proposition}

\subsection{Relation to the deformations of non-positive weight}\label{sectpinkham}
The motivation for the study of the moduli space of pairs comes from Pinkham's theory of deformations of singularities with $\bC^*$-action. The basic idea of this theory is that in the presence of a $\bC^*$-action a certain subspace  of the deformation space, {\it the deformations of negative weight}, can be globalized and interpreted as a moduli space of pairs. This modular interpretation gives effective tools for the study of the deformations of certain classes of singularities, such as the unimodal singularities (e.g. \cite{pinkham, pinkhammonodromy,looijenga,looijengatriangle,novaacta}). For the general theory, we refer to the work of  Pinkham \cite{pinkham0,pinkhamc*} (for a short exposition see \cite[Appendix]{looijengatriangle}). Here, we briefly recall the basic concepts of the theory, and explain the relation between the deformations of cones over $d$-tuples of points in $\bP^1$ and the variation  of GIT quotients for degree $d$ pairs.

Let $(Y_0,y_0)$ be the germ of a  singularity with good $\bC^*$-action in the sense of Pinkham (e.g. an isolated quasihomogeneous hypersurface singularity). Then there exists a formal versal deformation $\calY\to S$ such that $\bC^*$ acts equivariantly on $\calY$ and $S$. This determines subspaces $S_-$, $S_0$, $S_{\le 0}$ and  $S_{+}$  and  pull-back families (e.g. $\calY_-\to S_-$, {\it the deformations of negative weight}). Each of these spaces has an intrinsic interpretation (see \cite{pinkhamc*}). In particular, the zero weight deformations are those preserving the $\bC^*$-action, and the deformations of non-positive weight are the deformations that can be lifted to projective deformations of a completion  $\overline{Y}_0$ of the singularity. Specifically, we can assume $Y_0=\Spec A$ for a graded ring $A$. Then $\overline{Y}_0=\Proj A[t]$  (for $\deg(t)=1$) is a natural compactification of $Y_0$. The divisor $Y_{\infty} \cong \overline{Y}_0\setminus Y_0$ defined by $t$ is called the {\it divisor at infinity}. This construction works well in families and  identifies the deformations of non-positive weight to the projective deformations of the pair $(\overline{Y}_0,Y_\infty)$ (\cite[Thm. 2.9]{pinkhamc*}). The deformations of negative weight correspond to the deformations of $(\overline{Y}_0,Y_\infty)$, fixing the hyperplane at infinity $Y_{\infty}$. 

The key observation of Pinkham  is that, due to the $\bC^*$-action, the deformations of negative weight can be globalized in the following sense. The functor of deformations of negative weight is representable by an affine space $S_{-}$. The resulting weighted projective space $S_{-}\gquot \bC^*$  is then the moduli of  pairs $(\overline{Y},Y_{\infty})$, where $\overline{Y}$ is a projective deformation of $\overline{Y}_0$ and $Y_{\infty}$ is the fixed hyperplane section (for a functorial formulation see \cite[Appendix]{looijengatriangle}). The globalization does not hold in the zero weight direction, but we can interpret $S_0$ as corresponding to the deformations of the hyperplane at infinity $Y_{\infty}$. Also, there exists a natural $\bC^*$-equivariant map $S_{\le 0}\to S_0$, whose fibers correspond to the deformations of negative weight. 

The construction explained in the previous paragraphs is easily understood in the case of deformations of quasihomogeneous hypersurface singularities. For instance, in our situation, $Y_0$ is the affine cone over a $d$-tuple of points in $\bP^1$ and $\overline{Y}_0=:C_0$ is the natural projective cone. The hyperplane at infinity is simply a transversal line $L$ to $C_0$. It follows that the deformations of negative weight of the cone $Y_0$ correspond to the pairs $(C,L)$ with $C\cap L$ projectively equivalent to $C_0\cap L$. Similarly, the non-positive weight deformations correspond to the pairs $(C,L)$ such that $C\cap L$ is a small deformation of $C_0\cap L$. 

Since a degree $d$ pair $(C,L)$ such that $C\cap L$ is transversal and $C$ is not a cone is GIT stable at $t=\frac{d}{2}-\epsilon$ (cf. Thm. \ref{mainthm1}),  we can interpret $\calM(\frac{d}{2}-\epsilon)$  as a global object associated to the non-positive weight deformations $S_{\le 0}$. Furthermore, the variation of GIT map $\calM(\frac{d}{2}-\epsilon)\to \calM(\frac{d}{2})$ can be interpreted as a globalization of the map $S_{\le 0}\to S_0$ modulo the $\bC^*$-action. To be more precise, let us recall that a point $x_0\in \calM(\frac{d}{2})$ corresponds to the closed orbit of a pair $(C_0,L)$ consisting of a degree $d$ cone and a transversal line. The stabilizer $G_{x_0}$ of the pair $(C_0,L)$ is $\bC^*$ (at least generically, for $d\ge 5$). Let $S_{\le 0}$ be a normal slice to the orbit of $(C_0,L)$ and $S_0$ the invariant part. By Luna's slice theorem,  locally at $x_0$ in the \'etale  topology, we have the following commutative diagram: 
\begin{equation}\label{defsquare}
\begin{CD}
(S_{\le 0}\setminus S_0)/\bC^*@>>>\calM(\frac{d}{2}-\epsilon)\\
@VVV @VVV\\
S_0@>>> \calM(\frac{d}{2})
\end{CD}
\end{equation}
with the horizontal maps being \'etale. The stabilizer $G_{x_0}\cong \bC^*$ acts naturally on the tangent space $T_{S_{\le 0}}$. It is then a simple (almost tautological) computation to identify  $T_{S_{\le 0}}$ with the non-positive weight subspace of  $\Ext^1(\Omega_{Y_0},\calO_{Y_0})\cong\bC\{x,y\}/\langle J(f_d)\rangle$, where $f_d(x,y)$ is the equation of $Y_0$ and $J(f_d)$ is the Jacobian ideal. It follows that the germ of $S_{\le 0}$ at the origin is indeed the space of deformations of non-positive weight, and similarly for $S_0$. In conclusion, (\ref{defsquare}) corresponds indeed to a globalization of the natural map of deformations $S_{\le  0} \to S_0$. Moreover, the weighted projective fibers of $\calM(\frac{d}{2}-\epsilon)\to \calM(\frac{d}{2})$ (cf. Cor. \ref{structured2}) correspond to the negative weight deformations $S_{-}$ modulo $\bC^*$.

\begin{remark} The singularity $N_{16}$ is quasi-homogeneous with Milnor (and Tyurina) number $\mu=16$. The weights for $N_{16}$ are as follows:  one is positive, two are zero, and $13$ are negative. Since the only   positive weight is in the hessian direction, the positive weight deformations are topologically trivial and, as  is customary, we ignore them. 
\end{remark}

\section{The stability conditions for degree $5$ pairs}\label{sectdeg5} In this section, we do a detailed analysis of the stability conditions for degree $5$ pairs. The first step of this analysis is the determination of  the relevant critical slopes based on the algorithm described in \S\ref{sectalgo}.

\begin{lemma}
The critical slopes for degree  $5$ pairs are: 
$0$, $\frac{1}{7}$, $\frac{1}{4}$, $\frac{2}{5}$, $\frac{5}{8}$, $1$, $\frac{10}{7}$, $\frac{8}{5}$, $\frac{5}{3}$,
$\frac{7}{4}$, $\frac{13}{7}$, $2$, $\frac{11}{5}$, and $\frac{5}{2}$. 
\end{lemma}

The semi-stability of the pair $(C,L)$ for the slopes $0$, $1$, and $\frac{5}{2}$ is equivalent to the semi-stability of $C$, $C+L$, and $C\cap L$  respectively (Thm. \ref{mainthm1}). Thus, the stability condition for the slopes $0$, $1$, and $\frac{5}{2}$ is well understood.  In particular, as a corollary of the results of Shah on the stability of plane sextics (\cite[Thm. 2.3]{shah}), we obtain: 
 
\begin{corollary}\label{corquintic}
Let $(C,L)$ be a degree $5$ pair. Assume that the sextic $C+L$ is reduced. Then the pair $(C,L)$ is stable (semistable)  at $t=1$ if and only if $C+L$ has at worst simple (resp. simple elliptic or cusp) singularities. 
\end{corollary}

The description of the stability for the remaining critical slopes is a standard GIT computation. For degree $5$ pairs, the computation is simplified by the fact that the interval of stability can be determined by considering independently the worst singularity of $C$ and of the intersection $C\cap L$. It follows that the stability conditions are essentially determined by the results of \S\ref{sectbounds} and \S\ref{sectbounds2}. The following statement summarizes the results of this section.

\begin{theorem}\label{rule5}  Let  $(C,L)$ be a degree $5$ pair. The interval of stability $[\alpha,\beta]\subseteq [0,\frac{5}{2}]$ is determined by the following rules:
\begin{itemize}
\item[i)] Assume that $L$ does not pass through a singularity that makes $C$ unstable. Then $\alpha$ is determined by the worst singularity of $C$, i.e. 
$$\alpha=\max_{p\in \Sing(C)} \{t_p(C),0\}.$$ 
Similarly,  $\beta$ is determined by the worst intersection point of $C$ and $L$ (see \S\ref{stabred5} and \S\ref{nongen5}  respectively).
\item[ii)] Assume that $L$ passes through a singularity that makes $C$ unstable. Then the interval of stability of the pair is either empty or consists of a single point $t=1$ (see \ref{computationalpha}).
\end{itemize}
\end{theorem}

In particular, it follows that the change of the stability condition is very simple for all critical slopes except $t=1$.  Essentially, for such a slope, exactly one class of singularities becomes stable, and one type of degenerate intersection becomes unstable (see \S\ref{symquintic}). The change of stability at $t=1$ is more involved (see \S\ref{quintic1}), but it follows a pattern as noted by the following corollary:

\begin{corollary}\label{corquintic2}
Let $(C,L)$ be a degree $5$ pair. Assume that the intersection $C\cap L$ is transversal and that $C$ is not a cone. Then the interval of stability of the pair $(C,L)$ is $[\alpha,\frac{5}{2}]$ with:
\begin{itemize}
\item[i)] $\alpha\in[0,1)$ iff the quintic $C$ has at worst simple singularities;
\item[ii)] $\alpha=1$ iff the quintic $C$ has a simple elliptic or cusp singularity;
\item[iii)] $\alpha\in(1,\frac{5}{2})$ iff the quintic $C$ has a triangle singularity.
\end{itemize}
\end{corollary}
The division of singularities in  three large classes as above has significance both topologically (in terms of monodromy) and  algebraically (i.e. canonical, log canonical, or worse singularities -- see remark \ref{remarkt1}). The corollary says that the division makes sense also in terms of GIT stability. The coincidence of these three different points of view plays a key role in the second part of our paper (see also the discussion from Mumford \cite[\S3]{mumford}).

\subsection{The singularities of plane quintics and stability}\label{stabred5} The effect of a singular point $p\in C$ on the stability of a pair $(C,L)$ is measured by the stability threshold $t_p(C)$ as described in \S\ref{sectbounds}. Here we are interested in the possible values for $t_p(C)$ in the case of plane quintics. 

\subsubsection{Isolated singularities case} We start by noting the following list of possible singularities for  a plane quintic.
\begin{proposition}[Wall \cite{wallquintic}] Let $C$ be a reduced plane quintic, and $p\in C$ a singular point.  Then  the singularity at $p$ is one of the following types:
\begin{itemize}
\item[-] a simple singularity $A_k$ ($k\le 12$), $D_l$ ($l\le12$), or $E_m$ ($m=6,7,8$);
\item[-] a simple elliptic singularity $\widetilde{E}_r$ ($r=7,8$);
\item[-] a plane cusp singularity $T_{2,3,k}$ ($k=7,\dots,10$) or $T_{2,q,r}$  ($4\le p\le q\le 6$);
\item[-] a triangle singularity of type $Z_{11}$, $Z_{12}$, $W_{12}$, or $W_{13}$;
\item[-] a trimodal singularity of type $N_{16}$.
\end{itemize}
 \end{proposition}

The simplest singularities of the list give stable plane quintics (see \cite[pg. 80]{GIT}), and thus the corresponding stability threshold is less than $0$ (Thm. \ref{mainthm1}). 

\begin{lemma}
A plane quintic with only singularities of type $A_k$, $D_4$, or $D_5$ is GIT stable. \qed
\end{lemma}

For the analysis of the remaining singularities, we need to separate the case of singularity of analytic type $D_8$ into two subcases: either of type $D_8'$ or not. 
\begin{definition}
Let $p\in C$ be a singular point of a plane quintic. We say that $p$ has type $D_8'$ iff the quintic $C$ decomposes as a line plus a nodal quartic such that the line is tangent with multiplicity $4$ to the quartic and $p$ is both the node of the quartic and  the intersection point of the two components.
\end{definition}

For many classes of singularities of  plane quintics, the stability threshold can be computed directly in terms of the log canonical threshold.

\begin{lemma}
Assume that $p\in C$ is an isolated singularity of one of the following types: $D_6$, $D_8'$, $E_k$ ($k=6,7,8$),  
$\widetilde{E_7}$, $\widetilde{E_8}$, $Z_{11}$, $Z_{12}$, $W_{12}$ or $W_{13}$. Then the singularity at $p$ is linearly semi-quasihomogeneous. In particular, the stability threshold is computed from the log canonical threshold by the formula $t_p(C)=\frac{3}{\lc_p(C)}-5$.
\end{lemma}
\begin{proof} The first claim is a case-by-case analysis of the analytic type of the singularity at $p$ based on the Newton diagram. The second part follows from  \ref{logcan}.
\end{proof}

The remaining cases are handled by the following two lemmas.

\begin{lemma}
Assume that $p\in C$ is a singularity of type $T_{2qr}$ (with $4\le q\le r\le 6$) for the quintic $C$. Then we have $t_p(C)=1$.
\end{lemma}
\begin{proof}
The singularities of type $T_{2qr}$ with $4\le q\le r\le 6$ have multiplicity $4$, but no line in the tangent cone has multiplicity larger than $2$. The conclusion follows from \ref{boundsmult}. 
\end{proof}

\begin{lemma}
Let $p$ be an isolated singular point of the quintic curve $C$.
\begin{itemize}
\item[i)] If $p$ is a singular point of type $D_{k}$ ($k\ge 7$), but not of type $D_8'$, $t_p(C)=0$.
\item[ii)] If $p$ is a singular point of type $T_{23k}$ ($k\ge 7$), $t_p(C)=1$.
\end{itemize}
\end{lemma}
\begin{proof}
The statement follows easily from the fact that we can choose affine coordinates at $p$ such that the leading term (w.r.t. appropriate weights) of the defining equation of $C$ is $f_w(x,y)=x y^2+2x^3y+x^5$ and $f_w(x,y)=y^3+2y^2x^2+x^4y$ respectively (see Remark \ref{remarkcompute}). 
\end{proof}

In conclusion, we obtain: 
\begin{proposition}\label{tpfor5}
Let $p\in C$ be an isolated singular point of a plane quintic. The stability threshold of $p\in C$ is computed by 
$$\renewcommand{\arraystretch}{1.25}
t_p(C)=
\left\{\begin{array}{ll}
\le 0& \textrm{if $p$ is of type } A_k,\textrm{or } D_l, \textrm{but not of type } D_8'\\
\frac{1}{7}, \frac{1}{4}, \frac{2}{5}, \frac{5}{8}& \textrm{if $p$ is of type } E_6, D_8', E_7, \textrm{or } E_8 \textrm{ respectively}\\
1& \textrm{if $p$ is of type } \widetilde{E_r}, \textrm{or } T_{2,p,q}\\
\frac{10}{7}, \frac{8}{5}, \frac{5}{3}, \frac{13}{7}&  \textrm{if $p$ is of type } Z_{11}, Z_{12}, W_{12}, \textrm{or } W_{13} \textrm{ respectively}\\
\frac{5}{2} & \textrm{if $p$ is of type } N_{16}
\end{array}
\right.
$$
In particular, $t_p(C)<1$ iff $p$ is a simple singularity, $t_p(C)=1$ iff $p$ is a simple elliptic or cusp singularity, and $t_p(C)>1$ otherwise.
\end{proposition}

\subsubsection{The stability conditions for non-reduced quintics}\label{stabnonred5} The stability of pairs $(C,L)$ such that $C$ is a non-reduced quintic is determined by the following two lemmas.
\begin{lemma}
Let $(C,L)$ be a degree $5$ pair. If $C$ contains a triple line, then the pair $(C,L)$ is unstable for all slopes $t$.  \qed
\end{lemma}

If $C$ is non-reduced and contains no line with multiplicity $3$, we can write $C=2D+R$, where $D$  is either a line or a smooth conic and $R$ is the residual curve having no common component with $D$. We then have:

\begin{lemma}\label{lemmanonred}
Let $C=2D+R$ be a non-reduced quintic containing no triple line. Then $\alpha=\max_{p\in \Sing(C)} t_p(C)$ is given by table \ref{redquintic}. \qed 
\end{lemma}

\begin{table}[htb]
\begin{center}
\renewcommand{\arraystretch}{1.25}
\begin{tabular}{|l|l|c|}
\hline
The double component $D$& Geometry of $D\cap R$& $\alpha$\\
\hline\hline
$D$ is a smooth conic &  $R$ is secant to $D$ & $0$\\
\hline
$D$ is a smooth conic & $R$ is tangent to $D$ & $1$\\
\hline\hline
$D$ is a line         & $|D\cap R|\ge 2$ & $1$\\ 
\hline
$D$ is a line         & $D\cap R=\{p\}$ and $p\in R$ is smooth &$\frac{7}{4}$\\ 
\hline 
$D$ is a line         & $D\cap R=\{p\}$ and $p\in R$ is of type $A_1$ &$2$\\ 
\hline 
$D$ is a line         & $D\cap R=\{p\}$ and $p\in R$ is of type $A_2$ &$\frac{11}{5}$\\ 
\hline 
$D$ is a line         & $D\cap R=\{p\}$ and $p$ is a triple point of $R$&$\frac{5}{2}$\\ 
\hline 
\end{tabular}
\vspace{0.2cm}
\caption{The stability of pairs $(C,L)$ with  $C=2D+R$}\label{redquintic}
\end{center}
\end{table}

\subsection{The stability for non-generic intersections}\label{nongen5} As a consequence of the relation between the GIT stability and log canonicity (see Thm. \ref{mainthm2}), the interval of stability of a pair $(C,L)$ is determined by the worst singularity $p$ of $C$ and the point $p'$ of maximal multiplicity in the intersection $C\cap L$. As long as $p\neq p'$, the contributions of the singularities of $C$ and of the intersection $C\cap L$ to the interval of stability are easily quantified.  The situation $p=p'$ (i.e. the line $L$ passes through a bad singularity of $C$)  is geometrically more subtle, but in the case of quintics it is easily handled by the following result:

\begin{lemma}\label{computationalpha}
Let $(C,L)$ be a degree $5$ pair with interval of stability $[\alpha,\beta]$. Then,
$$\alpha=\max_{p\in \Sing(C)} \{t_p(C),0\}$$
except when there exists a point $p\in C\cap L$ that makes $C$ unstable, in which case either:
\begin{itemize}
\item[i)] $\alpha=1$ if $C$ is the union of a quartic $C'$ and a $4$-fold tangent line $L'$  through a double point of $C'$ with $L\neq L'$  and $p=C'\cap L'$;
\item[ii)] $[\alpha,\beta]=\emptyset$ otherwise. 
\end{itemize}
If $L$ does not pass through a singular point making $C$ unstable, then $\beta$ can be determined by considering only $1$-PS $\lambda$ such that $p_{\lambda}$ is the point with the highest multiplicity in $C\cap L$. 
\end{lemma}
\begin{proof}
Assume that $\alpha>\max_{p\in \Sing(C)} \{t_p(C),0\}$. From the definition of the interval of stability, it follows that there exists  $1$-PS $\lambda$ such that $\mu^t(x,\lambda)\ge 0$ iff $t\ge \alpha$, where $x=(c,l)\in X$ represents the pair $(C,L)$. Since $\alpha>0$ and $\mu^t(x,\lambda)=\mu(c,\lambda)+t\mu(l,\lambda)$, we get $\mu(c,\lambda)<0$, i.e. 
$p:=p_\lambda$ is a singular point that destabilizes $C$. Since $\alpha>t_p(C)$, we also get $\mu(l,\lambda)<||\lambda||$, i.e. $p\in L$. Thus $L$ passes through a singularity which destabilizes $C$. By a computer-aided analysis, it follows that the only case when  the interval of stability is non-empty is described by the lemma.  The computation for $\beta$ is similar. 
\end{proof}

The following three lemmas determine the endpoint $\beta$ in terms of the geometry of the intersection $C\cap L$. We note the following two facts that simplify the computation:
\begin{itemize}
\item[(1)] We can assume that either $L$ is a component of $C$, or there exists a (unique) point $p$ with $\mult_p(C)\ge 3$. Otherwise, from theorem \ref{mainthm1}, it follows that $\beta=\frac{5}{2}$.
\item[(2)] We can assume that $L$ passes only through mild singularities of $C$. Otherwise, we are in the situation covered by the previous lemma. 
\end{itemize}

The pairs which have $\beta>1$ play an important role in the following sections. 
\begin{lemma}\label{lemmasimplesing}
Let $(C,L)$ be a degree $5$ pair with interval of stability $[\alpha,\beta]$. Then, $\beta>1$ if and only if $L$ is not a component of $C$ and any point $p$ with $\mult_p(C\cap L)\ge 3$ defines a simple singularity for the total curve $C+L$. 
The geometric possibilities for $p\in C\cap L$ to define a simple singularity for $C+L$ and the corresponding value of $\beta$ are given in table \ref{simplesing}. \qed
\end{lemma}

 \begin{table}[htb]
\begin{center}
\renewcommand{\arraystretch}{1.25}
\begin{tabular}{|c|c|c|c|}
\hline
$\mult_p(C\cap L)$& Sing. at $p\in C$ & Sing. at $p\in (C+L)$& $\beta$\\
\hline\hline
$k\in\{1,\dots,5\}$ & smooth  & $A_{2k-1}$ & $\frac{5}{2}$, $\frac{5}{2}$, $\frac{11}{5}$,  $\frac{13}{7}$, and $\frac{5}{3}$ resp. \\
\hline
$k\in\{2,\dots,5\}$ & $A_1$   & $D_{2k}$ &$\frac{5}{2}$,  $2$, $\frac{8}{5}$, and $\frac{10}{7}$ resp. \\
\hline
$3$ & $A_2$   & $E_7$      &$\frac{7}{4}$\\ 
\hline
$2$ & $A_n$   & $D_{n+3}$  &$\frac{5}{2}$\\ 
\hline 
\end{tabular}
\vspace{0.2cm}
\caption{The case $p\in (C+L)$ defines a simple singularity}\label{simplesing}
\end{center}
\end{table}

The cases where $L$ is a component of $C$ are described by the following lemma.
\begin{lemma}\label{linecomponent}
Let $(C,L)$ be a degree $5$ pair and let $[\alpha,\beta]$ be its interval of stability. Assume that $L$ is a component of $C$, and let $R$ be the residual curve. Then
\begin{itemize}
\item[i)] If $L$ is a component of $R$ then the pair $(C,L)$ is unstable for all $t$.
\item[ii)] If $L$ is not a component of $R$ then  we are in one of the following situations:
\begin{itemize}
\item[-] if $\mult_p(L\cap R)\le 2$ for all $p\in R\cap L$ then $\beta=1$;
\item[-] if $L$ passes through a singular point $p$ of $R$ and  $\mult_p(L\cap R)\ge 3$ then $\beta\le 0$;
\item[-] if $L$ is $4$-fold tangent to $R$ is a smooth point then $\beta=\frac{1}{7}$;
\item[-] if $L$ is $3$-fold tangent to $R$ is a smooth point then $\beta=\frac{2}{5}$. \qed
\end{itemize}
\end{itemize}
\end{lemma}

The remaining cases are covered by the following lemma. 
\begin{lemma}\label{lemmaline}
Let $(C,L)$ be a degree $5$ pair with interval of stability $[\alpha,\beta]$. Assume that $L$ is not a component of $C$ and that there exists a point $p$ with $\mult_p(C\cap L)\ge 3$. Assume additionally that the singularity at $p$ does not destabilize  $C$, and that $p$ is a non-simple singularity for $C+L$. Then $\beta=1$ with the exception of the cases listed in table \ref{before1}. \qed
\end{lemma}

 \begin{table}[htb]
\begin{center}
\renewcommand{\arraystretch}{1.25}
\begin{tabular}{|c|c|c|}
\hline
 Sing. of $p\in C$ & Position of $L$ through $p$& $\beta$\\
\hline\hline
 $A_7$  & $C=L+R$, $L$ is $4$-fold tangent to the residual quartic $R$& $\frac{1}{7}$ \\
\hline
 $D_5$   & $L$ is a special tangent through $p$ & $\frac{1}{4}$ \\
\hline
 $A_5$   & $C=L+R$, $L$ is $3$-fold tangent to the residual quartic $R$&$\frac{2}{5}$\\ 
\hline
 $A_4$   & $L$ is $5$-fold tangent through $p$ &$\frac{5}{8}$\\ 
\hline 
\end{tabular}
\vspace{0.2cm}
\caption{Intersections that destabilize $(C,L)$ before slope $t=1$}\label{before1}
\end{center}
\end{table}

We close the discussion of the GIT analysis by noting that, while the results for degree $5$ generally follow some predictable patterns, there are some pathological examples. For example, there exist pairs which are strictly semistable for an entire interval (compare to \cite[Appendix]{dolgachevhu}). 

\begin{example}\label{thickwall} Let $(C,L)$ be the pair of equations $C: \left(x_1(x_0^2x_2-x_1^2)^2=0\right)$ and $L: (x_1=0)$. The quintic $C$ is a double conic together with a secant line. Thus, $C$ is semistable and $\alpha=0$ (cf. \ref{lemmanonred}).  The line $L$ is a component of $C$ (it coincides with the secant line), but it does not destabilize the pair $(C,L)$ until $t=1$ (cf. \ref{lemmaline}). Thus the pair $(C,L)$ has interval of stability $[0,1]$, but $(C,L)$ is never stable (e.g. it has a $\bC^*$-stabilizer). Note also that the orbit of  $(C,L)$ is a minimal orbit in $X^{ss}(t)$ for $t\in(0,1)$. At $t=1$,  the closure of the orbit of $(C,L)$ contains the orbit of the semistable pair $(C_0,L)$, where $C_0:(x_0^2x_1x_2^2=0)$.
\end{example}

On the other hand, the type of  example mentioned above does not occur for $t>1$.
\begin{lemma}
For any  $t>1$ which is not critical, we have $X^s(t)=X^{ss}(t)$. Thus,   $\calM(t)$ is a geometric quotient for all non-critical $t\in (1,\frac{5}{2})$. \qed
\end{lemma}

\subsection{Variation of GIT quotients for degree $5$ pairs}\label{symquintic} One of the essential aspects of the theory of variation of GIT quotients of \cite{thaddeus,dolgachevhu} is the fact that the birational transformations that occur at a wall crossing are very explicit. Specifically, let $t$ be a critical slope (wall). We have  $X^{ss}(t\pm \epsilon)\subsetneq X^{ss}(t)$. These inclusions induce contraction maps $\calM(t\pm \epsilon)\to \calM(t)$, whose common center $Z$ is the GIT quotient of $X^{ss}(t)\setminus(X^{ss}(t+\epsilon)\cup X^{ss}(t-\epsilon))$ by $G$. Similarly, the exceptional sets $E_{\pm}$ of these birational maps are $(X^{ss}(t\pm\epsilon)\setminus X^{ss}(t\mp \epsilon))\gquot G$. Furthermore,  precise results about the local structure of these maps can be obtained by applying Luna's slice theorem. In particular, the case when the stabilizer of the minimal orbits at $t$ is $\bC^*$ is well-behaved (see \cite[Thm. 5.6]{thaddeus} and \cite[Thm. 4.2.7]{dolgachevhu}).

Ours is a very simple situation. Namely, for all critical slopes $t$ except $0$, $1$ and $\frac{5}{2}$, the center  $Z\subset \calM(t)$ is a point and the two exceptional sets $E_{\pm}$ are weighted projective spaces (possibly modulo a finite group) of complementary dimension, i.e. $\dim E_++\dim E_-=\dim \calM(t)-1$.  Furthermore, $E_{\pm}$ can be given modular interpretation. This is due to the fact that, for each critical $t$ (except $0$, $1$ and $\frac{5}{2}$), there is only  one geometric situation relevant for the change of stability as illustrated by the following example: 

\begin{example}
The change of stability at $t=\frac{5}{3}$ can be described as: the pairs $(C,L)$ with $C$ having a singularity of type $W_{12}$ become stable, and those with $L$ $5$-fold tangent to $C$ become unstable (cf. \ref{tpfor5} and \ref{lemmasimplesing}). The center $Z$ of the birational transformations at $t=\frac{5}{3}$ corresponds to the unique pair $(C,L)$ with $C$ having a singularity of type $W_{12}$ and  $L$ being $5$-fold inflectional to $C$. Moreover, $E_+$ can be interpreted as the stratum of $W_{12}$ occurring in the deformation of $N_{16}$, and $E_-$ as the locus of quintics with a flex of order $5$. 
\end{example}

The following proposition concludes our discussion.
\begin{proposition}\label{orbitg}
For each critical slope $t\in\{\frac{1}{7},\frac{1}{4},\frac{2}{5},\frac{5}{8},\frac{10}{7},
\frac{8}{5},\frac{5}{3},\frac{7}{4},\frac{13}{7},2,\frac{11}{5}\}$ there exists a unique closed orbit $O((C_t,L))\subset X^{ss}(t)\setminus X^{ss}(t\pm \epsilon)$ with stabilizer $G^0_{(C_t,L)}\cong \bC^*$. The equations and the geometry of the 
pairs $(C_t,L)$ are described in  table \ref{special5} (see also the similar list in Wall \cite{wallquintic}). \qed
\end{proposition}

\begin{table}[htb]
\begin{center}
\renewcommand{\arraystretch}{1.25}
\begin{tabular}{|c|c|c|c|c|c|}
\hline
$t$           & Equation for $C$    &Sing. at $p$&Sing. at $p'$&$\mult_{p'}(C\cap L)$&List \cite{wallquintic}\\ 
\hline\hline
$\frac{1}{7}$ & $x_0^2x_2^3+x_0x_1^4=0$   &$E_6$       &$A_7$   &$L$ is a comp. of $C$  &(H2) \\ \hline
$\frac{1}{4}$ & $x_0^2x_1x_2^2+x_1^4x_2=0$&$D_8'$       &$D_5$   &$L$ is a comp. of $C$  &(H3) \\ \hline
$\frac{2}{5}$ & $x_0x_1^3x_2+x_0^2x_2^3=0$&$E_7$       &$A_5$   &$L$ is a comp. of $C$  &(H1)\\ \hline 
$\frac{5}{8}$ & $x_0^2x_2^3+x_1^5=0$      &$E_8$      &$A_4$   &$5$ &(H5) \\ \hline
$\frac{10}{7}$& $x_0x_1x_2^3+x_1^5=0$     &$Z_{11}$   &$A_1$   &$5$   &(H14) \\ \hline
$\frac{8}{5}$ & $x_0x_1x_2^3+x_1^4x_2=0$  &$Z_{12}$    &$A_1$   &$4$ &(H10)\\ \hline
$\frac{5}{3}$ & $x_0x_2^4+x_1^5=0$        &$W_{12}$    &smooth  &$5$   &(H13)\\ \hline
$\frac{7}{4}$ & $x_0^2x_2^3+x_1^3x_2^2=0$ &double line &$A_2$   &$3$   & -    \\ \hline
$\frac{13}{7}$& $x_0x_1x_2^3+x_1^4x_2=0$  &$W_{13}$    &smooth  &$4$&   (H9) \\ \hline
$2$           & $x_0x_1x_2^3+x_1^3x_2^2=0$&double line  &$A_1$   &$3$ & -   \\ \hline
$\frac{11}{5}$& $x_0x_2^4+x_1^3x_2^2=0$   &double line   &smooth  &$3$  & -   \\ \hline
\end{tabular}
\vspace{0.2cm}
\caption{The minimal orbits  for the critical slopes $t\neq 0,1,\frac{5}{2}$}\label{special5}
\end{center}
\end{table}

For the special critical slopes $0$, $1$ and $\frac{5}{2}$, we have the following structural results. From \ref{caset0} and \ref{structured2} it follows that the natural morphisms $\calM(\epsilon)\to \calM(0)$ and $\calM(\frac{5}{2}-\epsilon)\to \calM(\frac{5}{2})$ are fibrations: they are generically a $\bP^2$-bundle and a weighted projective bundle respectively.  The remaining case $t=1$ can be understood similarly to the cases when $t\neq 1$, but there are several issues that one has to consider. For example,  one complication is that the centers of the birational transformations that occur at $t=1$ are curves (see Prop. \ref{boundary2}), and thus we obtain weighted projective bundles for the exceptional loci (cf. \cite{dolgachevhu,thaddeus}). A more serious issue is that there exist  points in $\calM(1)$ such that the corresponding stabilizers are larger than $\bC^*$ (e.g. III(2) of \ref{boundary2}). One can still apply the Luna's slice theorem to understand the local structure at those points, but the situation is slightly more complicated than in the $\bC^*$ case.

\subsection{Stability of degree $5$ pairs at $t=1$}\label{quintic1}
In constrast to the other critical slopes, the change of stability at $t=1$ is quite involved, as seen by  inspecting the statements of \S\ref{stabred5} and \S\ref{nongen5}. If we restrict to the pairs $(C,L)$ with $L$ generic the change of stability at $t=1$ becomes more conceptual: it corresponds to the division of singularities in three classes as given by corollary \ref{corquintic2} (see also remark \ref{remarkt1}). 

For us, the most important aspect about the slope $t=1$ is the fact noted in \ref{corquintic}, namely that the pairs $(C,L)$ with $C+L$ defining a sextic with simple singularities are GIT stable. Thus, we can define a moduli space of such pairs as the geometric quotient: 
\begin{equation}\label{defm}
\calM:=\{(C,L)\mid C+L \textrm{ has at worst simple singularities}\}/(\textrm{proj. equiv.})
\end{equation}
The space $\calM$ is a quasi-projective variety, compactified by the GIT quotient $\calM(1)$. The boundary components of this compactification are given by the following result.

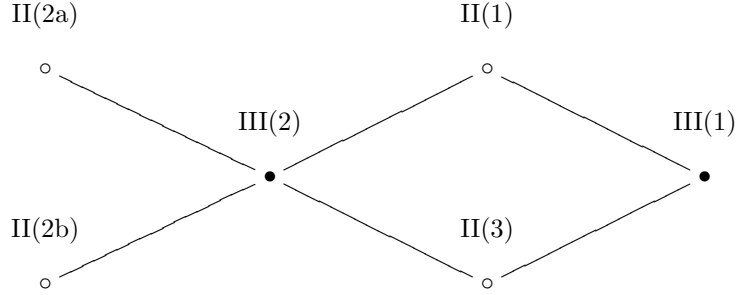
\begin{figure}[htb]
$$\xymatrix@R=.25cm{
{\textrm{II(2a)}}   &&                      &&{\textrm{II(1)}}    &&                 \\
{\circ}\ar@{-}[ddrr]&&                      &&{\circ}\ar@{-}[ddrr]&&                 \\
                    &&{\textrm{III(2)}}     &&                    &&{\textrm{III(1)}}\\
                    &&{\bullet}\ar@{-}[uurr]\ar@{-}[ddrr]&&                    &&{\bullet}        \\
{\textrm{II(2b)}}   &&                      &&{\textrm{II(3)}}    &&                 \\
{\circ}\ar@{-}[uurr]&&                      &&{\circ}\ar@{-}[uurr]&&                  \\
}
$$
\caption{Incidence diagram of the boundary components}\label{incidence}
\end{figure}

\begin{proposition}\label{boundary2}
The boundary $\calM(1)\setminus \calM$ consists of four one-dimensional components and two zero-dimensional
components  as  described below:
\begin{itemize}
\item[] Zero-Dimensional Components:
\begin{itemize}
\item[III(1)] The point corresponding to the closed orbit of the pair $(C,L)$ with equations $L: (x_0=0)$ and $C:\left(
x_2(x_0x_2-x_1^2)^2=0\right)$, and 
\item[III(2)] The point corresponding to the closed orbit of the pair $(C,L)$ with equations $L: (x_0=0)$ and $C:(
x_0x_1^2x_2^2=0)$.  
\end{itemize}
\item[] One-Dimensional Components:
\begin{itemize}
\item[II(1)]The rational curve  parameterizing  the orbits of the pairs $(C_\lambda,L)$ with equations $L: (x_0=0)$ and $C_\lambda:\left(x_2(x_0x_2-x_1^2)(x_0x_2-\lambda x_1^2)=0\right)$, where $\lambda\neq 0,1,\infty$,  
\item[II(2a)] The rational curve  parameterizing  the orbits of the pairs $(C_\lambda,L)$ with equations $L: (x_0=0)$ and 
$C_\lambda:\left(x_0 x_1x_2(x_2-x_1)(x_2-\lambda x_1)=0\right)$, where $\lambda\neq 0,1,\infty$,  
\item[II(2b)] The rational curve  parameterizing  the orbits of the pairs $(C_\lambda,L)$ with equations $L: (x_1=0)$ and 
$C_\lambda:\left(x_0^2 x_2(x_2-x_1)(x_2-\lambda x_1)=0\right)$, where $\lambda\neq 0,1,\infty$, and   
\item[II(3)] The rational curve  parameterizing  the orbits of the pairs $(C_\lambda,L)$ with equations $L: (x_1=0)$ and 
$C_\lambda:\left((\lambda x_0-(\lambda+1)x_1+ x_2)(x_0x_2-x_1^2)^2=0\right)$, where $\lambda\neq 0,1,\infty$.  
\end{itemize}
\end{itemize}
The incidence relations are given in figure \ref{incidence}. All the components except II(3) are strict GIT boundary components (i.e. parameterizes strictly semistable pairs). The  stabilizer for the type II boundary components (except II(3)) is $\bC^*$. \qed
\end{proposition}

\begin{remark} The labeling of the boundary components is in accordance with the similar list of Shah for plane sextics (\cite[Thm. 2.4]{shah}). Some of the boundary components in Shah's list do not occur in our situation (e.g. II(4) and IV), and the case II(2) splits in two subcases (depending on the relative position of the line $L$). 
\end{remark}

We note that there exists a close relationship between the boundary components of $\calM\subset \calM(1)$ and the simple elliptic and cusp singularities adjacent to $N_{16}$.   A similar situation was observed by Brieskorn \cite{brieskorn} for the case of the triangle singularities. He noticed that there exists a natural matching between the Baily-Borel compactification of a certain period domain $\calD/\Gamma$ and the simple elliptic and cusp singularities adjacent to the given triangle singularity. Furthermore, the incidence diagram of the boundary components coincides with the adjacency diagram for the corresponding singularities. This is also the case in our situation. Here we note this for the GIT compactification. By the results of section \ref{sectk3}, $\calM(1)$ can be interpreted as a Baily-Borel compactification of an appropriate $\calD/\Gamma$. Thus, the situation for $N_{16}$ is completely analogous to that for the triangle singularities.

\begin{figure}[htb]
$$\xymatrix@R=.25cm{
{\textrm{II(2a)}}             &&{\textrm{III(2)}}           &&{\textrm{II(1)}}             &&{\textrm{III(1)}}\\
{\circ}\ar@{-}[rr]&&{\bullet}\ar@{-}[rr]&&{\circ}\ar@{-}[rr]&&{\bullet}\\
{\Es}              &&{T_{2,p,q}}  &&{\Ee} &&{T_{2,3,p}}  
}
$$
\caption{The simple elliptic and cusp singularities adjacent to $N_{16}$}\label{simpleellipticcusp}
\end{figure}
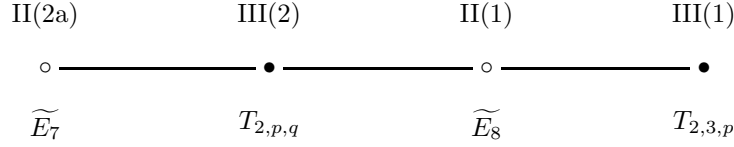

\begin{proposition}\label{strataquintic}
Assume that $(C,L)$ is a semistable pair at $t=1$ such that $C$ has at least one non-simple singularity. Assume also that $L$ is generic. Then, the image  $x\in \calM(1)\setminus \calM$ of the pair under the natural projection map $X^{ss}(1)\to \calM(1)$ satisfies: 
\begin{itemize}
\item[i)] If $C$ has a simple elliptic singularity of type $\widetilde{E_7}$ ($\widetilde{E_8}$), then $x$ belongs to the boundary component II(2a) (resp. II(1)).
\item[ii)] If $C$ has a cusp singularity of type $T_{2qr}$ with $4\le q\le r\le 6$ ($T_{23k}$ with $7\le k\le 10$), then $x$ is the point III(2) (resp. III(1)).
\item[iii)] If $C$ is non-reduced, then $x$ belongs to one of the components II(2b), II(3), III(1) or III(2).
\end{itemize}
In particular, there exists a natural matching between the simple elliptic and cusp singularities adjacent to $N_{16}$ and the boundary components of $\calM(1)$ as given in figure \ref{simpleellipticcusp}. \qed
\end{proposition}

\begin{remark}\label{remarkt1}
We close by noting the following explanation for the division of singularities given  by \ref{corquintic2}. According to theorem \ref{mainthm2}, the semistability at $t=1$ of the pair $(C,L)$ is roughly equivalent to the pair $(\bP^2,\frac{1}{2}(C+L))$ being log canonical. Let $S$ be the double cover of $\bP^2$ branched along $C+L$. From \cite[Prop. 5.20]{kollarmori}, it follows that $(\bP^2,\frac{1}{2}(C+L))$ is log canonical if and only if $S$ is log canonical. The classification of canonical and log canonical singularities is given by \cite[Ch. 4]{kollarmori}: the simple (du Val) singularities are canonical, and the simple elliptic and cusp singularities are strictly log canonical. The singularities of $S$ are in one-to-one correspondence (including the type) with the singularities of $C+L$. Since for $d=5$ and $t=1$, the converse of theorem \ref{mainthm2} also holds we obtain the characterization  of \ref{corquintic} and \ref{corquintic2} for the stability at $t=1$.
\end{remark}

\section{Moduli of pairs via K3 surfaces}\label{sectk3}
An alternative construction of the moduli space of degree $5$ pairs is obtained via the periods of $K3$ as follows. To a generic degree $5$ pair $(C,L)$ we associate the surface $S_{(C,L)}$ obtained as a double cover of $\bP^2$ along the sextic $C+L$. The minimal desingularization $\widetilde{S}_{(C,L)}$ of  $S_{(C,L)}$  is a degree two  $K3$ surface. It is easy to see that  $\widetilde{S}_{(C,L)}$ is an $M$-polarized $K3$ surface in the sense of Nikulin and Dolgachev (see \cite{mirrork3}), where $M$ is the rank $6$ hyperbolic lattice spanned by the polarization class and the exceptional divisors. It follows that the natural period map gives a birational isomorphism  between the moduli space of degree $5$ pairs and the moduli space $\calD/\Gamma$ of $M$-polarized $K3$ surfaces. This construction works in fact for all pairs $(C,L)$ such that $C+L$ has at worst simple singularities. As a consequence, we obtain:

\begin{theorem}\label{mainthm3}
The birational map which associates to a degree $5$ pair $(C,L)$ the periods of the $K3$ surface $\widetilde{S}_{(C,L)}$ extends to an isomorphism $\calP:\calM\to \calD/\Gamma$. 
\end{theorem}

We recall that $\calM$ denotes the moduli space of degree $5$ pairs $(C,L)$ satisfying the condition that the sextic $C+L$ has at worst simple singularities. The space $\calM$ is constructed as a geometric quotient, and it is compactified by the GIT quotient  $\calM(1)$. As mentioned above, $\calD/\Gamma$ is the moduli of $M$-polarized $K3$ surface. The precise definitions are given in \S\ref{quinticperiod} below. The proof of the theorem then follows from  standard results on $K3$ surfaces and the analysis of the geometric meaning of the notion of $M$-polarization. The details are given in \S\ref{sectproofs}. We note additionally that  it follows  automatically from work of Looijenga \cite{looijengacompact} that the isomorphism of \ref{mainthm3} extends  to the boundary.
 
\begin{theorem}\label{mainthm4}
The period map $\calP:\calM\to \calD/\Gamma$  extends to an isomorphism of projective varieties
 $\overline{\calP}:\calM (1)\xrightarrow{\cong} (\calD/\Gamma)^*$, where $(\calD/\Gamma)^*$ denotes the Baily-Borel compactification of $\calD/\Gamma$.
\end{theorem}

The Baily-Borel compactification and the matching of the boundary components given by the above theorem are discussed in \S\ref{sectbailyborel}.

\subsection{Standard notations and facts about lattices and $K3$ surfaces}
By a {\it lattice} we understand a free $\bZ$-module $L$ together with a symmetric bilinear form, which we denote  by $x.y$ for $x,y\in L$. The basic invariant of a lattice is its {\it signature}. In particular, if the signature is $(1,n)$ we call the lattice {\it hyperbolic}. A lattice is {\it even} if $x.x\equiv 0 \mod 2$ for every $x\in L$.  The direct sum of lattices is always assumed orthogonal. For a lattice embedding $M\subseteq L$, $M^{\perp}_L$ denotes the orthogonal complement of $M$ in $L$. 

The following standard lattices are used frequently: {\it the root lattices} $A_n$ (for $n\ge1$), $D_m$ (for $m\ge 4$), and $E_r$ (for $r=6,7,8$), {\it the hyperbolic plane} $U$, {\it the parabolic lattices} (negative semidefinite) $\widetilde{A}_n$, $\widetilde{D}_m$; and $\widetilde{E}_r$, and the hyperbolic lattices $T_{p,q,r}$ for $\frac{1}{p}+\frac{1}{q}+\frac{1}{r}<1$.
Given $L$ a lattice, $L(n)$ denotes the lattice with the same underlying $\bZ$-module as $L$ but with the bilinear form multiplied by $n$.

\begin{notation}\label{notationlattice}
Let $L$ be an even lattice, we define: 
\begin{itemize}
\item[-] $\Div(x)$:  the positive integer $d$ such that  $x.L\cong d\bZ$;
\item[-] $L^*:=\{y\in L\otimes \bQ\mid x.y\in \bZ \textrm{ for all } x\in L\} $ the {\it dual lattice};
\item[-] $A_L=L^*/L$: the {\it discriminant group} endowed with the induced finite quadratic form $q_L$;
\item[-] $O(L)$: the group of isometries of $L$;
\item[-] $O(q_L)$: the automorphisms of $A_L$ that preserve the quadratic form $q_L$;
\item[-] $O_{-}(L)$: the group of isometries of $L$ of spinor norm $1$;  
\item[-] $\widetilde{O}(L)$: the group of isometries of $L$ that induces the identity on  $A_L$;
\item[-] $O^{*}(L):=O_{-}(L)\cap \widetilde{O}(L)$;
\item[-] $\Delta(L)$: the set of roots of $L$, where $\delta\in L$ is a {\it root} if $\delta^2=-2$;
\item[-] $\widetilde{\Delta}(L)$: the set of generalized roots of $L$, where $\delta\in L$ is a {\it generalized root} if $\delta^2=-k$ and $\Div(\delta)=\frac{k}{2}$ for some even positive integer $k$;
\item[-] $W(L)$ (and $\widetilde{W}(L)$): {\it the (generalized) Weyl group}, i.e. the group of isometries generated by reflections $s_\delta$ in (generalized) roots $\delta$, where 
$$s_\delta(x)=x-2\frac{x.\delta}{\delta^2}\delta.$$
\end{itemize}
\end{notation}

\begin{definition}
Given two lattices $L$ and $L'$ and a lattice embedding $L\hookrightarrow L'$, we call it a {\it primitive embedding} iff $L'/L$ is a free $\bZ$-module. Equivalently, the embedding is primitive iff $\Sat(L):=\{y\in L'\mid n y\in L \textrm{ for some positive integer }n\}$ coincides with $L$.
\end{definition}

For a surface $S$, the intersection form gives a natural lattice structure on the torsion-free part of  $H^2(S,\bZ)$ and on the Neron-Severi group $\textrm{NS}(S)$. For a $K3$ surface, we have $H^1(S,\calO_S)=0$, and we identify $\Pic(S)\cong \textrm{NS}(S)$.  Both $H^2(S,\bZ)$ and $\Pic(S)$ are torsion free.  The natural map $\Pic(S)\xrightarrow{c_1}H^2(S,\bZ)$ is a primitive lattice embedding. By Hodge index theorem, $\Pic(S)$ is a hyperbolic lattice.

\begin{notation}
If $S$ is a $K3$ surface, we use $O(S)$, $W(S)$, $\Delta(S)$, etc. to denote the corresponding  objects associated to the lattice $\Pic(S)$. We also use:
\begin{itemize}
\item[-] $\Delta^+(S)$ for the set of effective $(-2)$ divisor classes in $\Pic(S)$;
\item[-] $V^+(S)\subset S_\bR$ the K\"ahler cone;
\item[-] $C^+(S)=\{x\in \Pic(S)\cap V^+(S) \mid x.\delta>0 \textrm{ for all }\delta\in \Delta^+(S)\}$ the ample cone.
\end{itemize} 
We note that, on a $K3$ surface, any $x\in \Pic(S)$ with $x.x\ge -2$ has the property that either $x$ or $-x$ is effective. In particular, $\Delta(S)=\Delta^+(S)\sqcup (-  \Delta^+(S))$.
\end{notation}

\begin{notation}\label{notationm0}
Unless specified otherwise, the symbols $\Lambda$, $M$ and $T$ will denote throughout the chapter the lattices $E_8^{\oplus 2}\oplus U^{\oplus 3}$, $D_4\oplus U(2)$ and $D_4\oplus E_8\oplus U\oplus U(2)$ respectively. The lattice $\Lambda$ is the unique even, unimodular lattice of signature $(3,19)$; $\Lambda$ is isometric to $H^2(S,\bZ)$ for any $K3$ surface $S$. The lattices $M$ and $T$ have signature $(1,5)$ and respectively $(2,14)$   and they can be embedded in $\Lambda$ such that they are mutually orthogonal.
\end{notation}

\begin{definition}
{\it A polarization} for a $K3$ surface is the class of a nef and big divisor $H$. The {\it degree of the polarization (and of  the surface)}  is $H^2$. 
\end{definition}

We use the following standard results on $K3$ surfaces: the global Torelli theorem, the surjectivity of the period for algebraic $K3$ surfaces, and the results on linear systems on a $K3$ surface. We recall the following theorem of Mayer (\cite[Thm 27]{friedmansurf}).

\begin{theorem}[Mayer's Theorem]\label{mayerthm}
Let $H$ be a nef and big divisor on the $K3$ surface $S$. Then $|H|$ has a base point iff $|H|$ has a fixed curve iff $H=kE+R$ (linearly equivalent), where $E$ is a smooth elliptic curve, $R$ is a smooth rational curve, $R.E=1$,  and $k\ge 2$.
\end{theorem}

We note also the following easy converse (\cite[Prop. 1, pg. 35]{dmorrison}).
\begin{proposition}\label{conversemayer}
Let $D$ be a big and nef divisor on a $K3$ surface. Assume that there exists a divisor $E$ such that $D.E=1$ and $E^2=0$, then the linear system $|D|$ has a fixed component.
\end{proposition}

\subsection{The $K3$ surface associated to a degree $5$ pair and the period map}\label{quinticperiod} 
Let $(C,L)$ be a degree $5$ pair such that $C+L$ is a sextic with at worst simple singularities, and $\widetilde{S}_{(C,L)}$ the associated degree two $K3$ surface.  We are concerned  here with establishing the basic properties of this type of surfaces, and finding a moduli space for them.   We proceed in three basic steps. First, by considering the generic case, we see that $\widetilde{S}_{(C,L)}$ is naturally $M$-polarized for a certain lattice $M$. Then, we establish some basic arithmetic properties about $M$ and its orthogonal complement $T:=M^\perp_\Lambda$ and determine the period domain $\calD/\Gamma$.  Finally, we show that the construction can be extended to the non-generic case. 

\subsubsection{The Picard lattice of $\widetilde{S}_{(C,L)}$}\label{picardlattice} 
In this subsection, we assume that $(C,L)$ is a generic degree $5$ pair. In particular, $C$ is smooth and $L$ is transversal.  The $K3$ surface $\scl$ is obtained as the desingularization of $S_{(C,L)}$, the double cover of $\bP^2$ along $C+L$. Let $\pi:\scl\to \bP^2$ be the natural projection, and $h=\pi^* l\in \Pic(\scl)$ the pullback of the class of a line from $\bP^2$. Since  the surface $S_{(C,L)}$ has five ordinary double points coming from the intersection $C\cap L$, it follows that  $\Pic(\scl)$ contains five additional classes $e_1,\dots,e_5$ corresponding to the exceptional divisors of $\scl\to S_{(C,L)}$. By construction we have that $h$ is a degree $2$ polarization for $\scl$ and the intersection numbers $h.e_i=0$ and $e_i.e_j= -2\delta_{ij}$. Finally, since $L$ is in the branch locus of $\pi$,  the surface $\scl$ contains  another rational curve, namely the inverse image  $L'$ of the line $L$. In conclusion, we obtain:

\begin{lemma}\label{lemmaemb} 
Let $\scl$ be the $K3$ associated to a generic degree $5$ pair $(C,L)$. Then $\scl$ contains $6$ irreducible $(-2)$-curves, whose classes $l',e_1,\dots,e_5$ satisfy $l'.e_i=1$ and $e_i.e_j=0$ for $i\neq j$. Additionally, the structural morphism $\pi:\scl\to \bP^2$ is given by the class $h:=2l'+e_1+\dots+e_5$. \qed
\end{lemma} 

As a consequence of the previous lemma, it follows that the Picard lattice of $\scl$ contains the sublattice spanned by the classes $l',e_1,\dots, e_5$. We will see that generically  $\Pic(\scl)$ coincides with this sublattice. Conversely, we will prove that this condition characterizes the surfaces $\scl$. To proceed, we need to fix the following notations (compatible with \ref{notationm0}).

\begin{notation}\label{notationm}
We denote by  $M$ the abstract rank $6$ lattice spanned by $l',e_1,\dots,e_5$ with the intersection form given by: $l'^2=e_i^2=-2$,  $l'.e_i=1$ and $e_i.e_j=0$ for $i\neq j$. The  basis $\{l',e_1,\dots,e_5\}$ of $M$ is assumed fixed. We denote
$$h:=2l'+e_1+\dots+e_5$$ 
and $f_i:=h-e_i$  for  $i=1,\dots,5$. In particular, $h^2=2$, $h.l'=1$, and $h. e_i=0$.
\end{notation}

\begin{proposition}\label{primitive}
Let $(C,L)$  be a generic degree $5$ pair and $j:M\to \Pic(\scl)$ the lattice embedding given by lemma \ref{lemmaemb}.  Then  $j$ is a primitive lattice embedding. 
\end{proposition}
\begin{proof}
Assume that $j$ is not primitive. Then the embedding $j$ factors as 
$$M\subsetneq \Sat(M)\hookrightarrow \Pic(\scl),$$
where $\Sat(M)$ is the saturation of the lattice $M$ in $\Pic(\scl)$. Thus, $\Sat(M)$ is a nontrivial overlattice of $M$ and as such is classified by a nontrivial isotropic subgroup of $A_M$ (\cite[Sect. 4]{nikulin}). The discriminant group $A_M= M^*/M\cong (\bZ/2\bZ)^4$ is generated by $\{f_i^*\mid i=1,\dots 5\}$, where $f_i^*$ is the class of $f_i/\Div(f_i)$ (N.B. $\Div(f_i)=2$ and the notations are those of \ref{notationlattice} and \ref{notationm}). It is easy to see that the only nontrivial isotropic elements of $A_M$ (for the induced quadric form) are precisely the elements $f_i^*$. It follows that $M\neq \Sat(M)$ is equivalent to some $f_i$ being $2$-divisible in $\Pic(\scl)$. We have $f_i=2f_i'$ for some $f_i'\in \Pic(\scl)$ and then
$$h=f_i+e_i=2f_i'+e_i.$$ 
By \ref{conversemayer} it follows that $h$ is not base point free. This is a contradiction to the fact that  $h$ defines the morphism $\pi:\scl\to \bP^2$. \end{proof}
\begin{remark}\label{strongprimitive}
More generally, let $S$ be any $K3$ surface. The same argument as above gives the following statement. {\it If $j:M\to \Pic(S)$ is a lattice embedding with $j(h)\in \Pic(S)$ a base point free polarization, then $j$ is primitive.} 
\end{remark}

\begin{proposition}\label{picequal}
Assume that $S$ is a $K3$ surface such that  $\Pic(S)$ is isometric to the lattice $M$. Then $S$ is the double cover of $\bP^2$ branched over a reducible sextic $C+L$. Moreover, $C$ is a smooth quintic and $L$ is a line intersecting $C$ transversely.   
\end{proposition}
\begin{proof}
By assumption there exist $h,l',e_1,\dots,e_5\in \Pic(S)$ satisfying the numerical conditions from \ref{notationm}. There is no loss of generality to assume that $h$ is nef (if not, this can be achieved by acting by $\pm W(S)$). By acting with the reflections $s_{e_i}$ (i.e. change $e_i$ to $-e_i$) we can further assume that the classes $e_i$ are effective. 

The class $h$ defines a degree two polarization for $S$. We claim that $h$ is base point free. If not, by  Mayer's theorem, we get that $h=2d+r$ with $d^2=0$ and $d.h=1$ for some $d,r\in \Pic(S)\cong M$. This gives a contradiction. Namely, we write:
$$d=a l'+b_1 e_1\dots +b_5 e_5$$ 
for some integers $a,b_1,\dots,b_5$. The condition $d.h=1$ gives $a=1$, and then
$$d^2=-2+2(b_1+\dots+b_5)-2(b_1^2+\dots+b_5^2)\equiv 2 \mod 4,$$ 
contradicting the assumption that $d^2=0$. In conclusion, the linear system defined by $h$ gives a degree two map $\pi:S\to \bP^2$ branched along a sextic $B$. 

Since $\langle h\rangle^\perp_M\cap \Delta(M)=\{\pm e_1,\dots,\pm e_5\}$, it follows that $e_1,\dots,e_5$ are classes of irreducible rational curves $E_1,\dots,E_5$. These curves are contracted by $\pi$ to five ordinary double points for the sextic $B$. Let $L'$ be the curve corresponding to $l'$ and $L:=\pi(L')$. Since $l'.h=1$, the projection formula gives that $L$ is a line.  Moreover, since $l'.e_i=1$, $L$ has to pass through all $5$ singular points of the branch curve $B$. By Bezout, the only possibility is that $L$ is a component of $B$. 
\end{proof}

From the previous two propositions, it follows that generically $\Pic(\scl)\cong M$. Here, generically should be understood in the sense of moduli, i.e. it holds on the complement of the union of a countable number of proper subvarieties.
 
\begin{corollary}\label{corpicequal}
Assume that $(C,L)$ is a sufficiently general degree $5$ pair, then  there exists a lattice isometry $M\cong \Pic(\scl)$. \qed
\end{corollary}

Without the genericity assumption we obtain the following statement. 
\begin{lemma}\label{coremb}
Let $(C,L)$ be a degree $5$ pair such that $C$ has at worst simple singularities and $L$ intersects $C$ transversely.  Then there exists a primitive embedding $j:M\hookrightarrow \Pic(\scl)$ such that $j(h)$ is a base point free degree $2$ polarization.
\end{lemma}
\begin{proof}
Since we assume transversal intersection, the construction of lemma \ref{lemmaemb} applies ad litteram. In particular, the class $j(h)$ defines the morphism $\pi:\scl\to \bP^2$, and thus it is base point free. The embedding $j$ primitive by \ref{strongprimitive}.
\end{proof}

\subsubsection{$M$-polarized $K3$ surfaces and the period map}\label{mpolarized} 
The statements of \ref{corpicequal} and \ref{coremb} say that the $K3$ surfaces $\scl$ associated to a degree $5$ pair $(C,L)$ are characterized by the fact that they are $M$-polarized $K3$ surfaces. 

\begin{definition}\label{defmpol}
Let $M$ be the lattice defined in \ref{notationm}. An {\it $M$-polarized $K3$ surface} is a pair $(S,j)$ such that $j:M\hookrightarrow \Pic(S)$ is a primitive lattice embedding. The embedding $j$ is called the {\it $M$-polarization} of $S$. If the polarization is understood, we simply say $S$ is an {\it $M$-polarized $K3$ surface}.
\end{definition}

The lattice $M$ admits a unique primitive embedding into the $K3$ lattice $\Lambda$.

\begin{lemma}\label{basicm}
Let $M$ be  as in \ref{notationm}. Then $M$ is isometric to the lattice $D_4\oplus U(2)$ and admits a unique primitive embedding $M\hookrightarrow \Lambda$ into the $K3$ lattice $\Lambda$. The orthogonal complement $T:=M^\perp_\Lambda$ with respect to this embedding is isometric to  $D_4\oplus E_8\oplus U\oplus U(2)$.
\end{lemma}
\begin{proof}
The isometry $M\cong D_4\oplus U(2)$ is given by taking  $\{ l',e_1,e_2,e_3,f_4,f_5\}$ as a basis for $M$. The existence of the primitive embedding $M\hookrightarrow \Lambda$ follows from \cite[Thm. 1.14.4]{nikulin}. The uniqueness is essentially equivalent to the fact that $M$ and its orthogonal complement $T$ are uniquely determined by their invariants (i.e. signature and discriminant group). Since $M$ and $T$ are indefinite lattices, the latter statement follows from  \cite[Thm. 1.10.1]{nikulin}. \end{proof}

It is a standard fact  (see \cite{mirrork3}) that the moduli space of $M$-polarized $K3$ surfaces is a quotient $\calD/\Gamma$ for a certain bounded symmetric domain $\calD$ of type IV and a certain arithmetic group $\Gamma$. Namely, the condition of the $M$-polarization determines a  tower of primitive embeddings
$$M\hookrightarrow \Pic(S) \hookrightarrow H^2(S,\bZ)\cong \Lambda.$$
This implies that the period of an $M$-polarized $K3$ surface belongs to the following subdomain of the period domain of $K3$ surfaces: 
$$\calD=\{\omega\in \bP(\Lambda\otimes \bC)\mid \omega.\omega=0,\ \omega.\bar{\omega}>0, \textrm{ and } \omega\perp M\}_0$$
(the (unique) embedding $M\hookrightarrow \Lambda$ is assumed fixed). Conversely, since $\Pic(S)\cong H^2(S,\bZ)\cap H^{1,1}(S)$, every point of $\calD$ corresponds to an $M$-polarized $K3$ surface. Let $T=M^{\perp}_\Lambda$ be the orthogonal complement of $M$ (called {\it the transcendental lattice}). It is convenient to identify $\calD$ to the domain
\begin{equation}\label{defcald}
\{\omega\in \bP(T\otimes \bC)\mid \omega.\omega=0,\ \omega.\bar{\omega}>0\}_0.
\end{equation}
In particular, note the natural action of  the groups  $O^*(T)$ and $O_{-}(T)$ on $\calD$.

To specify the moduli of $M$-polarized $K3$ surfaces it suffices to determine the arithmetic group $\Gamma$.  In the standard situation considered in \cite{mirrork3} one requires that the $M$-polarization is pointwise fixed by group $\Gamma$, and thus takes $\Gamma$ to be $O^*(T)$. In our geometric context we have to require less, namely $\Gamma$ should fix (not necessarily pointwise)  $M$ and the element $h\in M$. Thus, we make the following definition for $\calD$ and $\Gamma$. The reason for this choice is explained in proposition \ref{propmonodromy} below. 

\begin{notation}\label{defgamma}
Fix the primitive embedding $M\hookrightarrow \Lambda$ and let $T:=M^\perp_\Lambda$. We define $\calD$ to be the domain given by (\ref{defcald}), and $\Gamma:=O_{-}(T)$. 
\end{notation}

To explain the choice of $\Gamma$ and to see that indeed $\calD/\Gamma$ is a moduli space for degree $5$ pairs, we note the following  properties of the lattices $M$ and $T$.

\begin{lemma}\label{permutation}
Let $\{l',e_1,\dots,e_5\}$ be the standard basis of $M$ and regard $\Sigma_5$ as the subgroup of $O(M)$ which permutes the basis elements $\{e_1,\dots,e_5\}$. Then, the composition $\Sigma_5\to O(M)\to O(q_M)$ induces an isomorphism 
$\Sigma_5\to O(q_M)$. A similar statement holds for $T$. \qed
\end{lemma}

The lemma establishes the following relation between the arithmetic groups that occur in the construction of the moduli space for $M$-polarized $K3$ surfaces.
\begin{corollary}
Let $T$ be as above. Then $O_{-}(T)=O^*(T)\rtimes \Sigma_5$.
\end{corollary}

To understand the geometric meaning of the previous corollary, one has to investigate the relation  of these groups to the (generalized) Weyl group (see \ref{notationlattice}). 

\begin{proposition}\label{propmonodromy}
Let $T$ be the lattice of \ref{notationm0}. Then 
\begin{itemize}
\item[i)] $O^*(T)=W(T)$
\item[ii)]  $O_{-}(T)=\widetilde{W}(T)$
\end{itemize}
Geometrically, $O^*(T)$ is the local monodromy group of the singularity $N_{16}$, $O_{-}(T)$ is the monodromy group for the degree $5$ pairs $(C,L)$, and $\Sigma_5\cong O_{-}(T)/O^*(T)$ is the monodromy at infinity (i.e. it acts on the  intersection $C\cap L$).
\end{proposition}
\begin{proof}
It is easy to verify that the lattice $T$ is isometric to the Milnor lattice of the singularity $N_{16}$. For example any singularity of class $N_{16}$ has a $\mu$-constant deformation to the special case $x^5+y^5+z^2$. For this case, the Milnor lattice can be computed by the Thom-Sebastiani theorem (see \cite{agv2}). By a theorem of Ebeling (\cite[Thm. 5.5]{ebeling}) it follows that the local monodromy group coincides with $O^*(T)$ (see also \cite{pinkhammonodromy}).  The monodromy group  is generated by Picard-Lefschetz transformations (reflections) in the vanishing cycles. Thus, $O^*(T)=W(T)$. 

For the second part, we note $T\cong M\oplus E_8\oplus U$. The residual $\Sigma_5\cong O_{-}(T)/O^*(T)$ corresponds to the permutation of the basis elements $\{e_1,\dots,e_5\}$ in the $M$ summand. We note that the reflection in the generalized root $\delta=e_i-e_j$ transposes $e_i$ and $e_j$ and leaves the other basis elements invariants (including $l'$). Since $\widetilde{W}(T)\subseteq O_{-}(T)$ always holds, we conclude $O_{-}(T)=\widetilde{W}(T)$.

Let $U\subset X\cong \bP^N\times \bP^2$ be the open subset  parameterizing  degree $5$ pairs $(C,L)$ with $C$ smooth and $L$ transversal. Fix a base point $x\in U$ corresponding to a pair $(C,L)$. Then there exists a natural monodromy action of  $\pi_1(U,x)$ on $H^2(\scl,\bZ)$.  For obvious geometric reasons we have that the monodromy group $\Gamma$ satisfies:
$$\Gamma:=\textrm{Im}\left(\pi_1(U,x)\to \Aut(H^2(\scl,\bZ))\right)\subseteq O_{-}(T)$$
The argument above shows that, in fact, we have equality. Namely, the reflections in vanishing cycles corresponding to the degenerations of $C$  generate $O^*(T)\subseteq \Gamma$.  On the other hand, it is well known that $\Gamma$ acts as $\Sigma_5$ on the five points  of the intersection $C\cap L$ (e.g. \cite{harrismonodromy}). By construction, a permutation of the intersection points gives a permutation of the basis elements  $\{e_1,\dots,e_5\}$ of $M\hookrightarrow \Pic(\scl)$. Thus, indeed $\Gamma=O_{-}(T)$.
\end{proof}

We can now conclude that the moduli space of degree $5$ pairs $\calM$ constructed by GIT is birationally equivalent to the moduli of $M$-polarized $K3$ surfaces $\calD/\Gamma$.
\begin{corollary}\label{corbirational}
The period map that associates to a generic degree $5$ pair $(C,L)$ the periods of the $K3$ surface $\scl$ defines a birational map $\calP:\calM\dashrightarrow \calD/\Gamma$. 
\end{corollary}
\begin{proof}
Let $U$ be the open subset of $\calM$  parameterizing  the degree $5$ pairs $(C,L)$ with $C$ smooth and transversal intersection $C\cap L$.  Let $\widetilde{U}$ be the $\Sigma_5$-cover of $U$ that parameterizes triples $(C,L,\sigma)$, where $\sigma:\{1,\dots,5\} \to C\cap L$ is a labeling of the five points  of intersection of $C$ and $L$. By the results of \S\ref{picardlattice} (esp. \ref{coremb}), the surface $\scl$ carries a natural $M$-polarization $j:M\hookrightarrow \Pic(\scl)$. Therefore, by \cite{mirrork3}, there is a well defined map  $\widetilde{\calP}:\widetilde{U} \to \calD/O^*(T)$
sending  $(C,L,\sigma)$ to the periods of $(S_{(C,L)},j)$. The two main theorems for $K3$ surfaces (Torelli and surjectivity of the period map)  together with proposition \ref{picequal} give that $\widetilde{\calP}$ is a birational morphism. The map $\widetilde{\calP}$ is clearly $\Sigma_5$-equivariant (see \ref{propmonodromy}). Thus, it descends to a birational morphism $\calP:U \to \calD/\Gamma$ as needed.
\end{proof}

\subsubsection{$M$-Polarizations for non-generic intersections}\label{nongeneric}
In order to conclude that we have a period map $\calP:\calM\to \calD/\Gamma$ sending a degree $5$ pair $(C,L)$ to the periods of $\scl$, we need to address two points. Namely, we have to extend the construction of $M$ polarization to the case of non-transversal  intersection $C\cap L$, and we have to show that this construction fits in families.

The main problem in constructing an $M$-polarization in the case of non-transversal intersections is the fact that the Picard lattice $\Pic(\scl)$ acquires additional $(-2)$ classes. A priori it is not clear  which of those classes should be chosen as $l',e_1,\dots,e_5$ (see \ref{notationm}). For example, the reader is encouraged to consider the case when $L$ becomes simply tangent to $C$. It turns out that the right solution is to rigidify (or normalize) the notion of $M$-polarization.

\begin{definition}\label{goodmark}
Let $S$ be $K3$ surface, and $j:M\hookrightarrow \Pic(S)$ a primitive lattice embedding. We say that  $j$ is {\it a normalized $M$-polarization} if it satisfies the following two conditions:
\begin{itemize}
\item[(1)] $j(h)\in \Pic(S)$ is the class of a nef divisor;
\item[(2)] for $i=1,\dots,5$: if $\delta\in \Delta^+(S)\cap \langle j(h),\dots,j(e_{i-1})\rangle^\perp_{\Pic(S)}$ 
then $\delta.j(e_i)\le 0$.  
\end{itemize} 
\end{definition}

We now establish the main result of this section, the existence of an $M$-polarization in the non-generic case. The main idea is that we always have $5$ intersection points, but some of them are ``infinitely near points''.

\begin{proposition}\label{markinggen}
Let $(C,L)$ be a degree $5$ pair such that $C+L$ has at worst simple singularities and $\scl$ the associated $K3$ surface. Then $\scl$ carries a normalized $M$-polarization $j:M\hookrightarrow \Pic(\scl)$. 
\end{proposition}
\begin{proof}
Let $S_{(C,L)}$ be the double cover of $\bP^2$ branched along $C+L$. The surface $\scl$ is obtained from $S_{(C,L)}$ by taking a minimal resolution of the singularities. In the case of double covers, there exists a standard procedure of obtaining the desingularization (e.g. \cite[\S III.7, pg. 107]{barth}). Namely, there exists a commutative diagram:
\begin{equation}\label{resolution}
\begin{CD}
\scl@>\tau>>S_{(C,L)}\\
@V\pi'VV @VV\bar{\pi}V\\
S'@>\epsilon>>\bP^2
\end{CD}
\end{equation}
such that $\epsilon$ is a composition of blow-ups and $\pi'$ is a double cover with smooth  branch locus $B'\subset S'$.\
The surface $S'$ is obtained by an inductive process. Start with $S_0=\bP^2$ and $B_0=C+L$. Blow-up a singular point of $B_0$, and let $\epsilon_1:S_1\to S_0$ be the resulting surface. The new branch divisor $B_1$ is the strict transform of $B_0$ together with the exceptional divisor of $\epsilon_1$ reduced mod $2$. The process is repeated until the resulting divisor $B_N$ is smooth.  Let $S'=S_N$, $B'=B_N$ and $\epsilon=\epsilon_1\dots \epsilon_{N}$.  The double cover of $S'$ branched along $B'$ is  a minimal resolution of $S_{(C,L)}$. 

Since the surface $\scl$ does not depend on the order of the blow-ups, we choose to make  the first $5$ blow-ups in points belonging to (the strict transform of) the line $L$.  By abuse of notation, we denote by $L$ the strict transform of the line $L$ on all the surfaces $S_i$ (including $S'=S_N$). By construction, $L$ belongs to the branch divisor $B'$. The inverse image $L'$ of $L$ is a smooth rational curve on the $K3$ surface $\scl$. Thus, $(L')^2=-2$, which then gives $(L)^2=-4$. Since $L^2=1$ on $\bP^2$ and $L$ is a smooth curve, it follows that $L$ is affected by exactly $5$ blow-ups in the desingularization process. We choose these blow-ups as the first five steps of the sequence of blow-ups $S'\to\dots \to S_1\to \bP^2$ (this is possible). Let $p_i\in L$ be centers of these blow-ups and $E_i$ the exceptional divisors (N.B. $p_i\in S_{i-1}$ and $E_i$ is a divisor on $S_i$ for $i=1,\dots,5$). We then define the following divisors:
\begin{equation}\label{defjdelta}
D_i=\pi'^*\epsilon_{N}^*\dots \epsilon_{i+1}^*(E_i)
\end{equation}
for $i=1,\dots,5$.

The procedure described above produces $6$ divisors $L',D_1,\dots,D_5$ on the surface $\scl$ such that the polarization class of $\scl$ is
$$H:=(\epsilon \circ \pi')^*L=2L'+D_1+\dots+D_5.$$
We immediately see that $L',D_1,\dots,D_5$ satisfy the numerical conditions: $(L')^2=-2$, $D_i.D_j=-2\delta_{ij}$, and $L'.D_i=1$. For example, since $E_i$ is an exceptional divisor we have $E_i^2=-1$ on $S_i$. It follows that the pullback of $E_i$ on $S'$ will also have self-intersection $-1$. Thus, on the double cover $\scl\to S'$ we get $D_i^2=-2$. In conclusion, the linear map $j:M\to \Pic(S_{(C,L)})$ defined by sending $e_i\in M$ to be the class of $D_i$ given in (\ref{defjdelta}) and $l'\in M$ to the class of $L'$  defines a lattice embedding. Since $j(h)$ corresponds to the polarization given by $H$, we obtain that $j$ is a primitive lattice embedding (cf. \ref{strongprimitive}).

The final aspect is to note that the embedding $j$ is normalized in the sense of the definition \ref{goodmark} (i.e. satisfies the second requirement of the definition). This follows from the observation that the divisor $D_i$ on $\scl$ is  the fundamental cycle associated to the simple singularity of the curve $B_{i-1}$ in the point $p_i$.
\end{proof}

\begin{remark}\label{remarkgeom}
We note that the $M$-polarization $j:M\hookrightarrow \Pic(\scl)$ constructed in the previous proposition satisfies the following geometric properties:
\begin{itemize}
\item[i)] $j(h)$ is the class of the base point free polarization $H$;
\item[ii)] $j(l')$ is the class of the irreducible rational curve $L'$;
\item[iii)] $j(e_1), \dots, j(e_5)$ are the classes of the effective divisors $D_i$.
\end{itemize}
\end{remark}

We close by noting that the constructed $M$-polarizations for the surfaces $S_{(C,L)}$ fit well in families. Let $(\calC,\calL)\subset \bP^2_U\to U$ be the universal family of degree $5$ pairs satisfying the stability condition that $C+L$ has simple singularities. By taking a double cover, we obtain a family of surfaces $\calS\to U$ with only du Val singularities and a flat family of rational curves $\calL'\to U$.  After a finite base change, we can further assume that we have $5$ sections $\sigma_i$ corresponding to the $5$ points of intersection. By applying  Brieskorn's simultaneous resolution to the family $\calS$, we obtain a family of $K3$ surfaces $\widetilde{\calS}\to U$ (after a further finite base change of $U$). We note that the $M$-polarization in \ref{markinggen} is obtained by taking as the first $5$ steps of the simultaneous resolution process the blow-up of $\calS$ along the sections $\sigma_i$. In other words, the $M$-polarizations can be fitted together due to the fact that we can do the blow-up process of \ref{markinggen} in families (see the discussion of the simultaneous resolution from \cite[pg. 128--135]{kollarmori}).  We get a family $\widetilde{\calS}\to U$ of $M$-polarized $K3$ surfaces, which gives a period map $\widetilde{\calP}:U\to \calD/O^*(T)$ (see \cite{mirrork3}). This descends to a period map $\calP:\calM\to \calD/\Gamma$ by noting that the construction of $M$-polarization depends only on the choice of a labeling of the intersection points $C\cap L$ (i.e. a surjective map $\sigma:\{1,\dots,5\}\to C\cap L$ respecting the intersection multiplicities). Passing from $\calD/O^*(T)$ to $\calD/\Gamma$ amounts to forgetting the labeling of the intersection (see the proof of \ref{corbirational}).

\subsection{The proofs of the theorems \ref{mainthm3} and \ref{mainthm4}}\label{sectproofs}
In \S\ref{quinticperiod}, we have constructed a period map $\calP:\calM\to \calD/\Gamma$ by sending a degree $5$ pair $(C,L)$ to the periods of the $M$-polarized $K3$ surface $(\scl,j)$.  Since $\calM$ and $\calD/\Gamma$ are normal quasi-projective varieties and the period map is algebraic, to prove theorem \ref{mainthm3} it suffices to prove that $\calP$ is bijective. This amounts to showing that the data of $M$-polarization is rigid enough to recover uniquely the pair $(C,L)$ and, secondly, that any $M$-polarized $K3$ surface is of type $\scl$ for some degree $5$ pair $(C,L)$. Once theorem \ref{mainthm3} is established, we obtain immediately the stronger statement \ref{mainthm4} by applying some general results of Looijenga \cite{looijengacompact} (see \S\ref{proofmainthm4}).

\subsubsection{The surjectivity of the period map}\label{secsurjectivity} The surjectivity of the period map $\calP$ follows from the surjectivity of the period map for $K3$ surfaces and  proposition  \ref{surjectivity} below. To prove   \ref{surjectivity}, we need a series of technical lemmas. The first of those says that an $M$-polarized $K3$ surface is  a double cover of $\bP^2$. 

\begin{lemma}\label{pseudoample}
Let $(S,j)$ be an $M$-polarized $K3$ surface such that $j(h)$ is nef and $j(e_1),\dots,j(e_5)$ are classes of effective divisors. Then, the complete linear system defined by $j(h)$ is base-point free.  
\end{lemma}
\begin{proof}
By assumption, we can represent the classes $j(l'), j(\delta_1),\dots,j(\delta_5)$ by effective divisors  $L',D_1,\dots, D_5$. 
Let $H=2L'+D_1+\dots+D_5$. Assuming that $|H|$ is not base point free, we get $H\equiv 2E+R$ for some smooth elliptic curve $E$ with $H.E=1$  (cf. \ref{mayerthm}). In particular, note that $E$ is nef. Since $L',D_1,\dots, D_5$ are effective, it follows from $H.E=1$ that we can assume $E.L'=\dots=E.D_4=0$ and $E.D_5=1$. Let $F_5=H-E_5$ and $D=2E-F_5$ with classes $j(f_5)$ and $d$ respectively.  We then have $D.L'=\dots=D.D_5=0$ and $D^2=0$, which gives  $d\in M^{\perp}_{\Pic(S)}$ and $d^2=0$. Since both $M$ and $\Pic(S)$ are hyperbolic, the lattice $M^{\perp}_{\Pic(S)}$ is negative definite. In conclusion, we obtain $d=0$. Thus, we have $j(f_5)=2e$, where $e\in \Pic(S)$ is the class of the curve $E$.  There are two possibilities, either $e\in j(M)$ or $e\not \in j(M)$. The former case is not possible, since $f_5$ is not divisible in $M$ (see the proof of \ref{picequal}). The latter implies that the embedding $j:M\hookrightarrow \Pic(S)$ is not primitive (see the proof of \ref{primitive}), but this contradicts the assumptions. Therefore, the complete linear system $|H|$ is base point free.
\end{proof}

The following lemma shows that any $M$-polarization can be normalized.
\begin{lemma}\label{arrangegood}
Let $(S,j)$ be an $M$-polarized  $K3$ surface. Then there exists a $\phi\in \pm W(S)$ such that the composite map
$$M\xrightarrow{j} \Pic(S) \xrightarrow{\phi} \Pic(S)$$
defines a  normalized $M$-polarization for $S$. 
\end{lemma}
\begin{proof}
By acting with $\pm 1$, we can assume that $j(h)\in V^+(S)$. It is known that 
$$C(S)=\{x\in V^+(S)\mid x.\delta\ge 0 \textrm{ for all } \delta\in \Delta^+(S)\}$$ 
is a fundamental domain for the action of $W(S)$ on $V^+(S)$ (e.g. \cite[pg. 313]{barth}). Thus, there exists a $\phi\in W(S)$ such that changing the embedding by $\phi$, we have $\phi(j(h))$ is nef. It follows that we can assume $j(h)$ is nef. Acting with reflection in roots orthogonal to $j(h)$ preserves this condition. 

 Let $R$ be the sublattice of $\langle j(h)\rangle^\perp_{\Pic(S)}$ spanned by the roots. We have $j(e_i)\in R$ and $R$ is an even negative definite root lattice. In particular, $R$ decomposes as a direct sum $R_1\oplus\dots\oplus R_k$ of irreducible root systems of type $A$--$D$--$E$. We can assume that $j(e_1)\in R_1$. By acting with $W(R_1)$, we  can arrange that $j(e_1)$ is the highest root of $R_1$ (N.B. $\Delta^+(S)$ determines the set of positive roots for $R_1$). Since $j(h)$ is fixed by $W(R_1)$, $j(h)$ remains nef. By construction, $j(e_1)$ satisfies the second condition of definition \ref{goodmark}, i.e. 
$\delta\in  \langle h \rangle^{\perp}_{\Pic(S)}\cap \Delta^+(S)\ \Longrightarrow\ \delta.j(e_1)\le 0$. The claim follows by repeating  the process for $j(e_i)$ (for $i=2,\dots,5$), but acting only with reflections that stabilize $j(h), \dots, j(e_{i-1})$.
\end{proof} 

Note that a normalized $M$-polarization is essentially unique.
\begin{lemma}\label{lemmaunique}
Assume that $S$ is a $K3$ surface with two normalized $M$-polarizations $j$ and $j'$. Assume that $j(h)=j(h')$ and $j(l')=j'(l')$. Then up to a permutation of the labeling of $e_i$, we have $j=j'$ (i.e. $j(e_i)=j'(e_{\sigma(i)})$ for a permutation $\sigma$).
\end{lemma}
\begin{proof}
Note first that for root $\delta\in \Pic(S)$ such that $\delta.j'(h)\neq 0$ and $j'(h)\not\perp j'(M)$ we must have $\delta.j'(e_i)\neq 0$ for some $i$.  By the definition of normalized polarization, $j(e_1)$ is the highest root of an irreducible summand  in $\langle j(h)\rangle^\perp_{\Pic(S)}$ (w.r.t. $\Delta^+(S)$). Since $j(l')=j'(l')$ and $j(l').j(e_1)=1$, we conclude that $j(e_1).j'(e_i)\neq 0$ for some $i$. Thus, $j(e_1)$ and $j'(e_i)$ belong to the same irreducible root system in $\langle j(h)\rangle^\perp_{\Pic(S)}=\langle j'(h)\rangle^\perp_{\Pic(S)}$. Since $j(e_1)$ is the highest root, from the fact that $j'$ is also normalized we conclude that  $j(e_1)=j'(e_k)$ for some $k\le i$. We can assume $k=1$, and the argument can be repeated for the remaining $e_i$. 
\end{proof}

The last preliminary result shows that a normalized $M$-polarized  $K3$ surface $(S,j)$ satisfies the geometric properties listed in remark \ref{remarkgeom} (esp. ii). This is an important fact that allows us to pass from an $M$-polarization to a degree $5$ pair. We note that the basic idea of the proof of the lemma is that the condition of normalized polarization forces all the (-2)-curves orthogonal to $j(h)$, but not orthogonal to $j(M)$, to be components of $j(e_1),\dots,j(e_5)$. It follows then that $j(l')$ is irreducible.

\begin{lemma}\label{irredl}
Assume that $(S,j)$ is a normalized $M$-polarization. Then $j(l')$ is the class of an irreducible curve (thus, smooth and rational). 
\end{lemma}
\begin{proof}
As before, we represent the classes $j(l'),j(h),\dots,j(e_5)$ by effective divisors  $H,L',D_1, \dots ,D_5$. 
By lemma \ref{pseudoample}, we get that  $H$ is a base point free polarization. Assuming that $L'$ is not irreducible, we have a decomposition: 
\begin{equation}\label{decomposel}
L'=\sum n_i C_i+\sum m_i R_i
\end{equation}
where $n_i, m_i\ge 0$, and $C_i$ and $R_i$ are irreducible curves with $C_i^2\ge 0$ and $R_i^2=-2$.

We note first that there cannot be any non-rational curve occurring in  (\ref{decomposel}). From Hodge index theorem, we have $C_i^2\cdot H^2\le (H.C_i)^2$ with equality only if $C_i$ and $H$ are proportional. If $C_i^2>0$, we get $H.L'\ge H.C_i\ge 2$, contradicting the assumption $H.L'=1$. Similarly, if $C_i^2=0$ we get $C_i.H=1$, a contradiction to $|H|$ is base point free (see \ref{conversemayer}). It follows that $L'$ has the following decomposition:
\begin{equation}\label{decomposel2}
L'=R_0+\sum_{i=1}^k m_i R_i
\end{equation}
where $R_0,\dots,R_k$ are irreducible rational curves such that $H.R_0=1$ and $H.R_i=0$ for 
$i=1,\dots,k$. Since $R_0^2=L'^2=-2$, we get 
\begin{equation}\label{remainderl}
R_0.\left(\sum m_i R_i\right)=-\frac{1}{2}\left(\sum m_i R_i\right)^2\ge 1
\end{equation}
Thus, $R_0.L'\ge -1$. On the other hand, it is easy to see that $R_0.D_i\ge 0$. We conclude $1=R_0.H\ge 2R_0.L'$, which together with (\ref{remainderl}) gives $R_0.L'\in\{0,-1\}$.

Assume that $R_0.L'= 0$. Using $H.R_0=1$, we get $D_i.R_0=1$ for exactly one value of $i$, and  $D_i.R_0=0$ for the remaining values. Assume first that  $D_1.R_0=0$. Since $D_1.L'=1$, we get $D_1.\left(\sum m_i R_i\right)=1$. By assumption the embedding $j$ is normalized. In particular, $D_1$ is the highest root in $\langle j(h)\rangle^{\perp}$. Since  $R_i$ are effective roots in $\langle h\rangle^{\perp}$, we get $D_1.R_i\le 0$ for all $i=1,\dots k$. Thus, $1=D_1.\left(\sum m_i R_i\right)\le 0$, a contradiction. It follows that we must have $D_1.R_0=1$. The same argument as above gives $D_1.R_i=0$ for all $i$.  Then $R_i$ are effective roots in  $\langle j(h), j(e_1)\rangle^{\perp}$ and we obtain a contradiction to $D_2.R_0=0$.

The case $R_0.L'= -1$ is handled by a similar computation. We omit the details. In conclusion, the decomposition (\ref{decomposel2}) is trivial, i.e. $L'$ is irreducible. \end{proof}

We can now conclude that any $M$-polarized $K3$ surface is of type $\scl$.
\begin{proposition}\label{surjectivity}
Let $S$ be an $M$-polarized $K3$ surface.  Then, there exists a degree $5$ pair $(C,L)$ such that  $S\cong \scl$.
\end{proposition}
\begin{proof}
We are given a $K3$ surface $S$ together with a primitive lattice embedding $j:M\hookrightarrow \Pic(S)$. Without loss of generality, we can assume that $j$ is a normalized embedding (cf. \ref{arrangegood}). In particular, $j(h)$   is nef  and  $j(e_i)$ are classes of effective divisors.  It follows that $j(h)$ is base point free (cf. \ref{pseudoample}). Thus, $j(h)$ 
defines a generically $2:1$ morphism $\pi:S\to \bP^2$, with branch curve a sextic $B$. 

By lemma \ref{irredl}, the class $j(l')$ contains a (unique) smooth rational curve $L'$. It follows that $L=\pi_*L'$ is a line in $\bP^2$.  To prove that $L$ is a component of the branch locus $B$, it is enough to prove that class $j(l')$ of $L'$ is invariant under involution $i$, where $\tau$ is the natural sheet-exchanging  involution on $S$ and $i=\tau^*$ the induced involution on $H^2(S,\bZ)$. The polarization class  $j(h)$ is left invariant by $i$. We claim that also the classes $j(e_i)$ are invariant. Namely, the morphism $\pi$ will contract all the $(-2)$-curves orthogonal to $j(h)$ to singularities of the sextic $B$. From the normalization assumption, it follows that $j(e_1)$ is the fundamental cycle of a singularity of $B$. Thus, $j(e_1)$ is invariant. A similar argument works for all $e_i$. Since  $2l'=h-e_1-\dots-e_5$, we conclude that $j(l')$ is also invariant as needed.\end{proof}

\subsubsection{The proof of theorem \ref{mainthm3}}\label{secendproof}
First, the surjectivity follows from proposition \ref{surjectivity}, which assures us that the construction of \S\ref{quinticperiod} can be reversed. Given a point $w\in \calD/\Gamma$,
we choose a lift $\omega\in T_\bC\subset \Lambda_\bC$. This determines a lattice $P:=\Lambda\cap\langle \omega\rangle^{\perp}$ and a factorization of the fixed primitive embedding $M\hookrightarrow \Lambda$ into
$M\hookrightarrow P\hookrightarrow \Lambda$. To associate a $K3$ surface $S$ to $\omega$, we have to provide a choice for $V^+(P)$ and $\Delta^+(P)\subset \Delta(P)$. We define $V^+(P)$ to be the connected component of $V(P)$ that contains $j(h)$. Then, we define the partition $\Delta(P)=\Delta^+(P)\sqcup (-\Delta^+(P))$ by: 
\begin{itemize}
\item[i)] if $\delta\in\Delta(P)$ and  $\delta.j(h)\neq0$: $\delta\in  \Delta^+(S)$ if  $\delta.j(h)>0$ or  $-\delta\in  \Delta^+(P)$ otherwise;
\item[ii)] if $\delta\in \langle j(h),\dots,j(\delta_{i-1})\rangle^{\perp}_P\cap\Delta(P)$ and $\delta.j(\delta_i)\neq 0$:  $\delta\in  \Delta^+(P)$ if  $\delta.j(\delta_i)<0$ or  $-\delta\in  \Delta^+(P)$ otherwise (for $i=1,\dots,5$ and $\delta_0=h$);
\item[iii)] For $\Delta(P)\cap j(M)^\perp_P$ we choose an arbitrary Weyl chamber.
\end{itemize} 
 By the surjectivity of the period map for $K3$ surfaces, there exists a $K3$ surface $S$ with period point $\omega$ and such that $\Pic(S)=P$, $V^+(S)=V^+(P)$, and $\Delta^+(S)=\Delta^+(P)$. By construction $j:M\hookrightarrow \Pic(S)$ is a normalized $M$-polarization. Thus, there exists a degree $5$ pair $(C,L)$ such that $S\cong \scl$ (cf. \ref{surjectivity}).  

By global Torelli theorem, the surface $S$ obtained above is unique up to isomorphism. The lattice $M$ is left invariant by $O^*(T)$. It is easily seen that  the classes $h$ and $l'$ are left invariant by the full monodromy group $\Gamma=O_{-}(T)$ (see \ref{propmonodromy}).  It follows that the period point $\omega\in \calD/\Gamma$ uniquely determines the classes $j(h), j(l')\in \Pic(S)$.  The polarization $j(h)$ determines the double cover map $S\to \bP^2$ with branch curve $B$, and $j(l')$ determines a line component $L$ of $B$. Thus, $B=C+L$, and we conclude that a period point $\omega$ determines uniquely (up to projective isomorphism) a degree $5$ pair $(C,L)$. \qed

\subsubsection{The proof of theorem \ref{mainthm4}}\label{proofmainthm4}
E. Looijenga \cite{looijengacompact} has developed a general framework of comparing GIT compactifications to appropriate compactifications of the period space. Specifically,  \cite[Thm. 7.6]{looijengacompact}  says that once a period map gives an isomorphism $\calM\cong (\calD\setminus \calH)/\Gamma$ between a geometric quotient $\calM$ and the complement of an arithmetic arrangement $\calH$ of hyperplanes in a type IV domain, it automatically extends to an isomorphism $\overline{\calM}\cong \widetilde{\calD/\Gamma}$ between the GIT compactification of $\calM$ and the Looijenga's compactification associated to $\calH$.

 In our situation, we apply  \cite[Thm. 7.6]{looijengacompact} to the empty arrangement of hyperplanes. By theorem \ref{mainthm3}, we have an isomorphism $\calM\cong \calD/\Gamma$. Since the Looijenga's compactification associated to the empty arrangement of hyperplanes is precisely the Baily-Borel compactification $(\calD/\Gamma)^*$, we conclude that theorem \ref{mainthm4} holds once the technical assumptions of  \cite[Thm. 7.6]{looijengacompact} are verified. In our situation, this is easy. First, the codimension condition for the GIT quotient is clearly satisfied. The complement of the $G$-invariant open subset $U$  parameterizing  degree $5$ pairs $(C,L)$ such that $C+L$ has simple singularities  has high codimension in $X^{ss}(1)$. The second assumption in  \cite[Thm. 7.6]{looijengacompact} is that the isomorphism $\calM\cong \calD/\Gamma$ is an isomorphism of polarized varieties (both spaces are naturally polarized). Since $\calM$ is an open subset of the GIT quotient $\calM(1)$, its polarization is obtained  by restricting the polarization of the moduli space of plane sextics (see Thm. \ref{mainthm1}). Similarly,  the polarization on $\calD/\Gamma$ is obtained by restricting the polarization of the moduli space of degree two $K3$ surfaces. Thus the identification of polarizations on $\calM$ and $\calD/\Gamma$ follows from Looijenga's computation \cite[\S8]{looijengacompact}  (esp. \cite[Thm. 8.6]{looijengacompact}) for degree two $K3$ surfaces. \qed

\subsection{The Baily-Borel Compactification}\label{sectbailyborel}
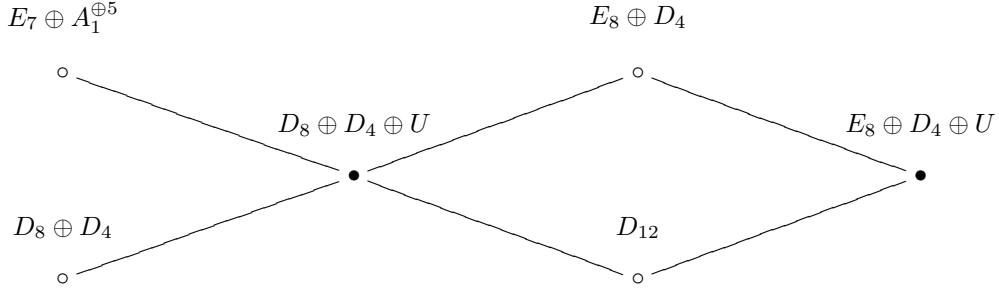
\begin{figure}[htb]
$$\xymatrix@R=.25cm{
{E_7\oplus A_1^{\oplus 5}}   &&                                   &&{E_8\oplus D_4}     &&                       \\
{\circ}\ar@{-}[ddrr]&&                                   &&{\circ}\ar@{-}[ddrr]&&                       \\
                    &&{D_8\oplus D_4\oplus U}            &&                    &&{E_8\oplus D_4\oplus U}\\
                    &&{\bullet}\ar@{-}[uurr]\ar@{-}[ddrr]&&                    &&{\bullet}              \\
{D_8\oplus D_4}     &&                                   &&{D_{12}}            &&                       \\
{\circ}\ar@{-}[uurr]&&                                   &&{\circ}\ar@{-}[uurr]&&                       \\
}
$$
\caption{The boundary components of the Baily-Borel compactification}\label{bbcomp}
\end{figure}

The quotient of a bounded symmetric domain by an arithmetic group  $\calD/\Gamma$ admits a canonical minimal compactification, the Baily-Borel compactification $(\calD/\Gamma)^*$. In the case of Type IV domains, the boundary components  of $\left(\calD/\Gamma\right)^*$ are either $0$-dimensional (type III) or $1$-dimensional (type II), and  are in  bijective correspondence with the  equivalence classes of the primitive isotropic sublattices of $T$ of rank $1$ and $2$ respectively. Thus, to determine the  Baily-Borel compactification  is a pure arithmetic question, which in our situation has the following answer:
\begin{theorem}\label{bailyborel}
The boundary of $\calD/\Gamma$ in the Baily-Borel compactification $\left(\calD/\Gamma\right)^*$ consists of two $0$-dimensional components and  four $1$-dimensional component. Their incidence graph is given in figure \ref{bbcomp}. 
\end{theorem}
\begin{proof}
The $0$-dimensional boundary components correspond to the classes of rank $1$ isotropic lattices of $T$. These are classified in \S\ref{sectinv} (esp.  \ref{brlemma}). Lemma  \ref{split} reduces the classification of isotropic rank $2$ lattices to the classification of isotropic vectors in the hyperbolic lattice $N:=E_8\oplus D_4\oplus U(2)$. We conclude by Vinberg's algorithm applied to  $N$ (see \S\ref{vinberg}). \end{proof}

According to theorem \ref{mainthm4}, the projective varieties $\calM(1)$ and $(\calD/\Gamma)^*$ are isomorphic. Since they compactify the same space $\calM\cong \calD/\Gamma$, the two sets of boundary components are isomorphic. Geometrically, the matching of the  boundary components (figures \ref{incidence} and \ref{bbcomp}) is obtained as follows. Let $(C_0,L_0)$ be a degree $5$ pair corresponding to a boundary point of $\calM(1)\setminus \calM$. We consider  a generic  family of degree $5$ pairs $(C_t,L_t)$ degenerating to $(C_0,L_0)$. This  produces a family of $K3$ surfaces $S_t$ (where $S_t=\widetilde{S}_{(C_t,L_t)}$) degenerating to $S_0$, the double cover corresponding to $(C_0,L_0)$. To the family of $K3$ surfaces $(S_t)_{t\in \Delta^*}$ there is associated a canonical limit mixed Hodge structure $H^2_{\lim}$, which determines a boundary point in the Baily-Borel compactification. Since the period map extends to the boundary, the resulting point does not depend on the choice of  the degenerating family. Thus, to compute the limit mixed Hodge structure $H^2_{\lim}$ it is enough to do semistable reduction and apply the Clemens-Schmid sequence (e.g. \cite{clemensschmid}) for a suitable degeneration to $(C_0,L_0)$. One distinguishes two cases, type II and III, depending on the monodromy of the family. Using  the incidence relation of the boundary components, the matching for type III case is obtained from that for type II. Finally, the computation in the type II case is very similar to that for degree two $K3$ surfaces (e.g. \cite[Rem. 5.6]{friedmanannals}). We conclude:

\begin{corollary}\label{matchboundary}
Via the extended period map $\overline{\calP}:\calM(1) \to (\calD/\Gamma)^*$, the boundary components of the GIT quotient $\calM(1)$  map to boundary components of the Baily-Borel compactification  $(\calD/\Gamma)^*$ as given in table \ref{comparebound}. \qed
\end{corollary}

\begin{table}[htb]
\begin{center}
\begin{tabular}{ll}
\hline
{\it GIT boundary (see figure \ref{incidence})}& {\it Baily-Borel boundary (see  figure \ref{bbcomp})}\\
\hline
II(1)&$E_8\oplus D_4$\\
II(2a)&$E_7\oplus A_1^{\oplus 5}$\\
II(2b)&$D_4\oplus D_8$\\
II(3)&$D_{12}$\\
\hline
III(1)&$E_8\oplus D_4\oplus U$ \\
III(2)& $D_8\oplus D_4\oplus U$\\
\hline
\end{tabular}
\end{center}
\caption{The matching of the boundary components}\label{comparebound}
\end{table}

\subsubsection{Invariants of isotropic sublattices}\label{sectinv}
Explicit computations of the Baily-Borel compactification for type IV domains are given in \cite{brieskorn, scattone, sterk}. The computations in our situation are similar and we only sketch the arguments. 

To compute the Baily-Borel compactification of the period space $\calD/\Gamma$ we have to classify the isomorphism classes of isotropic sublattices $E$ of the transcendental lattice $T\cong D_4\oplus E_8\oplus U(2)\oplus U$ modulo $O(T)$ (N.B. $\Gamma=O_{-}(T)$ and $\bZ/2\bZ\cong O(T)/O_-(T)$ correspond to the choice of component for $\calD$). We start by noting the following invariants for an isotropic lattice $E$ (see  \cite{brieskorn,scattone}):
\begin{itemize}
\item[-] {\it the rank $k\in\{1,2\}$ of $E$};
\item[-] {\it the isotropic subgroup $H_E:=E^{\perp\perp}_{T^*}/E$} of the discriminant group $A_T$;
\item[-] {\it the isomorphism class of the lattice  $E^{\perp}/E$}.
\end{itemize} 
The main conclusion of our computations below is that these invariants completely determine the isomorphism classes of the isotropic sublattices  of $T$. 

For the moment, let us describe the possible values for these invariants. Since the discriminant group $A_T$ is isomorphic to $(\bZ/2\bZ)^4$, it follows easily that $H_E$ is either trivial or isomorphic to $\bZ/2\bZ$. Furthermore, in the latter case,  all possible subgroups $H_E$ are conjugate by $O(T)$. Fixing the rank $k$ and the group $H_E$ is equivalent to specifying the signature and the discriminant group for the lattice $E^{\perp}/E$ (N.B. $A_{E^{\perp}/E}\cong  H^{\perp}_E/H_E$). In particular, since the lattice $E^{\perp}/E$ is $2$-elementary, it follows that $k$ and $H_E$  determine the genus of $E^{\perp}/E$ (\cite[Thm. 3.6.2]{nikulin}). In the rank $1$ case,  the lattice $E^{\perp}/E$ is indefinite and its isomorphism class is completely determined by the genus. Thus, if $E$ is isotropic of rank $1$, then  $E^{\perp}/E$ is  isomorphic to either  $D_4\oplus D_8\oplus U$ or $D_4\oplus E_8\oplus U$. In the rank $2$ case, the two possibilities for $H_E$ determine two genera for $E^{\perp}/E$, namely $\calG(D_4\oplus D_8)$ and $\calG(D_4\oplus E_8)$, which in turn  give $4$ possibilities for the isomorphism classes.

\begin{lemma}\label{lemmagenus} 
The genera of $\calG(D_4\oplus D_8)$ and $\calG(D_4\oplus E_8)$ contain exactly two isomorphism classes of even negative definite lattices, namely
\begin{itemize}
\item[i)] $\calG(D_4\oplus E_8)$ consists of the root lattices $D_4\oplus E_8$ and $D_{12}$;
\item[ii)] $\calG(D_4\oplus D_8)$ consists of the lattice $D_4\oplus D_8$ and  an index $2$ overlattice 
of $E_7\oplus A_1^5$. 
\end{itemize}
\end{lemma}
\begin{proof}
By  \cite[Prop. 6.1.1.]{scattone}, every lattice in the genus of $L=D_4\oplus E_8$ can be obtained as the orthogonal complement $L^{\perp}_\Lambda$ for a suitable primitive embedding  of $L$ into  an  even unimodular negative definite lattice $\Lambda$ of rank $24$. Since $L$ is a root lattice, the embedding $L\hookrightarrow \Lambda$ factors as $L\hookrightarrow \underline{\Lambda} \subseteq \Lambda$, where $\underline{\Lambda}$ is the sublattice spanned by the roots of $\Lambda$.  The possibilities for $\Lambda$ (and $\underline{\Lambda}$) are classified by Niemer's theorem. Finally, since the embeddings of root lattices are well understood, we conclude by a case-by-case analysis. The item ii is similar.\end{proof}

It is easy to see that each choice of admissible invariants  actually corresponds to some isotropic sublattice $E$. For example, there exists a rank $2$ isotropic sublattice $E$ such that $E^\perp/E\cong D_{12}$ by the following argument. Since the lattice $T$ is determined by its invariants, there exists an isomorphism $T\cong D_{12}\oplus U\oplus U(2)$. By choosing an isotropic vector in each of the hyperbolic summands, we obtain a rank $2$ isotropic sublattice $E$ having the right invariants. In conclusion, we have obtained $6$ boundary components for the Baily-Borel compactification satisfying the incidence relations of figure \ref{bbcomp} (N.B. the label in figure \ref{bbcomp} corresponds to the isomorphism class of  $E^\perp/E$). The remaining part for the proof of theorem \ref{bailyborel} is to see the converse. There exists at most one class of isotropic sublattices having some prescribed invariants. For the rank $1$ case, this is automatic by the following result of Brieskorn
 \cite[Kor. 2, pg. 87]{brieskorn}:
\begin{lemma}\label{brlemma}
Let $T$ be an even  lattice containing at least two hyperbolic planes such that the natural map $O(T)\to O(q_T)$ is surjective. Then the classes of rank $1$ isotropic sublattices modulo $O(T)$ are in bijection with the classes of isotropic vectors in $A_T$ modulo $\pm O(q_T)$.
\end{lemma}

In our situation, $T\cong D_4\oplus D_8\oplus U^2$ and  $O(T)\to O(q_T)$ is surjective, thus we can apply the previous 
lemma. The bijection of the lemma is given by associating to an isotropic lattice $E$ the generator of $H_E$. Thus, we obtain indeed only two classes of isotropic rank $1$ sublattices (distinguished by the invariant $H_E$).

\subsubsection{Rank $2$ isotropic sublattices}\label{vinberg} Let $N\cong E_8\oplus D_4\oplus U(2)$ be the sublattice of $T$ obtained by splitting off a hyperbolic plane from $T$. We claim that the classification of the isotropic rank $2$ sublattices of $T\cong N\oplus U$ is essentially equivalent to the classification of the isotropic  rank $1$ sublattices in $N$. The idea is simple. Namely, given $E$ a rank $2$ isotropic sublattice we can choose a rank $1$ sublattice $E'$ in $E$ and use the classification of $E'$ given by lemma \ref{brlemma} to reduce to $N$. More precisely, we have:
  
\begin{lemma}\label{split}
Let $E$ be a rank $2$ isotropic sublattice of $T$. Then there exists a basis $\{b_1,b_2\}$ of $E$ and isometry $\phi:T\to U\oplus N$ such that   $\phi(b_1)\in U$ and $\phi(b_2)\in N$. \qed
\end{lemma}

Thus, to classify $E$ it is enough to  classify the isotropic  rank $1$ sublattices in $N$ modulo $O(N)$. Note that going back from an isotropic vector in $N$ to an isotropic rank $2$ sublattice in $T$ might give some repetitions, but these  are easy to detect. In conclusion, we are done once we classify  the isotropic vectors in $N$.

\begin{figure}[htb]
\begin{center}
\scalebox{.75}{
$$\xymatrix@R=.25cm{
                    &&                    &&{\bullet}\ar@{-}^{\infty}[rrrr]&&&&{\bullet}           &&      &&\\
                    &&                    &&                                          &&&&                       &&      &&\\
                    &&                    &&{\bullet}\ar@{-}^{\infty}[rrrr]&&&&{\bullet}           &&      &&\\
                    &&                    &&                                          &&&&                       &&      &&\\
{\bullet}\ar@{-}[rr]\ar@{-}[ddddddd]&&{\bullet}\ar@{-}[rruuuu]\ar@{-}[rruu]\ar@{-}[rr]\ar@{-}[rrdd]\ar@{-}[rrdddd]&&{\bullet}\ar@{-}[rrrr]^{\infty}&&&&{\bullet}\ar@{-}[rr]&&{\bullet}\ar@{-}[rr]\ar@{-}[lluuuu]\ar@{-}[lluu]\ar@{-}[lldd]\ar@{-}[lldddd]
&&{\bullet}\ar@{-}[ddddddd]\\
                    &&                    &&                                          &&&&                       &&      &&\\
                    &&                    &&{\bullet}\ar@{-}^{\infty}[rrrr]&&&&{\bullet}           &&       &&\\
                    &&                    &&                                         &&&&                        &&       &&\\
                    &&                    &&{\bullet}\ar@{-}^{\infty}[rrrr]&&&&{\bullet}           &&       &&\\
                    &&                    &&                                          &&&&                       &&       &&\\
                    &&                    &&                                          &&&&                       &&       &&\\
{\bullet}\ar@{-}[rrr]&&&{\bullet}\ar@{-}[rrr]&&&{\bullet}\ar@{-}[rrr]\ar@{-}[dd]&&&{\bullet}\ar@{-}[rrr]&&&{\bullet}\\
                    &&                    &&                                          &&&&                       &&       &&\\
                                        &&                    &&                                          &&{\bullet}&&                       &&       &&\\
}
$$}
\caption{Vinberg Diagram for $N$}\label{vindiag}
\end{center}
\end{figure}
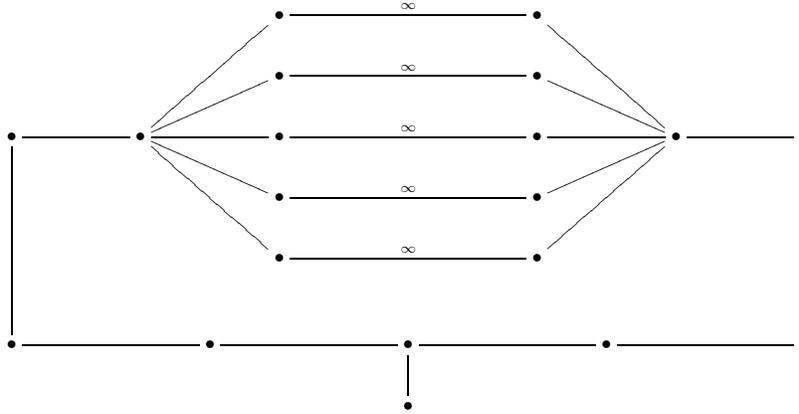 

Since $N$ contains only one hyperbolic plane,  lemma \ref{brlemma} cannot be applied. Instead, we use the fact that $N$ is a hyperbolic lattice, i.e. it has  signature $(1,n)=(1,13)$. For such lattices  there exists a classification algorithm for the isotropic vectors due to Vinberg \cite[\S 1.4]{vinberg}  (see also \cite[\S 4.3]{sterk}). The algorithm starts by fixing an element $h\in N$ of positive square. Then, the algorithm consists of  inductively choosing roots $\delta_1,\delta_2,\dots$ such that  a certain distance function to $h$ is minimized. The algorithm stops with the choice of the root $\delta_N$ if the following stop condition (\cite[Thm. 2.6bis]{vinberg}) is satisfied: {\it every connected parabolic subdiagram of the Dynkin diagram $\Sigma$ associated to the roots $\delta_1,\dots,\delta_N$ is a connected component of some parabolic subdiagram of rank $n-1$ (i.e. the maximal rank)}. If the algorithm stops, the classes of isotropic vectors in $N$ are determined by the rule: {\it every isotropic line  $E$ in $N$ is equivalent to the null space  of the lattice corresponding to a parabolic subdiagram of rank $n-1$ of the Dynkin diagram $\Sigma$}. We recall that by a parabolic diagram we understand the extended Dynkin diagram of a root system, and the rank is the rank of the corresponding root system.

In our situation, a straightforward application of Vinberg's algorithm for $N$ produces the Dynkin diagram $\Sigma$ given in figure \ref{vindiag}. The computation is simplified by the observation $N\cong M\oplus E$. This gives the natural  choice $h\in M$ (the polarization class used in \S\ref{quinticperiod}) to start the algorithm. Note also that the Dynkin diagram $\Sigma$ satisfies the stopping criterion. Up to the obvious symmetries,  there are only $4$ distinct possibilities for the maximal rank parabolic subdiagrams  of the Dynkin diagram of figure \ref{vindiag}, namely $\widetilde{A}_1^{\oplus 5}\oplus \widetilde{E}_7$, $\widetilde{D}_8\oplus \widetilde{D}_4$, $\widetilde{E}_8\oplus \widetilde{D}_4$, and $\widetilde{D}_{12}$. Together with lemma \ref{lemmagenus}, this implies that there are exactly $4$ classes of isotropic rank $2$ sublattices of $T$ modulo $O(T)$. This concludes the proof of theorem \ref{bailyborel}.

\section{Applications to the deformations of $N_{16}$}\label{sectn16}
In this section,  we return to our original motivation, the study of the deformations of $N_{16}$. We consider the germ $(Y_0,y_0)$ of a singularity of type $N_{16}$ and let   $\calY\to S$ be the semi-universal deformation. Since $Y_0$ is a hypersurface singularity, the base $S$ is the germ of a smooth $16$-dimensional space ($\mu=\tau=16$ for $N_{16}$). The main question that one would like to understand is the natural stratification of $S$ in terms of the singularities of the nearby fibers $Y_s$. We recall that there exists an open dense stratum of $S$ which parameterizes the smooth fibers $Y_s$, whose complement  is called the discriminant hypersurface $\Sigma\subset S$. The discriminant $\Sigma$ is stratified in terms of the complexity of the singularities occurring for $Y_s$ ($s\in \Sigma$). For instance, we can discuss about the equisingular stratum,  the stratum of simple singularities, the stratum of simple elliptic singularities, etc. The basic question, then, is to say something about the structure of this stratification. The easiest and the best understood situation is that of the deformations of the simple singularities. In that case, $S$ can be taken to be the quotient of affine space $\bA^\mu$ by a finite group $W$, and the discriminant corresponds to $\calH/W$, where $\calH$ is an arrangement of hyperplanes in $\bA^\mu$. Furthermore, the stratification of the discriminant corresponds to the natural stratification given by the intersections of hyperplanes from $\calH$ (e.g. \cite[Ch. 5]{arnold}). By work of Pinkham \cite{pinkham,pinkhammonodromy}, Looijenga \cite{looijenga,looijengatriangle}, Brieskorn \cite{novaacta}, and others, a rather similar situation was shown to hold also for the unimodal (e.g. triangle) singularities. In our situation, putting together the results of the previous sections, we show that the structure of the deformations of $N_{16}$ is very similar to the structure of the deformations of the triangle singularities. The only major difference is that one has to work in a relative setting over the equisingular deformations. 

In our study of the deformations of $N_{16}$, we make the following two standard restrictions (see \cite{looijengatriangle}, \cite[\S10]{looijengacompact}). First, we consider only the case of quasi-homogeneous $N_{16}$. This restriction implies the existence of a natural $\bC^*$-action on the base $S$ of the versal deformation, which in turn determines subspaces $S_{\le 0}$, $S_0$, $S_{+}$, etc. (see \S\ref{sectpinkham}). In our situation, the space of the deformations of positive weight $S_{+}$ (the equisingular deformations moving away from the $\bC^*$-action) is $1$-dimensional. The second restriction is to  ignore these deformations, i.e. we consider only the codimension $1$ subspace $S_{\le 0}$ of $S$. This type of restriction is harmless, at least from a topological point of view. By results of Wirthm\"uller and Damon, for non-simple hypersurface germs, the universal deformation is topologically trivial along the Hessian direction (i.e. the positive weight direction in our situation).

Due to the assumptions made in the previous paragraph, we are in the situation covered by the theory of deformations with $\bC^*$-action of Pinkham \cite{pinkham0}. It follows that the non-positive weight deformations of $N_{16}$ can be identified with the projective deformations of the pair $(C_0,L_0)$, where $C_0\subset \bP^2$ is the cone over $5$ points in $\bP^1$ and $L_0$ is a generic line. By gluing together these deformations in a stack, we obtain a modular interpretation for the deformations of non-positive weight of $N_{16}$. As discussed in section \ref{sectgen} (esp. \S\ref{sectpinkham}), we obtain a diagram:
\begin{equation}\label{defsquare2}
\begin{CD}
(S_{\le 0}\setminus S_0)/\bC^*@>>>\calM_{\le 0}@>>>\calM(\frac{5}{2}-\epsilon)\\
@VVV @VVV@VVV\\
S_0@>>>\calM_0@>>> \calM(\frac{5}{2})
\end{CD}
\end{equation}
where $\calM_0$ is the moduli of $5$ distinct points in $\bP^1$, $\calM_{\le 0}$ is the moduli space of pairs $(C,L)$ such that $C$ is not a cone and $C\cap L$ is transversal, and the map $\calM_{\le 0}\to \calM_0$ is the natural forgetful map sending $(C,L)$ to the intersection $C\cap L$. The left side square of (\ref{defsquare2}) is the usual process of passing from deformations to the (coarse) moduli space, and the right side square is the compactification given by the variation of GIT quotients. Note that everything is done in a relative version over the $0$ weight deformations (the equisingular deformations preserving the $\bC^*$-action). By restricting to fibers, we obtain the deformations of negative weight (the smoothing directions) and the usual globalization of Pinkham. For comparison, we note that for triangle singularities, $S_0$ is trivial (i.e. a point).  

As mentioned above, we are concerned here with the natural stratification of $\calM_{\le 0}$. For this, we introduce the following notations.
\begin{notation}
Let $\calM_{r}$ be {\it the regular locus}, i.e. the moduli space of pairs $(C,L)$ such that $C$ is smooth and $L$ is transversal, and $\calM_{f}$ {\it the simple singularity locus} ($C$ has at worst simple singularities). Both spaces are open subsets of $\calM_{\le 0}$, and we have: $\calM_r\subset\calM_f\subset \calM_{\le 0}$. Note also that $\calM_r\subset\calM_f\subset \calM\cong \calD/\Gamma$ (Thm. \ref{mainthm3}).
\end{notation}

To describe the stratification of $\calM_{\le 0}$, we make use of the construction  of section \ref{sectk3} of the moduli space of degree $5$ pairs as the quotient of a bounded symmetric domain by an arithmetic group. We recall the isomorphism $\calM(1)\cong (\calD/\Gamma)^*$ (cf. \ref{mainthm3} and \ref{mainthm4}). Since  $\calM(1)$ and $\calM(\frac{5}{2}-\epsilon)$ are birational having $\calM_f$ a common open subset (cf. \ref{corquintic}), from the $\calD/\Gamma$ description 
of the moduli of pairs we obtain a good understanding of the regular and simple singularities strata $\calM_r$ and $\calM_f$.  The main conclusion, as in the case of simple and unimodal singularities, is that these spaces can be described as complements of arithmetic arrangements of hyperplanes (Thm. \ref{thmstructure}). Additionally, the possible combinations of simple singularities for a nearby fiber can be obtained in an algorithmic way  (see \S\ref{secsimple}). 

The structure of the strata $\calM_{\le 0}$ corresponding to the non-simple singularities adjacent to $N_{16}$ is obtained by  following the birational modifications that transform $\calM(1)$ into $\calM(\frac{5}{2}-\epsilon)$. The nature of the birational map $\calM(1)\dashrightarrow \calM(\frac{5}{2}-\epsilon)$ is both explicit and simple. As noted in \S\ref{symquintic}, the flips occurring  at the wall crossings between $1$ and $\frac{5}{2}$ are of the simplest type possible: at the critical slope $t\in (1,\frac{5}{2})$ the flip replaces a weighted projective subspace in $\calM(t-\epsilon)$ by a weighted projective subspace in $\calM(t+\epsilon)$ of complementary dimension. Geometrically, the introduced space corresponds to the stratum of a non-simple singularity. To be more precise, the strata of the triangle singularities adjacent to $N_{16}$ (i.e. $Z_{11}$, $Z_{12}$, $W_{12}$, and $W_{13}$) are introduced in this way (one at a time) at the slopes $\frac{10}{7}$, $\frac{8}{5}$, $\frac{5}{3}$, and $\frac{13}{7}$  respectively. The strata of simple elliptic and cusp singularities are introduced by the birational morphism $\calM(1+\epsilon)\to \calM(1)$. The center of this birational map is over the boundary of $\calD/\Gamma\cong \calM$ in the Baily-Borel compactification $(\calD/\Gamma)^*\cong \calM(1)$ (see \S\ref{quintic1} and \S\ref{sectbailyborel}). The structure of these strata can be made quite explicit by invoking Luna's slice theorem. We only note here that the simple elliptic strata (corresponding to $\widetilde{E_7}$ and  $\widetilde{E_8}$) are  weighted projective bundles over rational curves (type II boundary components) that parameterize the modulus of the corresponding simple-elliptic singularity.  

We close by noting the strong arithmetic nature of the stratification of $\calM_{\le 0}$. For the stratum of simple singularities, this is clear (see \ref{thmstructure} below). For the simple elliptic and cusp singularities, we have noted in \S\ref{quintic1} (esp. figure \ref{simpleellipticcusp}) the relation to the Baily-Borel compactification, analogous to what was observed by Brieskorn \cite{novaacta,brieskorn} for triangle singularities. Finally, the introduction of the strata corresponding to the triangle singularities is very similar to the compactification procedure of Looijenga \cite{looijengacompact} (esp. \cite[\S10]{looijengacompact}). This is discussed in \S\ref{sectriangle} below.

\begin{remark}
It is probably worthwhile to note what is general and what is special about the case of $N_{16}$. The results of section \ref{sectgen} and the connection between the variation of GIT quotients and the deformations with $\bC^*$-action are quite general. They hold for any degree $d$ cones and even in higher dimensions. The basic facts behind this are: the singularities are naturally stratified by the log canonical threshold, the log canonical threshold is closely related to the GIT stability condition, the cones are the worst singularities and, finally, a basic construction of the versal deformation space is as a slice of an appropriate Hilbert scheme. What is special about $N_{16}$ is the relation to the $K3$ surfaces. This is due to the fact that $N_{16}$ is a simple elliptic singularity, and as such it shares quite a few common characteristics (esp. of cohomological nature) with the triangle singularities. 
\end{remark}

\subsection{Discriminants}\label{secdiscriminants}
We start by identifying the complements (the discriminants)  of the regular and simple singularities locus  $\calM_r$ and  $\calM_f$ respectively in $\calD/\Gamma$.

\begin{lemma}
Let $(C,L)$ be a degree $5$ pair such that $C+L$ has at worst simple singularities, and let $j:M\hookrightarrow \Pic(\scl)$ be the $M$-polarization constructed in \ref{markinggen}. The following conditions are equivalent:
\begin{itemize}
\item[i)] the curve $C$ is smooth and the intersection $C\cap L$ is transversal.
\item[ii)] for every $\delta\in \Pic(\scl)\setminus j(M)$ we have $\delta.j(h)\neq 0$. 
\end{itemize}
\end{lemma}
\begin{proof}
By construction, the polarization $j(h)$ defines the double cover map $\scl\to \bP^2$ with branch locus the sextic $B=C+L$. It is well known that the singularities of $B$ are in one-to-one correspondence with the irreducible summands of the root sublattice of $\langle j(h)\rangle^\perp_{\Pic(\scl)}$. Thus, the condition that $C$ is smooth and the intersection is transversal is equivalent to saying that there is no root $\delta\in \Pic(\scl)$ orthogonal to $j(h)$ except those coming from $M$.
\end{proof}

The previous lemma says that the non-generic pairs $(C,L)$ are detected by the existence of a root $\delta\in \Pic(\scl)$ orthogonal to the polarization class $h$. Arithmetically, one distinguishes two possibilities: either $\delta$ is orthogonal to $M$, or not.  The next lemma gives a geometric meaning to this division. 
\begin{lemma}\label{lemmalternative}
Let $(C,L)$ be a pair and $(\scl,j)$ the associated $K3$ surface. Then 
\begin{itemize}
\item[i)]  the quintic $C$ is singular if and only if there exists $\delta\in \Delta(\scl)\cap \langle j(h)\rangle^\perp$ such that $\delta\not \in j(M)$ and $\delta\perp j(M)$;
\item[ii)] the intersection $C\cap L$ is not transversal  if and only if there exists $\delta\in \Delta(\scl)\cap \langle j(h)
\rangle^\perp$ such that $\delta\not \in j(M)$ and $\delta\not \perp j(M)$.
\end{itemize}
\end{lemma}
\begin{proof}
The ``only if'' part follows by construction. Namely, the case when $C$ is singular away from $L$, or $L$ is tangent to $C$ is immediate. The only thing to check is that if $L$ passes through a singular point of $C$, there exist both a root which is orthogonal to $j(M)$, and another which is not. This follows easily from the construction of \ref{markinggen}. 

Conversely, assume that there exists $\delta\in \Delta(\scl)\cap \langle j(h)\rangle^\perp$ such that $\delta\not\in j(M)$. By the previous lemma, either $C$ is singular, or the intersection $C\cap L$ is degenerate. The claim is that we can distinguish between the two cases based on the fact that either $\delta$ is orthogonal to $j(M)$ or not. Since $\delta.j(h)=0$, the case  that $\delta$ is not orthogonal to $j(M)$ implies $\delta.j(e_i)\neq 0$ for some $i$. Thus, the irreducible summand of the root sublattice of $\langle j(h)\rangle_{\Pic(\scl)}^\perp$ containing $j(e_i)$ is larger than $A_1$. Since by construction $j(e_i)$ corresponds to a point of intersection $p\in C\cap L$, we obtain that  the intersection at $p$ is not transversal. 

Assume now that there exists a root $\delta\in \Pic(\scl)$ such that $\delta\perp j(M)$. Assume $C$ is smooth.  We obtain a contradiction as follows. Since all the singularities of $C+L$ come from the intersection $C\cap L$,   we must have that the intersection $C\cap L$ is degenerate at some point $p$ and that $\delta$ belongs to the root system corresponding to the singularity at $p$. Since $C$ is smooth,  the only geometric possibility is that $L$ is tangent  with multiplicity $k\ge 2$ to $C$ in the point $p$. Thus, the singularity at $p$ for $C+L$ is $A_{2k-1}$, which gives an embedding $A_{2k-1}\hookrightarrow \langle j(h)\rangle_{\Pic(\scl)}^\perp$.  Without loss of generality, we can assume $j(e_1),\dots,j(e_k)$ belong to this $A_{2k-1}$ root system  (see \ref{markinggen}). By assumption, $\delta$ also belongs to $A_{2k-1}$. We obtain that $\delta,\dots,j(e_k)$ span a sublattice of $A_{2k-1}$ isometric to $A_1^{\oplus (k+1)}$. This is not possible since there is no lattice embedding of $A_1^{\oplus (k+1)}$ into $A_{2k-1}$. \end{proof}

\begin{notation}\label{defdiscr} 
Fix a primitive embedding $M\hookrightarrow \Lambda$ (see \S\ref{mpolarized}). We denote $\Delta_h:=\{\delta\in \Delta(\Lambda)\mid \delta.h=0\}$ the set of roots orthogonal to the polarization class $h\in M$. We then define a partition $\Delta_h(\Lambda)= \Delta_\infty\sqcup \Delta_f$ by setting
$$\Delta_\infty=\{\delta\in\Delta(\Lambda)\mid \delta.h=0 \textrm{ and }  \delta.\delta_i\neq 0 \textrm{ for some }  i=1,\dots,5\}$$
and $\Delta_f=\{\delta\in\Delta(\Lambda)\mid \delta\perp M\}=\Delta(T)$ respectively. We denote by $H_\delta\subset \calD$ the hyperplane orthogonal to a root $\delta$, i.e. $H_\delta=\{\omega\in \bP(\Lambda\otimes \bC)\mid \omega.\omega=0,\ \omega.\overline{\omega}>0,\ \omega.\delta=0, \textrm{ and } \omega\perp M\}$. The sets roots $\Delta_{\infty}$ and $\Delta_f$ are stable under the monodromy group $\Gamma$, and they define two arithmetic arrangements of hyperplanes (in the sense of \cite{looijengacompact}): $\calH_\infty=\cup_{\delta\in \Delta_\infty} H_\delta$ and $\calH_f=\cup_{\delta\in \Delta_f}H_\delta$ respectively.
\end{notation}

With these preliminaries, we obtain the following result describing the structure of the regular and the simple singularity part in the deformation space of $N_{16}$. This is analogous to the situation for the triangle singularities (\cite[Thm. 6.4]{looijengatriangle} and \cite{novaacta}).

\begin{theorem}\label{thmstructure}
With  notations as above, we have the isomorphisms:
\begin{itemize}
\item[i)] $\calM_r\cong \left(\calD\setminus (\calH_\infty\cup \calH_f)\right)/\Gamma$
\item[ii)] $\calM_f\cong \left(\calD\setminus \calH_\infty\right)/\Gamma$
\end{itemize}
\end{theorem}
\begin{proof} This follows from theorem \ref{mainthm3} and the fact that the period point $\omega\in \calD/\Gamma$ associated to the $K3$ surface $\scl$ satisfies: 
\begin{itemize}
\item[i)] $C$ is singular if and only if $\omega\in \calH_f/\Gamma$;
\item[ii)] the intersection $C\cap L$ is not transversal if and only if $\omega\in \calH_\infty/\Gamma$ 
\end{itemize}
(cf. lemma \ref{lemmalternative}). \end{proof}

\subsection{The simple singularities locus}\label{secsimple}
One important application of the construction of section \ref{sectk3} is the description of the possible combinations of simple singularities for a nearby fiber in the deformation of $N_{16}$. This type of application was extensively considered in the case of unimodal singularities (e.g. for triangle singularities see  \cite[\S6]{looijengatriangle} and \cite{urabetriangle}). In our situation, we obtain the following purely arithmetic characterization of the nearby singularities. 

\begin{theorem}\label{simplecombinations}
A configuration of simple singularities occurs as the singular locus of a fiber in the universal deformation of $N_{16}$ if and only if the root lattice $R$ associated to the configuration satisfies the following property: there exists an overlattice $N$ of $M\oplus R$ such that
\begin{itemize}
\item[i)] the composition $M\hookrightarrow M\oplus R\subseteq N$ is primitive;
\item[ii)] the root sublattice of $\langle h\rangle^\perp_N$ coincides with  $R\oplus A_1^{\oplus 5}$;
\item[iii)] there exists a primitive embedding of $N$ into the $K3$ lattice $\Lambda$.
\end{itemize}
\end{theorem}
\begin{proof}
Assume that the given configuration is realized for a nearby fiber. By Pinkham's compactification procedure, this is equivalent to the existence of degree $5$ pair $(C,L)$ such that $C$ has the given combination of simple singularities, and $L$ is transversal to $C$.  It is immediate to see that $M\oplus R$ embeds in the $\Pic(\scl)$, and its saturation  $N:=\Sat_{\Pic(\scl)}(M\oplus R)$ satisfies the required properties.

Conversely, assume that there exists a lattice $N$ satisfying the above properties. By iii) and the surjectivity of the period map for the $K3$ surfaces there exists a $K3$ surface $S$ with $\Pic(S)=N$. The assumption i) says that $S$ is $M$-polarized. We can arrange that the $M$-polarization is normalized (without affecting the conditions i--iii). By the arguments of proposition \ref{surjectivity}, it follows that $h\in M$ defines a double cover $\pi:S\to \bP^2$ with branch locus a reducible sextic $B=C+L$. The singularities of the sextic $B$ are in bijective correspondence with the irreducible summands of root sublattices of  $\langle h\rangle^{^\perp}_{\Pic(S)}$. From ii), it then follows that $L$ is transversal to $C$ and the configuration of singularities of $C$ is described by $R$. 
\end{proof}

It is well known that the combinations of singularities that occur in the deformation of a given simple singularity are described by the proper subgraphs of the corresponding Dynkin diagram. From this fact and the previous theorem, it follows immediately that the singularities of a nearby fiber with only simple singularities can be deformed independently. 

\begin{corollary}\label{corsimple}
Assume that a given configuration of simple singularities occurs as the singular locus for a nearby fiber in the deformation of $N_{16}$. Then, any configuration obtained by deforming this configuration also occurs. \qed
\end{corollary}

A list of combinations with maximal Milnor number ($12$ or $13$) was obtained by Wall  \cite{wallquintic} by a case-by-case analysis of the degenerations of plane quintics. Such a list can be also obtained in an algorithmic way by applying \ref{simplecombinations} and the results of Nikulin \cite{nikulin}. Specifically, given a configuration of simple singularities, we consider the associated root lattice $R$ (the direct sum of the corresponding root lattices of type $A$--$D$--$E$). We can check that the configuration occurs  in three basic steps:
\begin{itemize}
\item[] Step I: {\it Find all the overlattices $N$ of $M\oplus R$}. These are classified by the isotropic subgroups of the corresponding discriminant group (\cite[Sect. 4]{nikulin}).   
\item[] Step II: {\it Given $N$ from Step I, check the properties i) and ii) of  \ref{simplecombinations}}. As the embedding $M\subset N$ constructed above is  explicit, this is routine. 
\item[] Step III: {\it Check that $N$ can be primitively embedded in $\Lambda$} (see \cite[Thm. 3.6.2]{nikulin}).
\end{itemize}

We illustrate theorem \ref{simplecombinations} and the above algorithm with two examples. 

\begin{example}
The worst $A_n$ singularity that occurs for a nearby fiber is $A_{12}$.  First, we note that $A_{13}$ does not occur. We apply \ref{simplecombinations} for $R=A_{13}$. Since $A_{M\oplus R}\cong (\bZ/2\bZ)^{ 4}\times \bZ/14\bZ$, it follows there is no proper overlattice of $M\oplus R$ satisfying i). Thus, we must have $N=M\oplus R$. The lattice $N$ cannot be embedded into $\Lambda$, since it does not satisfy the necessary condition $l(A_N)\le (\rk(\Lambda)-\rk(N))$ (here $l(A_N)=5$). On the other hand,  for $A_{12}$ the embedding of $N:=M\oplus R$ in $\Lambda$ exists by \cite[Thm. 3.6.2]{nikulin}. Note also that  all the singularities $A_n$ with $n\le 12$ occur (cf. \ref{corsimple}).
\end{example}

\begin{example}
We  note that typically in \ref{simplecombinations} we need to consider proper overlattice of $M\oplus R$. For example, the maximal number of nodes that can occur is $10$.  Let $R=A_1^{\oplus 10}$. Since $l(A_{M\oplus R})>\rk(\Lambda)-\rk(M\oplus R)$, the lattice $M\oplus R$ cannot be primitively embedded  into the $K3$ lattice $\Lambda$. Thus we have to consider a proper overlattice $N$ of $M\oplus R$.  The overlattices of $M\oplus R$ are classified by the isotropic subgroups $H$ (w.r.t. the induced quadratic form $q_M$) of $A_{M\oplus R} \cong (\bZ/2\bZ)^{\oplus 14}$. Since $A_N\cong H^{\perp}/H$, we obtain $l(H)\ge 4$. Together with the conditions i--ii, this  restriction gives (up to permutations) only one possibility for $H$. The resulting overlattice $N$  embeds into $\Lambda$. Geometrically, the divisibility conditions given by taking the overlattice $M\oplus R\subset N$ translate into collinearity conditions for the nodes. Thus, the need to consider the overlattice $N$ can be interpreted as saying that  $10$ nodes occur only for quintics consisting of $5$ generic lines. 
\end{example}

\subsection{Triangle singularities adjacent to $N_{16}$ and flips}\label{sectriangle}
In \S\ref{secdiscriminants} we have seen that the birational modifications that occur when passing from $\calM\cong\calD/\Gamma$ to $\calM_{\le 0}$ occur over the arrangement of hyperplanes $\calH_{\infty}$. Here we show that, in fact, the variation of GIT quotients  between the slopes $1$ and $\frac{5}{2}$ is essentially dictated by the arithmetic properties of the arrangement $\calH_{\infty}$. The situation is very similar to that considered by Looijenga \cite[\S10]{looijengacompact}, but due to the fact that we work in a relative version, we prove less than \cite[10.1]{looijengacompact}. As in \cite{looijengacompact}, the intersections of the various hyperplanes from the arrangement give a natural stratification of $\calH_{\infty}$. We show then that the variation of GIT quotients flips the strata of this stratification, starting with the highest codimension one. Geometrically, these birational transformations can be interpreted as removing the equisingular stratum $\Sigma^+_\calT\subset \calM_{\le 0}$ of  the triangle singularities of type $\calT$ and replacing it by the deformations of negative weight  $\Sigma^-_\calT\subset \calD/\Gamma$ of the singularities of type $\calT$ (N.B. everything is modulo a $\bC^*$-action). 

The situation described above is essentially identical to that for the deformation of triangle singularities considered in Looijenga \cite[\S10]{looijengacompact}. However, there is a major difference:  we have to distinguish two cases for the intersection of hyperplanes from the arrangement $\calH_{\infty}$, either ``transversal'' or ``tangential''. The birational transformations that occur in our situation affect only the strata coming from tangential intersections. In other words, we perform only some of the steps in Looijenga's compactification procedure. It is likely that an appropriate modification of Looijenga's construction will select precisely  these steps and give the natural fibration $\calM(\frac{5}{2}-\epsilon)\to \calM(\frac{5}{2})$ (see the following remark). 
 
 \begin{remark}\label{remrelation}
 We noted in the proof of \ref{lemmadiscriminant} the existence of two natural discriminant divisors  $\Sigma_0$ and $\Sigma_1$  parameterizing  the pairs $(C,L)$ with singular $C$, and the pairs with non-transversal $C\cap L$ respectively. Additionally,  it is easily seen that $\Sigma_0$ and $\Sigma_1$ are the natural polarizations on the spaces $\calM(0)$ and $\calM(\frac{5}{2})$ (or more precisely they are nef divisors on appropriate birational models that give the Mori fiber structures $\calM(\epsilon)\to \calM(0)$ and $\calM(\frac{5}{2}-\epsilon)\to \calM(\frac{5}{2})$ respectively). It was noted in Thaddeus \cite[(3.4)]{thaddeus} that a $1$-parameter  variation of GIT quotients is equivalent to running a directed minimal model program (MMP). In our situation, this means that we are running a MMP to change the polarization $\Sigma_0$ into $\Sigma_1$\footnote{Since all the birational morphisms occurring between $0$ and $\frac{5}{2}$ are small, the divisors $\Sigma_0$ and $\Sigma_1$ are well defined Weil divisors in $\calM(t)$. Moreover, for noncritical $t\in (1,\frac{5}{2})$,  the normal variety $\calM(t)$ is a geometric quotient, in particular $\bQ$-factorial. Thus, for such $t$, $\Sigma_0$ and $\Sigma_1$ are $\bQ$-Cartier.}. On the other hand, the proper transforms of $\Sigma_0$ and $\Sigma_1$ in $\calM(1)\cong (\calD/\Gamma)^*$ can be identified with the discriminants $\calH_f/\Gamma$ and $\calH_\infty/\Gamma$ respectively (cf. \ref{lemmalternative}). Since Looijenga's construction \cite{looijengacompact} is essentially the same as running a directed MMP program for the divisor $\calH_{\infty}/\Gamma$, we see that Looijenga's construction should coincide with the variation of GIT construction. The technical details to make this identification more precise are beyond the aim of this paper, so we choose to prove only some weaker statements (e.g. Thm. \ref{thmstrata} below). 
  \end{remark}
 
\subsubsection{The stratification of $\calH_{\infty}$ and flips} We make the basic observation  that we can characterize in arithmetic terms the locus of degree $5$ pairs which are stable at $t=1$, but cease to be stable before $t=\frac{5}{2}$. Namely, all such pairs have the property that there exists a point $p$ such that $\mult_p(C\cap L)\ge 3$ and the singularity at $p$ of $C+L$ is a simple singularity (cf. \ref{lemmasimplesing}). This type of conditions can be translated into statements about the Picard lattice of $\scl$ as in \S\ref{quinticperiod}. 

\begin{definition}
Let $\delta_1, \delta_2\in \Delta_\infty$ be two roots such that the corresponding hyperplanes in $\calH_{\infty}$ satisfy $\calH_{\delta_1}\neq \calH_{\delta_2}$ and $\calD\cap \calH_{\delta_1}\cap \calH_{\delta_2}\neq \emptyset$. We say that the intersection of the hyperplanes $\calH_{\delta_1}$ and $\calH_{\delta_2}$ is {\it tangential} iff $e_i.\delta_j\neq 0$ for $j=1,2$ and some fixed $i\in\{1,\dots,5\}$. Otherwise, we say that the intersection is {\it transversal}.
\end{definition}

Geometrically, a hyperplane $\calH_{\delta}$ parameterizes the pairs $(C,L)$ with $L$ tangent to $C$. The tangential (transversal) intersection of such hyperplanes corresponds to the case when $L$ is inflectional (resp. bitangent) to $C$. As mentioned above, the birational modification $\calM(1)\dashrightarrow \calM(\frac{5}{2}-\epsilon)$ occurs precisely over the codimension $2$ locus in $\calD/\Gamma$ given by the tangential intersections of hyperplanes from $\calH_{\infty}$. The following result gives a finer stratification of this locus in terms of the factorization of  the birational map $\calM(1)\dashrightarrow \calM(\frac{5}{2}-\epsilon)$ in a sequence of simple flips. 

\begin{theorem}\label{thmstrata}
For $t\in\{\frac{10}{7}, \frac{8}{5}, \frac{5}{3}, \frac{7}{4}, \frac{13}{7}, 2,\frac{11}{5}\}$, let $\Sigma_t$ be the closure in $\calM\cong \calD/\Gamma$ of the locus of degree $5$ pairs $(C,L)$ with the interval of stability $[\alpha, t]$ (for some $\alpha<1$). Then $\Sigma_t$ is an irreducible intersection of hyperplanes from $\calH_\infty$, which has the property that a generic point $\omega\in \Sigma_t$ determines an $M$-polarized $K3$ surface with Picard lattice isometric to $M_t$ as given in table \ref{tablelattice}. Additionally, an embedding $\Sigma_t\subset \Sigma_t'$ coresponds to a tower of primitive embeddings $M\hookrightarrow M_{t'}\hookrightarrow M_t\hookrightarrow \Lambda$ (where $M$ and $\Lambda$ are as in \ref{notationm0}). The possible inclusions are described by table \ref{tablelattice}.
\end{theorem}
\begin{proof}
By the results of section \ref{sectdeg5} (esp. \ref{lemmasimplesing} and table \ref{simplesing}), the pairs belonging to $\Sigma_t$ are characterized by the worst singularity at infinity, i.e. the type of the simple singularity occurring  for $C+L$  in the point of maximal multiplicity for the intersection $C\cap L$. The list of possibilities is given in table \ref{tablelattice}. By applying the procedure given by \ref{markinggen},  we obtain lattice embeddings $M\hookrightarrow M_t\subseteq \Pic(\scl)$ (an explicit example is computed in lemma \ref{trianglez11} below) such that $M\hookrightarrow \Pic(\scl)$ is a normalized embedding (cf. Def. \ref{goodmark}). By arguments  similar to those in \S\ref{quinticperiod}, one verifies that the embedding $M_t\subseteq \Pic(\scl)$ is primitive, that generically we have equality, and that all embeddings $M\hookrightarrow M_t\hookrightarrow \Lambda$ are conjugate by $\Gamma$ (N.B. $M\hookrightarrow \Lambda$ is assumed fixed). The fact that $\Sigma_t$ is an intersection of hyperplanes from $\calH_\infty$ is equivalent to saying that $M_t$ is generated by elements of $M$ and $\rank(M_t)-\rank(M)$ roots orthogonal to $h$.  Finally, the statement about the inclusion $\Sigma_t\subset \Sigma_t'$ is a  compatibility statement  for the construction of \ref{markinggen}. This is easily verified by using the obvious embeddings of root lattices corresponding to the singularities at infinity (e.g. $D_8\oplus A_1\hookrightarrow D_{10}$).
\end{proof}

\begin{remark}
Note that the embeddings $M\hookrightarrow M'\hookrightarrow \Lambda$ determine embeddings $T'\hookrightarrow T$ for the transcendental lattices (i.e. the orthogonal complements in $\Lambda$), which in turn determine embeddings of bounded symmetric domains $\calD'\subset \calD$. Thus each $\Sigma_t$ has the structure  of a locally symmetric variety $\calD'/\Gamma'$ for some $\calD'$ and $\Gamma'$ (or more precisely $\calD'/\Gamma'\to \Sigma_t\subset \calD/\Gamma$ is a normalization morphism).
\end{remark}

\begin{table}[htb]
\begin{center}
\renewcommand{\arraystretch}{1.25}
\begin{tabular}{|c|c|c|c|c|c|c|}
\hline
Case & Wall & Sing. at infinity& Picard lattice $M_t$ &Codim.&Specialization of\\ 
\hline\hline
(1)&$\frac{10}{7}$& $D_{10}$ &$T_{2,3,8}$&5& (2), (3)\\ \hline
(2)&$\frac{8}{5}$  & $D_8+A_1$&$T_{2,4,6}$& 4 &(5), (6)\\ \hline
(3)&$\frac{5}{3}$ & $A_9$&$T_{2,5,5}$ & 4 &(5)\\ \hline
(4)&$\frac{7}{4}$ & $E_7+2A_1$ &$E_8\oplus U$&  4 & (6)\\ \hline
(5)&$\frac{13}{7}$&$A_7+A_1$ &$T_{3,4,4}$  &  3&  (7)\\ \hline
(6)&$2$           &$D_6+2A_1$ &$E_7\oplus U$   &  3&(7)\\ \hline
(7)&$\frac{11}{5}$&$A_5+2A_1$ &$E_6\oplus U$ & 2  &- \\ \hline
\end{tabular}
\vspace{0.2cm}
\caption{The stratification of $\calH_{\infty}$}\label{tablelattice}
\end{center}
\end{table}

\begin{lemma}\label{trianglez11}
Let $S$ be a $K3$ surface such that $\Pic(S)\cong T_{2,3,8}$. Then there exists $h\in\Pic(S)$ a degree two base point free polarization which defines a double cover map $\pi:S\to \bP^2$ with branch locus a reducible sextic $B=C+L$ such that 
\begin{itemize}
\item[i)] $L$ is a line meeting the residual quintic $C$ in a single point $p$;
\item[ii)] $p$ is an ordinary double point for $C$.
\end{itemize}
Conversely, for a sufficiently general degree $5$ pair $(C,L)$ satisfying i) and ii), we have $\Pic(\scl)\cong T_{2,3,8}$.
\end{lemma}
\begin{proof}
Let $S$ be a $K3$ surface such that $\Pic(S)\cong T_{2,3,8}$. We claim that: {\it there exists a normalized embedding $j:M\hookrightarrow \Pic(S)$  such that  $\langle j(h)\rangle^\perp_{\Pic(S)}\cong D_{10}$}. Assuming the claim, from \S\ref{secsurjectivity}, it follows that $j(h)$ defines a base point free degree $2$ polarization. Again from \S\ref{secsurjectivity}, the resulting double cover $S\to \bP^2$ is branched along a reducible sextic  $B=C+L$. The condition $\langle j(h)\rangle^\perp_{\Pic(S)}\cong D_{10}$ implies that $B$ has a $D_{10}$ singularity at some point $p$. For numerical reasons, $j(l')$ has to meet the lattice $D_{10}$. It follows that  $p\in L$ (note lemma \ref{irredl}).  The only geometric possibility is the situation given by i) and ii). To prove the claim on the existence of a normalized embedding, it suffices to produce an embedding $M\hookrightarrow T_{2,3,8}$ such that the sublattice orthogonal to $h$ is $D_{10}$ (cf.  \ref{arrangegood}). 
This is done below (see (\ref{explicitformula})).

Conversely, assume that the pair $(C,L)$ satisfies the conditions i) and ii). Since $\mult_p(C\cap L)=5$ and $p$ is a node for $C$, it follows that $p$ is a singularity of type $D_{10}$ for $C+L$. The surface $\scl$ is obtained by desingularizing the double cover of $\bP^2$ branched along $C+L$. It follows that $\Pic(\scl)$ contains $10$ exceptional curves $e_1,\dots,e_{10}\in \Pic(S_{(C,L)})$, which together with the class $l'$ span the $T_{2,3,8}$ sublattice:  
$$\xymatrix@R=.25cm{
 {l'}              &{e_9}         &{e_8}                &{e_7}        &{e_6}    &{e_5}&{e_4}                &{e_3}        &{e_2}    &{e_1}\\
{\circ}\ar@{-}[r]&{\bullet}\ar@{-}[r]&{\bullet}\ar@{-}[r]\ar@{-}[dd]&{\bullet}\ar@{-}[r]&{\bullet}\ar@{-}[r]&{\bullet}
\ar@{-}[r]&{\bullet}\ar@{-}[r]&{\bullet}\ar@{-}[r]&{\bullet}\ar@{-}[r]&{\bullet}\\
&                                      &                   &                              &                   &  \\
&                                      &{\bullet}          &                              &                   &  \\
&                                      &{e_{10}}         &                              &    &  
}$$ 
The $M$-polarization $j:M\hookrightarrow \Pic(\scl)$ constructed by \ref{markinggen} factors through $T_{2,3,8}$:
\begin{eqnarray}
j(l')&=&l'\notag\\
j(e_1)&=&e_1+2(e_2+\dots+e_8)+e_9+e_{10}\notag\\
j(e_2)&=&e_3+2(e_4+\dots+e_8)+e_9+e_{10}\notag\\
j(e_3)&=&e_5+2(e_6+e_7+e_8)+e_9+e_{10}\label{explicitformula}\\
j(e_4)&=&e_7+2e_8+e_9+e_{10}\notag\\
j(e_5)&=&e_9\notag.
\end{eqnarray}
The conclusion now follows as in \S\ref{picardlattice}.
\end{proof}

\begin{remark}
In principle, it is possible to classify the possible intersection strata  of the hyperplanes from $\calH_\infty$ in purely  arithmetic terms. Appropriately considered (e.g. ignore the transversal intersections), the resulting stratification of $\calH_\infty$ should coincide with the stratification of theorem \ref{thmstrata}. For triangle singularities, the analogous analysis  was done by Looijenga \cite[\S5]{looijengatriangle}. The arithmetic of our situation is considerably more involved (e.g. the occurrence of high codimension strata). 
\end{remark}

\subsubsection{The geometric interpretation of the flips} We close by taking a closer look at the flips that occur at the walls $\frac{10}{7}$, $\frac{8}{5}$, $\frac{5}{3}$,  and $\frac{13}{7}$.   These are the walls where we introduce strata corresponding to the triangle singularities adjacent to $N_{16}$ into the universal deformation space $\calM_{\le 0}$. We recall  the relevant adjacency diagram:
$$\xymatrix@R=.25cm@C=.25cm{
{Z_{11}}&&\ar@{>}[ll]{Z_{12}}&&\\
&&&&\\
{W_{12}}\ar@{>}[uu]&&\ar@{>}[uu]\ar@{>}[ll]{W_{13}}&&\ar@{>}[ll]{N_{16}}
}$$ 

For simplicity, we focus on the  first wall $\frac{10}{7}$ and the associated singularity $Z_{11}$. The variation of GIT quotients gives a commutative diagram of birational maps:
$$\xymatrix@R=.25cm@C=.25cm{
&&{\calM(\frac{10}{7}-\epsilon)}\ar@{>}[ddll]_{\pi}\ar@{>}[ddrr]^{\pi_{-}}\ar@{<-->}[rrrr]&&&&{\calM(\frac{10}{7}+\epsilon)}\ar@{>}[ddll]_{\pi_{+}}\ar@{<-->}[rrrr]&&&&{\calM(\frac{5}{2}-\epsilon)\supset \calM_{\le0}}\\
&&&&\\
{\calM(1)}\ar@{<-->}[rrrr]&& &&{\calM(\frac{10}{7})}
}
$$  
The morphism $\pi$ is a small blow-up of the boundary $\calD/\Gamma$ in $(\calD/\Gamma)^*$, and leaves untouched  the stratum $\Sigma_{\frac{10}{7}}$ given by theorem \ref{thmstrata}.  The structure of the morphisms $\pi_{-}$ and $\pi_+$ is very simple (see  \S\ref{symquintic}). They contract  subspaces $E_{\pm}$ in $\calM(\frac{10}{7}\pm \epsilon)$ to a point $z_0$ in $\calM(\frac{10}{7})$ corresponding to the orbit of the pair $(C_0,L_0)$, where
$$C_0:x_0x_1x_2^3+x_1^5=0,\ \ L_0:x_0=0$$
(such a pair is semistable only at $t=\frac{10}{7}$).  The curve $C_0$ has a singularity of type $Z_{11}$ at $x_0=(1:0:0)$, and $L_0$ intersects $C$ with multiplicity $5$ in a unique point $p=(0:0:1)$.  Additionally, the stability condition for degree $5$ pairs identifies  the two exceptional loci $E_{\pm}$:
\begin{itemize}
\item[($-$)] $E_{-}=U_{-}/G$, where $U_{-}\subset X^s(\frac{10}{7}-\epsilon)$ is the $G$-invariant set of pairs $(C,L)$ such that there exists a point $p$ with $\mult_p(C\cap L)=5$ and $p$ is an ODP for $C$.  Furthermore, $C$  has at worst simple, simple elliptic or cusp singularities.
\item[($+$)] $E_{+}=U_+/G$ , where $U_+\subset X^s(\frac{10}{7}+\epsilon)$ is the $G$-invariant set of pairs $(C,L)$ such that $C$ has a singularity of type $Z_{11}$. Furthermore, the intersection $C\cap L$ is not too degenerate (better than $L$ is $5$-fold tangent to $C$ in a node).
\end{itemize}
In particular, it follows that the open stratum $E_+^\circ$ of $E_+$  parameterizing  pairs with $C\cap L$ transversal is precisely the $Z_{11}$ equisingular stratum in $\calM_{\le 0}$. Similarly, $E_-^\circ$  parameterizing  pairs with $C$ having simple singularities is  the stratum $\Sigma_{\frac{10}{7}}\subset \calD/\Gamma$.

It is interesting to note that Luna's slice theorem and the arguments of \S\ref{sectpinkham} identify $E_{-}$ to the quotient $S_{-}\gquot \Aut(C_0,L_0)$, where $S_{-}$ denotes the deformations of negative weight of the $Z_{11}$ singularity of $C_0$ at $p$. The deformations of negative weight of $Z_{11}$ were studied by Pinkham \cite{pinkhamduality} and Looijenga \cite{looijengatriangle}. In particular, they identify the simple singularities locus modulo the $\bC^*$-action with the moduli space of $T_{2,3,8}$-polarized $K3$ surfaces. By lemma \ref{trianglez11}, a $T_{2,3,8}$-polarized $K3$ surfaces is $M$-polarized, and the corresponding locus in $\calD/\Gamma$ is $\Sigma_{\frac{10}{7}}$. We conclude that the identification $\Sigma_{\frac{10}{7}}\cong E_-^\circ$ is precisely Pinkham's construction for $Z_{11}$ (N.B. the normalization of $\Sigma_{\frac{10}{7}}$ is isomorphic to $\calD'/\Gamma'$ for appropriate choices for $\calD'$ and $\Gamma'$). The role of the lattice $T_{2,3,8}$ for $Z_{11}$ is explained by Pinkham \cite{pinkhamduality}, and a posteriori  it offers an explanation for the occurrence of $T_{2,3,8}$ in table \ref{tablelattice}.

The situation for the other triangle singularities adjacent to $N_{16}$ is similar. One only has to note that an adjacency of triangle singularities determines a flip. For example, $Z_{12}\to Z_{11}$ gives the embedding $\Sigma_{\frac{10}{7}}\subset \Sigma_{\frac{8}{5}}$. To pass from $\Sigma_{\frac{8}{5}}\subset \calM(1)$ to $E_{-}\subset \calM(\frac{8}{5}-\epsilon)$, we first need to flip the locus $\Sigma_{\frac{10}{7}}$ (roughly speaking, this introduces the $Z_{11}$ stratum in the deformation of $Z_{12}$). This is easily seen to coincide with Looijenga's construction \cite[\S10]{looijengacompact} for the triangle singularities. The locus $\Sigma_{\frac{10}{7}}\subset \Sigma_{\frac{8}{5}}$ is precisely the locus of ``critical embeddings'' of Looijenga \cite[\S5, \S6]{looijengatriangle}. Equivalently, the restriction of arrangement $\calH_{\infty}$ to the subdomain $\calD''$ corresponding to the singularity $Z_{12}$ (or $\Sigma_{\frac{8}{5}}$) is the arrangement considered by Looijenga in \cite[\S10]{looijengacompact}.

\bibliography{references}
\end{document}